\pdfoutput=1

\documentclass[a4paper,10pt,reqno]{amsart}

\usepackage[utf8]{inputenc}
\usepackage[english]{babel}

\usepackage{amsmath}
\usepackage{amssymb}
\usepackage{BOONDOX-cal}
\usepackage{stmaryrd}
\usepackage{thmtools}

\usepackage{graphicx}
\usepackage{caption}
\usepackage[all,2cell]{xy}
\UseAllTwocells
\usepackage{tikz}
\usetikzlibrary{arrows.meta,decorations.pathmorphing}

\usepackage{hyperref}
\newcommand{\doi}[1]{\href{https://doi.org/#1}{\nolinkurl{https://doi.org/#1}}}

\makeatletter
\newcommand*\bigcdot{\mathpalette\bigcdot@{.5}}
\newcommand*\bigcdot@[2]{\mathbin{\vcenter{\hbox{\scalebox{#2}{$\m@th#1\bullet$}}}}}
\newcommand{\leqnomode}{\tagsleft@true}
\newcommand{\reqnomode}{\tagsleft@false}
\makeatother

\newcommand{\bb}[1]{\ensuremath{\lvert #1 \rvert}}
\newcommand{\vs}[1]{\mathbin{\downarrow}#1}

\newsavebox{\xymor}
\newcommand{\xynode}{\makebox[0ex]{}}
\savebox{\xymor}{\ensuremath{%
\xymatrix@1@C=19pt{\xynode \ar@{>}[r] & \xynode }}}
\newcommand{\mor}{\usebox{\xymor}}

\theoremstyle{plain}
\newtheorem{theorem}{Theorem}[section]
\newtheorem{proposition}[theorem]{Proposition}

\newtheorem{corollary}[theorem]{Corollary}
\newtheorem{lemma}[theorem]{Lemma}
\newtheorem{remark}[theorem]{Remark}
\theoremstyle{definition}
\newtheorem{definition}[theorem]{Definition}
\newtheorem{example}[theorem]{Example}

\theoremstyle{remark}

\numberwithin{equation}{section}

\begin{document}
\title[Riguet and Generalized Congruences on a Category]{Riguet and Generalized Congruences on a Category: Relationships and Applications}

\author[Climent]{J. Climent Vidal}
\address{Universitat de Val\`{e}ncia\\
         Departament de L\`{o}gica i Filosofia de la Ci\`{e}ncia
         }
\email{Juan.B.Climent@uv.es}
\author[Cosme]{E. Cosme Ll\'{o}pez}
\address{Universitat de Val\`{e}ncia\\
         Departament de Matem\`{a}tiques\\ and 
         Nantong University\\
         School of Mathematics and Statistics
         }
\email{Enric.Cosme@uv.es}
\author[Ruiz]{R. Ruiz Mora}
\address{Universitat de Val\`{e}ncia\\
         Departament de Matem\`{a}tiques
         }
\email{Raul.Ruiz-Mora@uv.es}


\begin{abstract}

We investigate Riguet congruences and generalized congruences on a category, focusing on their interrelations from both lattice-theoretic and category-theoretic perspectives.
We also characterize functors that are full and surjective on objects in terms of regular epimorphisms, extremal epimorphisms and in terms of strong and regular generalized congruences.
On the lattice-theoretic side, we prove that for a category $\mathsf{C}$, the set $\mathrm{RCgr}(\mathsf{C})$ of all Riguet congruences, ordered by inclusion, is a bounded directed-complete ordered set, while the set $\mathrm{GCgr}(\mathsf{C})$ of all generalized congruences is an algebraic lattice. 
We establish a bridge between these structures via a Scott continuous morphism.
From a category-theoretic standpoint, we lift these results to relative adjunctions  between the categories $\mathsf{RCgr}(\mathsf{C})$ and $\mathsf{GCgr}(\mathsf{C})$ associated to the above ordered sets, as well as between the categories $\mathsf{RCCat}$, of Riguet classified categories, and $\mathsf{GCCat}$, of generalized classified categories.
Furthermore, within Manes' framework of categories of $\mathsf{K}$-objects with structure, we investigate the relationship between the wide subcategory $\mathsf{RCCat}_{\mathrm{full}}$ of $\mathsf{RCCat}$, whose morphisms are the full morphisms of $\mathsf{RCCat}$, and $\mathsf{GCCat}$, relating these constructions to the Grothendieck theory of fibrations. 
Finally, we present applications of Riguet congruences across various mathematical  fields.

\end{abstract}

\maketitle

\section{Introduction}\label{S:intro}

The origin of this work lies in the field of many-sorted universal algebra. In the single-sorted case, the arity of an operation is an element of $\mathbb{N}$. The set $\mathbb{N}$ itself may be viewed, among other interpretations, as the set of finite cardinal numbers (finite initial ordinal numbers); as the underlying set of the free monoid on $1 = \{0\}$; and as the underlying set of the free abelian monoid on $1$ (for both constructions, see~\cite{nB70A})---note that the free monoid on $1$ and the free abelian monoid on $1$ are isomorphic. These three perspectives give rise to corresponding notions of arity for many-sorted algebras with a set of sorts $A$: as the set of finite 
$A$-sorted cardinal numbers, i.e., families $(n_{a})_{a\in A}$, where, for every $a\in A$, $n_{a}$ is a natural number and $\mathrm{card}(\{a\in A\mid n_{a}\neq 0\})<\aleph_{0}$; as the underlying set of the free monoid on $A$; and as the underlying set of the free abelian monoid on $A$---note that if $\mathrm{card}(A)\geq 2$, the free monoid on $A$ and the free abelian monoid on $A$ are not isomorphic. In much of the literature on many-sorted algebras, the second option is adopted, in one form or another. The reason is that, when writing the arguments of a many-sorted operation, one must impose an order on them. Accordingly, an element $\mathbf{a}$ of the free monoid on $A$ specifies that the $i$-th argument is of the sort given by  the $i$-th symbol of $\mathbf{a}$. On the other hand, as it is well known, the \emph{skeletal} category $\mathsf{Card}_{\mathrm{f}}$, of finite cardinal numbers, is equivalent to---and a full subcategory of---the category $\mathsf{Set}_{\mathrm{f}}$, of finite sets in $\mathbf{U}$ (see the Conventions below). Thus, for a set of sorts $A$, a natural question arises: for the category $\mathsf{Set}^{A}_{\mathrm{f}}$, of finite $A$-sorted sets (defined in Subsection~\ref{S:freemonoid}), which categories---whether skeletal or non-skeletal---related in a natural way to $A$ and not necessarily subcategories of $\mathsf{Set}^{A}_{\mathrm{f}}$, play a role analogous to that of $\mathsf{Card}_{\mathrm{f}}$ for $\mathsf{Set}_{\mathrm{f}}$, beyond its skeleton $\mathsf{Card}^{A}_{\mathrm{f}}$ of finite $A$-sorted cardinal numbers? 

The answer to this question, for skeletal categories, was obtained using a type of congruences introduced by Riguet in~\cite{JR60}. In~\cite{JR60}---a work essentially devoted to categorizing and generalizing Bourbaki's notion of species of structure, with or without morphisms, as introduced in~\cite{nB70}, and which, to the best of our knowledge, is not indexed in either MathSciNet or Zentralblatt MATH---Riguet introduces, prior to defining his own notion of congruence on a category (which we shall call, in his honour, a \emph{Riguet congruence}), the concept of congruence on a category as considered by Cartan and Eilenberg in~\cite{CE58}. This latter notion, which has become standard in category theory and is usually referred to simply as a \emph{congruence} on a category in the literature, is designated by Riguet as a \emph{relation de congruence forte} [strong congruence relation]; a brief treatment of this type of congruence may be found in~\cite{sM98}. By contrast, Riguet refers to his own notion simply as a \emph{relation de congruence} [congruence relation]. 
Riguet's congruences classify not only morphisms, as is customary, but also objects, with both classifications required to satisfy specific compatibility conditions. 

To delve, for now, a little deeper into their meaning: for a category $\mathsf{C}$, a Riguet congruence $\Phi$ on $\mathsf{C}$ consists of an equivalence relation on the objects of $\mathsf{C}$, denoted by $\Phi^{\mathrm{ob}}$, together with a family of relations $\Phi^{\mathrm{ob}} = \left(\Phi^{\mathrm{fl}}_{\scalebox{0.7}{$\left(\begin{smallmatrix}a&b\\a'& b'\end{smallmatrix}\right)$}}\right)_{\scalebox{0.7}{$\left(\begin{smallmatrix}a&b\\a'& b'\end{smallmatrix}\right)$} \in \Phi^{\mathrm{ob}}\times \Phi^{\mathrm{ob}}}$ in which $\scalebox{0.7}{$\left(\begin{smallmatrix}a&b\\a'& b'\end{smallmatrix}\right)$} \in \Phi^{\mathrm{ob}}\times \Phi^{\mathrm{ob}}$ means that $(a,a')$ and $(b,b') \in \Phi^{\mathrm{ob}}$ and where each $\Phi^{\mathrm{fl}}_{\scalebox{0.7}{$\left(\begin{smallmatrix}a&b\\a'& b'\end{smallmatrix}\right)$}}$ is a subset of $\mathrm{Hom}_{\mathsf{C}}(a,b) \times \mathrm{Hom}_{\mathsf{C}}(a',b')$, satisfying the compatibility conditions stated in Definition~\ref{DRiguet} (note that the rows of the matrix refer to the hom-sets, while the columns refer to $\Phi^{\mathrm{ob}}$-equivalent objects). Among these, the quadrangular completion condition---missing from Riguet's definition in~\cite{JR60}---which asserts that for every $(a,a'), (b,b') \in \Phi^{\mathrm{ob}}$ and every $f \in \mathrm{Hom}_{\mathsf{C}}(a,b)$, there exists an $f' \in \mathrm{Hom}_{\mathsf{C}}(a',b')$ such that $(f,f') \in \Phi^{\mathrm{fl}}_{\scalebox{0.7}{$\left(\begin{smallmatrix}a&b\\a'& b'\end{smallmatrix}\right)$}}$, plays a pivotal role, among others, in defining the quotient category $\mathsf{C}/{\Phi}$ and characterizing, via the notion of a Riguet-full functor, the regular epimorphisms between categories. In connection with the latter, let us point out that in~\cite{BBP99}, which does not cite~\cite{JR60}, Bednarczyk \textit{et al.} introduced the notion of generalized congruence on the categories in $\mathsf{cat}$, the category of $\mathbf{U}$-small categories (see the Conventions below), to characterize, among other things, the class of all regular epimorphisms in $\mathsf{cat}$ as those extremal epimorphisms which induce regular congruences---a particular type of generalized congruence---; Riguet congruences allow us, through the quadrangular completion condition, to define a subclass of functors, which we call Riguet-full functors, to characterize the class of all regular epimorphisms as those extremal epimorphisms which are Riguet-full. On the other hand, and curiously enough, some of the defining conditions of a Riguet congruence are analogous to some of the defining conditions of a double category; and, in fact, as we will see in this paper, to a particular Riguet congruence---on a category associated with the free monoid on a set---we will assign a double category. 

The central objective of this work is to elucidate the structural interplay between Riguet congruences and generalized congruences. The latter may be viewed as a natural extension in which equivalence is imposed not only on individual morphisms, but also on nonempty finite sequences of morphisms.

We begin by revisiting and refining the notion of Riguet congruence on a category, defining its associated quotient category and stating its universal property. We then characterize these congruences via Riguet-full functors.
Furthermore, we characterize functors which are full and surjective on objects both as Riguet-full extremal epimorphisms and as Riguet-full regular epimorphisms.

On this basis, we clarify the relationship between the two types of congruences from both lattice-theoretic and category-theoretic perspectives. Regarding the former, we prove that for a category $\mathsf{C}$, the poset $(\mathrm{RCgr}(\mathsf{C}),\subseteq)$ is a bounded directed-complete ordered set, while $(\mathrm{GCgr}(\mathsf{C}),\subseteq)$ is an algebraic lattice. We establish a bridge by constructing a Scott continuous morphism, mapping Riguet congruences to (strong) generalized congruences.

From a category-theoretic standpoint, we lift these order-theoretic results to the level of functors and relative adjunctions. Specifically, we prove the existence of a relative adjunction between the category $\mathsf{RCgr}(\mathsf{C})$, determined by $(\mathrm{RCgr}(\mathsf{C}),\subseteq)$, and the category $\mathsf{GCgr}(\mathsf{C})$, determined by $(\mathrm{GCgr}(\mathsf{C}),\subseteq)$. We also establish a relative adjunction between the category $\mathsf{RCCat}$, of Riguet classified categories, and $\mathsf{GCCat}$, of generalized classified categories.

Moreover, within Manes' framework of categories of $\mathsf{K}$-objects with structure~\cite{em76}, we investigate the wide subcategory $\mathsf{RCCat}_{\mathrm{full}}$ of $\mathsf{RCCat}$, whose morphisms are the full morphisms of $\mathsf{RCCat}$, together with the category $\mathsf{GCCat}$. For $\mathsf{RCCat}_{\mathrm{full}}$, we take $\mathsf{K} = \mathsf{Cat}_{\mathrm{full}}$, the wide subcategory of $\mathsf{Cat}$ whose morphisms are the full functors, while for $\mathsf{GCCat}$ we take $\mathsf{K} = \mathsf{Cat}$. In both cases, these categories are further related to the Grothendieck theory of fibrations~\cite{Gro71}.

Finally, we present a series of examples in different areas of mathematics where the notion of Riguet congruence can be naturally applied. Specifically: 
(1) on the category $\mathsf{C}(\mathbf{A}^{\star})$ associated to $\mathbf{A}^{\star}$, the free monoid on a set $A$, we define a Riguet congruence $\Phi$ such that the resulting (skeletal) quotient $\mathsf{C}(\mathbf{A}^{\star})/\Phi$ is isomorphic to $\mathsf{Card}^{A}_{\mathrm{f}}$; 
(2) on the fundamental groupoid $\pi_{1}(\mathbf{X}, x_{0})$ of a path-connected topological space $\mathbf{X}$ we define a Riguet congruence $\Phi$ such that the quotient category $\mathsf{\Pi}_{1}(\mathbf{X})/\Phi$ is a category with a single object whose endomorphisms are isomorphic to the fundamental group $\pi_{1}(\mathbf{X}, x_{0})$; 
(3) given a Grothendieck category $\mathsf{A}$ and a localizing subcategory $\mathsf{S}$, we construct a Riguet congruence $\Phi$ on $\mathsf{A}$ that characterizes the Serre localization $\mathsf{A}/ \mathsf{S}$; 
(4) given a deterministic finite automaton $\mathcal{A}$, we define a Riguet congruence $\Phi$ on $\mathsf{C}(\mathcal{A})$, the category associated to $\mathcal{A}$, and prove that the quotient is isomorphic to the action groupoid  
$\mathsf{C}(\mathcal{A})/\Phi \cong \mathsf{Act}(\mathbf{T}(\mathcal{A}/{\sim}), Q/{\sim})$; and 
(5) on the category $\mathsf{Rel}$, of sets and relations, we define a Riguet congruence $\Phi$ such that the quotient category $\mathsf{Rel}/\Phi$ is isomorphic to the category $\mathsf{Card}_{\mathrm{orb}}$, whose objects are the cardinal numbers $\kappa \in \mathbf{U}$, and for any cardinals $\kappa, \lambda \in \mathbf{U}$, the morphisms from $\kappa$ to $\lambda$ are the orbits of the set $\mathrm{Sub}(\kappa \times \lambda)$ under the action of the product of symmetric groups $\boldsymbol{\mathfrak{S}}_{\kappa} \times \boldsymbol{\mathfrak{S}}_{\lambda}$.

Before outlining the contents of the paper, we note that the results presented here were obtained independently of most of the existing literature dealing, in one way or another, with quotient categories. In fact, our study primarily drew on the work of Riguet~\cite{JR60}. It was only after carrying out our investigation into Riguet congruences and establishing some examples thereof that we became aware of the  work by Bednarczyk \textit{et al.}~\cite{BBP99} and decided to investigate the relationships between them. For completeness, and to place this work in a broader historical context, we briefly mention some of the most relevant contributions concerning various types of congruences on categories, without attempting to provide an exhaustive account. 
Among them, one may mention the works of Grothendieck~\cite{Gro57} and Gabriel~\cite{pG62}, who studied quotients of abelian categories by certain full subcategories; Mersch~\cite{JM65}, who investigated the problem of quotients of categories and proved, in particular, that the usual conditions for a congruence relation in algebra are not sufficient in categories; Ehresmann~\cite{Ehr65}, who defined, for single-sorted categories, the notions of quotient category by a compatible equivalence relation, by a subcategory, and by an ideal, as well as that of strict quotient category; Isbell~\cite{JI57}, who considered the notion of an identification category (derived from a category with an equivalence relation); Isbell~\cite{JI68}, who constructed the smallest congruence containing a given precongruence; Higgins~\cite{pjH71}, who defined quotients of groupoids by normal subgroupoids, thereby classifying both objects and morphisms; B\"{o}rger~\cite{RB77}, who discussed various types of congruences on categories, including those arising from subcategories of the square of a category satisfying the usual axioms of an equivalence relation on both objects and morphisms, and provided examples distinguishing them; Fiore \textit{et al.}~\cite{FPP08}, who defined two kinds of congruences on categories---the standard one and the one previously singled out by B\"{o}rger~\cite{RB77}---and two kinds of congruences on double categories; and, more recently, Borceux \textit{et al.}~\cite{BCGT23}, who built on~\cite{RB77,JI68} to establish some atypical properties of certain coequalizers (quotient functors) in $\mathsf{Cat}$. 

In what follows, we present a comprehensive exposition of the main results developed in the subsequent sections.

In Section~\ref{S:riguet}, we revisit and refine the notion of a Riguet congruence on a category $\mathsf{C}$. We then define the associated quotient category and establish its universal property. 
Furthermore, we characterize functors which are full and surjective on objects as those that are Riguet-full and extremal epimorphisms, as well as those that are Riguet-full and regular epimorphisms.

Thereafter, we clarify---both from the lattice-theoretic and the category-theoretic viewpoints---the relationship between Riguet congruences and generalized congruences on a category. More concretely, from the lattice-theoretic side, given a category $\mathsf{C}$, we state that $(\mathrm{RCgr}(\mathsf{C}),\subseteq)$, the set of all Riguet congruences on $\mathsf{C}$, ordered by inclusion, is a bounded directed-complete ordered set and that $(\mathrm{GCgr}(\mathsf{C}),\subseteq)$, the set of all generalized congruences on $\mathsf{C}$, ordered by inclusion, is an algebraic lattice. Then, after introducing the notion of strong generalized congruence on a category,
we characterize functors that are full and surjective on objects in terms of extremal epimorphisms and generalized kernels that are both strong and regular generalized congruences.
We also prove, and this constitutes a fundamental result, that there exists a Scott continuous morphism $(\bigcdot)^{\natural}$ from 
$(\mathrm{RCgr}(\mathsf{C}),\subseteq)$ to $(\mathrm{GCgr}(\mathsf{C}),\subseteq)$ whose image is included in $\mathrm{SGCgr}(\mathsf{C})$, the set of strong generalized congruences on $\mathsf{C}$.  Moreover, we provide a constructive description of $\Phi^{\natural\mathrm{fl}}$, the morphism part of $\Phi^{\natural}$, the strong generalized congruence on $\mathsf{C}$ associated to $\Phi$.

From the category-theoretic side, we show that, for any category $\mathsf{C}$, there exist an adjunction from $\mathsf{RCgr}(\mathsf{C})$, the category canonically associated to $(\mathrm{RCgr}(\mathsf{C}),\subseteq)$, to $\mathsf{SGCgr}(\mathsf{C})$, that arising from $(\mathrm{SGCgr}(\mathsf{C}),\subseteq)$, the ordered set of strong generalized congruences on $\mathsf{C}$. In addition, after defining the categories $\mathsf{RCCat}$, of Riguet classified categories, and $\mathsf{GCCat}$, of generalized classified categories, we prove that there exists a functor $(\bigcdot)^{\natural}$ from $\mathsf{RCCat}$ to $\mathsf{GCCat}$---derived from the above Scott continuous morphism. This functor sends a Riguet classified category $(\mathsf{C},\Phi)$, where $\Phi$ is a Riguet congruence on $\mathsf{C}$, to the generalized classified category $(\mathsf{C},\Phi^{\natural})$, where $\Phi^{\natural}$ is the strong generalized congruence on $\mathsf{C}$ associated to $\Phi$. We further prove that there are functors $Q_{\mathrm{R}}$, $U_{\mathrm{R}}$ from $\mathsf{RCCat}$ to $\mathsf{Cat}$, the category of all $\mathbf{U}$-categories (see the Conventions below), and $Q_{\mathrm{G}}$, $U_{\mathrm{G}}$ from $\mathsf{GCCat}$ to $\mathsf{Cat}$, together with natural transformations $\mathrm{pr}_{\mathrm{R}}$ from $U_{\mathrm{R}}$ to $Q_{\mathrm{R}}$, $\mathrm{pr}_{\mathrm{G}}$ from $U_{\mathrm{G}}$ to $Q_{\mathrm{G}}$ and $\alpha$ from $Q_{\mathrm{R}}$ to $Q_{\mathrm{G}}\circ (\bigcdot)^{\natural}$. Furthermore, we prove that the corestriction of the functor $(\bigcdot)^{\natural}$ to $\mathsf{SGCCat}$, the full subcategory of $\mathsf{GCCat}$ determined by the strong generalized classified categories, has a right adjoint $(\bigcdot)^{\flat}$ and that there are functors $Q_{\mathrm{S}}$, $U_{\mathrm{S}}$ from $\mathsf{SGCCat}$ to $\mathsf{Cat}$ together with natural transformations $\mathrm{pr}_{\mathrm{S}}$ from $U_{\mathrm{S}}$ to $Q_{\mathrm{S}}$ and $\beta$ from $Q_{\mathrm{R}}\circ (\bigcdot)^{\flat}$ to $Q_{\mathrm{S}}$. 
Following this, we prove that, for a category $\mathsf{C}$ and a Riguet congruence $\Phi$ on it, if there exists a choice function for the family $(\mathrm{Hom}_{\mathsf{C}}(a,b))_{(a,b) \in \Phi^{\mathrm{ob}}}$ satisfying a natural compatibility condition with respect to $\Phi$, then the functor $\alpha_{(\mathsf{C},\Phi)}$ from $\mathsf{C}/\Phi$ to $\mathsf{C}/\Phi^{\natural}$ is an equivalence. 

Following this, within Manes' framework of categories of $\mathsf{K}$-objects with structure, we investigate the wide subcategory $\mathsf{RCCat}_{\mathrm{full}}$ of $\mathsf{RCCat}$, whose morphisms are the full morphisms of $\mathsf{RCCat}$, and the category $\mathsf{GCCat}$. For $\mathsf{RCCat}_{\mathrm{full}}$ we take $\mathsf{K} = \mathsf{Cat}_{\mathrm{full}}$, the wide subcategory of $\mathsf{Cat}$ whose morphisms are the full functors, while for $\mathsf{GCCat}$, we take $\mathsf{K} = \mathsf{Cat}$. In both cases, these categories are further related to the Grothendieck theory of fibrations, by considering the restriction of $U_{\mathrm{R}}$ to $\mathsf{RCCat}_{\mathrm{full}}$, the functor $U_{\mathrm{G}}$, the restriction of $(\bigcdot)^{\natural}$ to $\mathsf{RCCat}_{\mathrm{full}}$ and the canonical embedding of $\mathsf{Cat}_{\mathrm{full}}$ into $\mathsf{Cat}$.

Section~\ref{S:examples} is devoted to presenting several illustrative applications of Riguet congruences across a diverse range of mathematical contexts. Specifically, in Subsection~\ref{S:freemonoid},  we consider the category $\mathsf{C}(\mathbf{A}^{\star})$ associated to the free monoid on a set $A$, and define a Riguet congruence $\Phi$ on $\mathsf{C}(\mathbf{A}^{\star})$ that identifies objects with identical letter frequencies and morphisms that commute with the canonical permutations associated to the pairs of identified words. We prove that the resulting skeletal quotient $\mathsf{C}(\mathbf{A}^{\star})/\Phi$ is isomorphic to $\mathsf{Card}^{A}_{\mathrm{f}}$, the skeleton of the category of finite $A$-sorted sets. This Riguet congruence replaces the free monoid with its abelianization.

In Subsection~\ref{S:groupoid}, we consider the fundamental groupoid $\mathsf{\Pi}_{1}(\mathbf{X})$ of a path-connected space $\mathbf{X}$ and define a Riguet congruence collapsing the object set via transport isomorphisms. This yields a single-object quotient category whose endomorphisms are isomorphic to the fundamental group $\pi_{1}(\mathbf{X}, x_{0})$, thereby eliminating dependence on the base point and isolating the global algebraic structure. The construction relies on path-connectedness: since all objects are isomorphic, local data can be canonically transported to a fixed base point, reducing the groupoid to its algebraic core.

In Subsection~\ref{S:grothendieck}, we consider a Grothendieck category $\mathsf{A}$ with a localizing subcategory $\mathsf{S}$ and define a Riguet congruence $\Phi$ on $\mathsf{A}$ that characterizes the Serre localization $\mathsf{A}/\mathsf{S}$. We prove that the resulting skeletal quotient $\mathsf{A}/\Phi$ is isomorphic to the standard skeleton of $\mathsf{A}/\mathsf{S}$, effectively internalizing the quotient structure within $\mathsf{A}$. This Riguet congruence is enabled by the existence of the section functor $R$ (the right adjoint to the localization functor $Q$), whose fully faithful nature ensures that $Q$ is full. In this setting, the Riguet congruence operates by lifting every ``abstract fraction'' of the localization back to a concrete morphism in $\mathsf{A}$. The fullness of $Q$ guarantees that the quadrangular completion condition is satisfied, allowing morphisms between $\mathsf{S}$-closed objects to provide a complete internal description of the quotient. By identifying objects via their $\mathsf{S}$-injective hulls and morphisms through the unit of the adjunction $\eta \colon \mathrm{Id}_{\mathsf{A}} \mor RQ$, we reduce the localized structure to its canonical skeletal representation within the category of $\mathsf{S}$-closed objects.

In Subsection~\ref{S:automata}, we consider the free category $\mathsf{C}(\mathcal{A})$ associated to a deterministic finite automaton $\mathcal{A} = (Q, \Sigma, \delta, q_{0}, F)$, and define a Riguet congruence $\Phi$ on $\mathsf{C}(\mathcal{A})$ based on the Myhill-Nerode congruence $\sim$ on $\mathcal{A}$ and the kernel $\mathbf{K}$ of the canonical  homomorphism from $\boldsymbol{\Sigma}^{\star}$, the free monoid on $\Sigma$, to $\mathbf{End}(Q/{\sim})$, the monoid of self-mappings of the quotient of $Q$ by $\sim$, whose image is $\mathbf{T}(\mathcal{A}/{\sim})$,  the transition monoid of $\mathcal{A}/{\sim}$. We prove that the resulting quotient is isomorphic to $\mathsf{Act}(\mathbf{T}(\mathcal{A}/{\sim}), Q/{\sim})$, the action category.

In Subsection~\ref{S:rel}, we consider the category $\mathsf{Rel}$ of sets and binary relations, and define a Riguet congruence on it that identifies objects by their cardinality and morphisms by their orbits under the action of symmetric groups. This yields an isomorphism with $\mathsf{Card}_{\mathrm{orb}}$, reducing relations to their structural essence independently of labeling. The Riguet congruence acts via symmetric group actions: passing to orbits removes element-level identities and recovers the underlying combinatorial skeleton.

\subsection{Conventions}\label{S:conventions}
In this work, we fix once and for all a Grothendieck universe $\mathbf{U}$ (see \cite{sM98}, pp.~21--24). A set is called $\mathbf{U}$-small if it is an element of $\mathbf{U}$, and $\mathbf{U}$-large if it is a subset of $\mathbf{U}$. Throughout, we use standard concepts and constructions from set theory and category theory. However, concerning set theory, we adopt the following conventions. We denote by $\mathbb{N}$ the set of all natural numbers and, for every $n \in \mathbb{N}$, we set $n = \{0, \dots, n-1\}$. For a set $A$, we denote by $\mathrm{Sub}(A)$ the set of all subsets of $A$, by $\Delta_{A}$ the diagonal of $A$, and, when convenient, by $\nabla_{A}$ the relation $A \times A$. If $f$ is a mapping from $A$ to $B$ and $X \subseteq A$, we let $f[X]$ denote the direct image of $X$ under $f$. Moreover, if $Y \subseteq B$, we let $f^{-1}[Y]$ denote the inverse image of $Y$ under $f$. 
With regard to category theory, we assume, as in~\cite{sM98}, that the hom-sets of the categories are pairwise disjoint. Moreover, we let $\mathsf{Set}$ denote the category of sets in $\mathbf{U}$ and mappings between them, and $\mathsf{Set}_{\mathrm{f}}$ the full subcategory of $\mathsf{Set}$ consisting of the finite sets in $\mathbf{U}$. We denote by $\mathsf{cat}$ the category of all $\mathbf{U}$-small categories, i.e., the category whose objects are all categories $\mathsf{C}$ such that both $\mathrm{Ob}(\mathsf{C})$ and $\mathrm{Mor}(\mathsf{C})$ are elements of $\mathbf{U}$, and whose morphisms are the functors between $\mathbf{U}$-small categories. Finally, we denote by $\mathsf{Cat}$ the category of all $\mathbf{U}$-categories, i.e., the category whose objects are all categories $\mathsf{C}$ such that $\mathrm{Ob}(\mathsf{C}) \subseteq \mathbf{U}$ and, for every $x,y \in \mathrm{Ob}(\mathsf{C})$, $\mathrm{Hom}(x,y) \in \mathbf{U}$, and whose morphisms are the functors between $\mathbf{U}$-categories.

\section{Riguet congruences and their relationship with generalized congruences}\label{S:riguet}

We begin by revisiting and refining the notion of a Riguet congruence on a category $\mathsf{C}$, as in~\cite{JR60}. We prove that the set $\mathrm{RCgr}(\mathsf{C})$ of all Riguet congruences on $\mathsf{C}$, ordered by inclusion, is a bounded directed-complete ordered set. Furthermore, we define the associated quotient category and establish its universal property.

\begin{definition}\label{DRiguet}
A \emph{Riguet congruence} on a category $\mathsf{C}$ is an ordered pair
$\Phi=\left(\Phi^{\mathrm{ob}},\Phi^{\mathrm{fl}}\right)$ in which 
\begin{enumerate}
\item $\Phi^{\mathrm{ob}}$ is an equivalence relation on $\mathrm{Ob}(\mathsf{C})$ and 
\item $\Phi^{\mathrm{fl}} = \left(\Phi^{\mathrm{fl}}_{\scalebox{0.7}{$\left(\begin{smallmatrix}a&b\\a'& b'\end{smallmatrix}\right)$}}\right)_{\scalebox{0.7}{$\left(\begin{smallmatrix}a&b\\a'& b'\end{smallmatrix}\right)$} \in \Phi^{\mathrm{ob}}\times \Phi^{\mathrm{ob}}}$ (where ``$\mathrm{fl}$'' stands for the French word \emph{fl\`{e}che} and $\scalebox{0.7}{$\left(\begin{smallmatrix}a&b\\a'& b'\end{smallmatrix}\right)$} \in \Phi^{\mathrm{ob}}\times \Phi^{\mathrm{ob}}$ means that $(a,a')$ and $(b,b') \in \Phi^{\mathrm{ob}}$) a choice function for $(\mathrm{Sub}(\mathrm{Hom}_{\mathsf{C}}(a,b) \times \mathrm{Hom}_{\mathsf{C}}(a',b')))_{\scalebox{0.7}{$\left(\begin{smallmatrix}a&b\\a'& b'\end{smallmatrix}\right)$} \in \Phi^{\mathrm{ob}}\times \Phi^{\mathrm{ob}}}$ 
such that the following conditions are satisfied 
\begin{enumerate}
\item for every $a$, $b\in \mathrm{Ob}(\mathsf{C})$, if $(a,b) \in \Phi^{\mathrm{ob}}$, then $(\mathrm{id}_{a},\mathrm{id}_{b})\in \Phi^{\mathrm{fl}}_{\scalebox{0.7}{$\left(\begin{smallmatrix}a&a\\b&b\end{smallmatrix}\right)$}}$;
\item for every $a$, $b\in \mathrm{Ob}(\mathsf{C})$, $\Delta_{\mathrm{Hom}_{\mathsf{C}}(a,b)} \subseteq \Phi^{\mathrm{fl}}_{\scalebox{0.7}{$\left(\begin{smallmatrix}a&b\\a&b\end{smallmatrix}\right)$}}$, in other terms, if $f\in \mathrm{Hom}_{\mathsf{C}}(a,b)$, then $(f,f)\in \Phi^{\mathrm{fl}}_{\scalebox{0.7}{$\left(\begin{smallmatrix}a&b\\a&b\end{smallmatrix}\right)$}}$;
\item for every $(a,a'), (b,b') \in \Phi^{\mathrm{ob}}$, $\left(\Phi^{\mathrm{fl}}_{\scalebox{0.7}{$\left(\begin{smallmatrix}a&b\\a'& b'\end{smallmatrix}\right)$}}\right)^{-1}= \Phi^{\mathrm{fl}}_{\scalebox{0.7}{$\left(\begin{smallmatrix}a'&b'\\a&b\end{smallmatrix}\right)$}}$, in other terms, if $(f,f')\in \Phi^{\mathrm{fl}}_{\scalebox{0.7}{$\left(\begin{smallmatrix}a&b\\a'& b'\end{smallmatrix}\right)$}}$, then  $(f',f)\in \Phi^{\mathrm{fl}}_{\scalebox{0.7}{$\left(\begin{smallmatrix}a'&b'\\a& b\end{smallmatrix}\right)$}}$;
\item for every $(a,a')$, $(a',a'')$, $(b,b')$, $(b',b'') \in \Phi^{\mathrm{ob}}$, 
$
\Phi^{\mathrm{fl}}_{\scalebox{0.7}{$\left(\begin{smallmatrix}a'&b'\\a''&b''\end{smallmatrix}\right)$}} \circ \Phi^{\mathrm{fl}}_{\scalebox{0.7}{$\left(\begin{smallmatrix}a&b\\a'& b'\end{smallmatrix}\right)$}} \subseteq \Phi^{\mathrm{fl}}_{\scalebox{0.7}{$\left(\begin{smallmatrix}a&b\\a''&b''\end{smallmatrix}\right)$}}
$, in other terms, if $(f,f')\in \Phi^{\mathrm{fl}}_{\scalebox{0.7}{$\left(\begin{smallmatrix}a&b\\a'& b'\end{smallmatrix}\right)$}}$ and $(f',f'')\in \Phi^{\mathrm{fl}}_{\scalebox{0.7}{$\left(\begin{smallmatrix}a'&b'\\a''& b''\end{smallmatrix}\right)$}}$, then $(f,f'')\in \Phi^{\mathrm{fl}}_{\scalebox{0.7}{$\left(\begin{smallmatrix}a&b\\a''& b''\end{smallmatrix}\right)$}}$;
\item  for every $(a,a'),(b,b'), (c,c') \in \Phi^{\mathrm{ob}}$,
\[
\text{ if }
\left\lbrace
\begin{array}{c}
(f,f')\in \Phi^{\mathrm{fl}}_{\scalebox{0.7}{$\left(\begin{smallmatrix}a&b\\a'& b'\end{smallmatrix}\right)$}}\\[2pt]
{\text{and}}\\[2pt]
(g,g')\in \Phi^{\mathrm{fl}}_{\scalebox{0.7}{$\left(\begin{smallmatrix}b&c\\b'& c'\end{smallmatrix}\right)$}}
\end{array}
\right.,
\text{ then }
(g\circ f, g'\circ f')\in \Phi^{\mathrm{fl}}_{\scalebox{0.7}{$\left(\begin{smallmatrix}a&c\\a'& c'\end{smallmatrix}\right)$}};
\]
\item for every $(a,a'), (b,b') \in \Phi^{\mathrm{ob}}$ and every $f\in \mathrm{Mor}(\mathsf{C})$, if $f \in \mathrm{Hom}_{\mathsf{C}}(a,b)$, then there exists an $f' \in \mathrm{Hom}_{\mathsf{C}}(a',b')$ such that $(f,f') \in \Phi^{\mathrm{fl}}_{\scalebox{0.7}{$\left(\begin{smallmatrix}a&b\\a'& b'\end{smallmatrix}\right)$}}$. Sometimes we will refer to $(\mathrm{f})$ as the \emph{quadrangular completion condition}. 
\end{enumerate} 
\end{enumerate}
Conditions $(\mathrm{d})$ and $(\mathrm{e})$ are illustrated in the diagram shown in Figure~\ref{FRiguet}. For an object $a$ in $\mathsf{C}$, we let $[a]_{\Phi^{\mathrm{ob}}}$ stand for the $\Phi^{\mathrm{ob}}$-equivalence class of $a$. Moreover, we denote by $\mathrm{RCgr}(\mathsf{C})$ the set of all Riguet congruences on $\mathsf{C}$.
\end{definition}

\begin{figure}
\centering
\begin{tikzpicture}
[
ACliment/.style={-{To [angle'=45, length=5.75pt, width=4pt, round]}},
scale=.85
]
            
\node[] (a1) 	at (0,0) [] {$a$};
\node[] (b1) 	at (2,0) [] {$b$};
\node[] (c1) 	at (4,0) [] {$c$};
\node[] (a2) 	at (0,-1.5) [] {$a'$};
\node[] (b2) 	at (2,-1.5) [] {$b'$};
\node[] (c2) 	at (4,-1.5) [] {$c$};
\node[] (a3) 	at (0,-3) [] {$a''$};
\node[] (b3) 	at (2,-3) [] {$b''$};
\node[] (a4) 	at (0,-4.5) [] {$a$};
\node[] (b4) 	at (2,-4.5) [] {$b$};
\node[] (a5) 	at (0,-6) [] {$a''$};
\node[] (b5) 	at (2,-6) [] {$b''$};
\node[] (a6) 	at (6,0) [] {$a$};
\node[] (c3) 	at (8,0) [] {$c$};
\node[] (a7) 	at (6,-1.5) [] {$a'$};
\node[] (c4) 	at (8,-1.5) [] {$c'$};

\node[] (aux1) at (1,-3) [] {};
\node[] (aux2) at (1,-4.5) [] {};
\node[] (aux3) at (4,-0.75) [] {};
\node[] (aux4) at (6,-0.75) [] {};

\draw[ACliment, decorate, 
decoration={coil, aspect=0, amplitude=.2mm, segment length=1.5mm, post length=2mm}
] (aux1) to node [] {} (aux2);
\draw[ACliment, decorate, 
decoration={coil, aspect=0, amplitude=.2mm, segment length=1.5mm, post length=2mm}
] (aux3) to node [] {} (aux4);
            
\draw[ACliment]  (a1) to node [above]	{$\scriptstyle f$} (b1); 
\draw[ACliment]  (a2) to node [midway, fill=white] {$\scriptstyle f'$} (b2); 
\draw[ACliment]  (b1) to node [above]	{$\scriptstyle g$} (c1); 
\draw[ACliment]  (b2) to node [below] {$\scriptstyle g'$} (c2); 
\draw[ACliment]  (a3) to node [above] {$\scriptstyle f''$} (b3);  
\draw[ACliment]  (a4) to node [below]	{$\scriptstyle f$} (b4);  
\draw[ACliment]  (a5) to node [below]	{$\scriptstyle f''$} (b5);
\draw[ACliment]  (a6) to node [above]	{$\scriptstyle g\circ f$} (c3);
\draw[ACliment]  (a7) to node [below]	{$\scriptstyle g'\circ f'$} (c4);  
\foreach \x in {-0.07, 0, 0.07} {
\draw ([xshift=\x cm] a1.south) -- ([xshift=\x cm] a2.north);
\draw ([xshift=\x cm] b1.south) -- ([xshift=\x cm] b2.north);
\draw ([xshift=\x cm] c1.south) -- ([xshift=\x cm] c2.north);
\draw ([xshift=\x cm] a2.south) -- ([xshift=\x cm] a3.north);
\draw ([xshift=\x cm] b2.south) -- ([xshift=\x cm] b3.north);
\draw ([xshift=\x cm] a4.south) -- ([xshift=\x cm] a5.north);
\draw ([xshift=\x cm] b4.south) -- ([xshift=\x cm] b5.north);
\draw ([xshift=\x cm] a6.south) -- ([xshift=\x cm] a7.north);
\draw ([xshift=\x cm] c3.south) -- ([xshift=\x cm] c4.north);
}
\end{tikzpicture}
\caption{Conditions~$(\mathrm{d})$ and $(\mathrm{e})$ (vertical and horizontal compatibility) in Definition~\ref{DRiguet}.}
\label{FRiguet}
\end{figure}

\begin{remark}\label{rtccltcc}
Although condition $(\mathrm{a})$ of Definition~\ref{DRiguet} is not included in Riguet~\cite{JR60}, it is a natural compatibility condition ensuring that $\mathrm{id}_{a}$ and $\mathrm{id}_{b}$ become identified whenever $a$ and $b$ are identified.   
Likewise, condition $(\mathrm{f})$ of Definition~\ref{DRiguet} is absent from Riguet's work \cite{JR60}. However, this condition is required to define the composition of morphisms and, consequently, obtain a quotient category; this point is discussed in detail in Remark~\ref{CompOp}. 
Moreover, the quadrangular completion condition of Definition~\ref{DRiguet} is equivalent to each of the following conditions:  
\begin{enumerate}
\item[$(\mathrm{rtcc})$](Right triangular completion condition) For every $a$, $b$, $b'\in \mathrm{Ob}(\mathsf{C})$ and every $f\in \mathrm{Mor}(\mathsf{C})$, if $(b,b')\in \Phi^{\mathrm{ob}}$ and $f\in \mathrm{Hom}_{\mathsf{C}}(a,b)$, then there exists an $f' \in \mathrm{Hom}_{\mathsf{C}}(a,b')$ such that $(f,f')\in \Phi^{\mathrm{fl}}_{\scalebox{0.7}{$\left(\begin{smallmatrix}a&b\\a&b'\end{smallmatrix}\right)$}}$.
\item[$(\mathrm{ltcc})$](Left triangular completion condition) For every $b$, $b'$, $c\in \mathrm{Ob}(\mathsf{C})$ and every $g\in \mathrm{Mor}(\mathsf{C})$, if $(b,b')\in \Phi^{\mathrm{ob}}$ and $g\in \mathrm{Hom}_{\mathsf{C}}(b,c)$, then there exists a $g' \in \mathrm{Hom}_{\mathsf{C}}(b',c)$ such that $(g,g') \in \Phi^{\mathrm{fl}}_{\scalebox{0.7}{$\left(\begin{smallmatrix}b&c\\b'&c\end{smallmatrix}\right)$}}$. 
\end{enumerate}
In addition, on p.~9 of \cite{JR60}, Riguet describes $\Phi^{\mathrm{ob}}$ as a binary relation on $\mathrm{Ob}(\mathsf{C})$; however, he should have said that it is an equivalence relation on $\mathrm{Ob}(\mathsf{C})$.
\end{remark}

\begin{remark}\label{RNempHom}
If $\Phi$ is a Riguet congruence on a category $\mathsf{C}$, then, by condition $(\mathrm{f})$ in Definition~\ref{DRiguet}, for every $(a,b) \in \Phi^{\mathrm{ob}}$, $\mathrm{Hom}_{\mathsf{C}}(a,b)$ is nonempty. Moreover, if $\Phi$ is such that $\Phi^{\mathrm{ob}} = \Delta_{\mathrm{Ob}(\mathsf{C})}$, then condition $(\mathrm{f})$ reduces to the requirement that $\mathrm{Hom}_{\mathsf{C}}(a,a)$ is nonempty, which holds trivially for any object $a$ of $\mathsf{C}$.
\end{remark}

\begin{remark}
For any objects $a$, $a'$, $b$, $b'$ of $\mathsf{C}$ such that $\scalebox{0.7}{$\left(\begin{smallmatrix}a&b\\a'& b'\end{smallmatrix}\right)$} \in \Phi^{\mathrm{ob}}\times \Phi^{\mathrm{ob}}$, the component relation $\Phi^{\mathrm{fl}}_{\scalebox{0.7}{$\left(\begin{smallmatrix}a&b\\a'& b'\end{smallmatrix}\right)$}}$ is itself a difunctional relation as defined by Riguet in \cite{JR48}. Formally:
$\Phi^{\mathrm{fl}}_{\scalebox{0.7}{$\left(\begin{smallmatrix}a&b\\a'& b'\end{smallmatrix}\right)$}}\circ 
\left(\Phi^{\mathrm{fl}}_{\scalebox{0.7}{$\left(\begin{smallmatrix}a&b\\a'& b'\end{smallmatrix}\right)$}}\right)^{-1}\circ \Phi^{\mathrm{fl}}_{\scalebox{0.7}{$\left(\begin{smallmatrix}a&b\\a'& b'\end{smallmatrix}\right)$}}\subseteq \Phi^{\mathrm{fl}}_{\scalebox{0.7}{$\left(\begin{smallmatrix}a&b\\a'& b'\end{smallmatrix}\right)$}}$.
\end{remark}

In the following proposition we state that the set of all Riguet congruences on a category, ordered by inclusion, is a bounded directed-complete ordered set. The proof is immediate from the definition of a Riguet congruence and is omitted.

\begin{proposition}\label{RCgrAlgLatt}
Let $\mathsf{C}$ be a category. Then 
\begin{enumerate} 
\item $\Delta_{\mathsf{C}} = \left(\Delta_{\mathrm{Ob}(\mathsf{C})}, (\Delta_{\mathrm{Hom}_{\mathsf{C}}(a,b)})_{(a,b)\in \mathrm{Ob}(\mathsf{C})^{2}})\in\mathrm{RCgr}(\mathsf{C}\right)$ and it is the smallest Riguet congruence on $\mathsf{C}$ (where, to simplify the notation, we have written $(a,b)\in \mathrm{Ob}(\mathsf{C})^{2}$ instead of 
$\left(\begin{smallmatrix}a&b\\a& b\end{smallmatrix}\right)\in\Delta_{\mathrm{Ob}(\mathsf{C})}^{2}$).
\item $\Theta_{\mathsf{C}} = \left(\Theta_{\mathrm{Ob}(\mathsf{C})},(\mathrm{Hom}_{\mathsf{C}}(a,b) \times \mathrm{Hom}_{\mathsf{C}}(a',b'))_{\scalebox{0.7}{$\left(\begin{smallmatrix}a&b\\a'& b'\end{smallmatrix}\right)$}\in\Theta_{\mathrm{Ob}(\mathsf{C})}}\right)$, where $\Theta_{\mathrm{Ob}(\mathsf{C})}$ is the set of all $(a,b)\in \mathrm{Ob}(\mathsf{C})\times \mathrm{Ob}(\mathsf{C})$ such that $\mathrm{Hom}_{\mathsf{C}}(a,b)\neq \varnothing$ and $\mathrm{Hom}_{\mathsf{C}}(b,a)\neq \varnothing$, is in $\mathrm{RCgr}(\mathsf{C})$ and it is the greatest Riguet congruence on $\mathsf{C}$.
\item If $\mathcal{F}$ is a nonempty directed subset of $\mathrm{RCgr}(\mathsf{C})$, then $\bigcup \mathcal{F}\in\mathrm{RCgr}(\mathsf{C})$.
\end{enumerate}
Therefore, $\mathrm{RCgr}(\mathsf{C})$, ordered by inclusion, is a bounded directed-complete ordered set (bounded dcpo for short). That is, it is a dcpo possessing both a minimum and a maximum element. 
\end{proposition}

\begin{remark}
For the Riguet congruences on a category $\mathsf{C}$, the set operation $\bigcup$ and the inclusion relation are defined coordinatewise.
\end{remark}

We next state the relationship between the congruences on a category, as defined, e.g., by Mac Lane in~\cite{sM98}, and those of Riguet.

\begin{proposition}\label{MacLaneRiguet}
Let $\mathsf{C}$ be a category. Then there exists an embedding, i.e., an injective and isotone morphism, from $(\mathrm{Cgr}(\mathsf{C}),\subseteq)$, the algebraic lattice of all congruences on $\mathsf{C}$, into $(\mathrm{RCgr}(\mathsf{C}),\subseteq)$, which induces an embedding, i.e., a functor injective on objects and faithful, from $\mathsf{Cgr}(\mathsf{C})$, the category associated to $(\mathrm{Cgr}(\mathsf{C}),\subseteq)$, into $\mathsf{RCgr}(\mathsf{C})$, the category associated to $(\mathrm{RCgr}(\mathsf{C}),\subseteq)$. Hence 
$\mathsf{Cgr}(\mathsf{C})$ is isomorphic to its image in $\mathsf{RCgr}(\mathsf{C})$.
\end{proposition}

\begin{proof}
The mapping that assigns to a congruence $(\Phi_{a,b})_{(a,b)\in \mathrm{Ob}(\mathsf{C})^{2}}$ on $\mathsf{C}$ the Riguet congruence $(\Delta_{\mathrm{Ob}(\mathsf{C})},(\Phi_{a,b})_{(a,b)\in \mathrm{Ob}(\mathsf{C})^{2}})$ on $\mathsf{C}$ (where, to simplify, we have written $(a,b)\in \mathrm{Ob}(\mathsf{C})^{2}$ instead of $\left(\begin{smallmatrix}a&b\\a& b\end{smallmatrix}\right)\in\Delta_{\mathrm{Ob}(\mathsf{C})}^{2}$) is the desired embedding. 
\end{proof}

Before defining the quotient category associated with a category and a Riguet congruence on it, and stating its universal property, we first establish a preliminary result that will be needed to carry this out. The result follows immediately from the definitions, and we omit the proof.

\begin{proposition}\label{EquivDer}
Let $\mathsf{C}$ be a category and $\Phi$ a Riguet congruence on it. For every $a,b \in \mathrm{Ob}(\mathsf{C})$, let $\Phi^{\mathrm{fl}}_{[a]_{\Phi^{\mathrm{ob}}},[b]_{\Phi^{\mathrm{ob}}}}$ be $\bigcup_{\substack{a',a''\in [a]_{\Phi^{\mathrm{ob}}}\\b',b''\in [b]_{\Phi^{\mathrm{ob}}}}}\Phi^{\mathrm{fl}}_{\scalebox{0.7}{$\left(\begin{smallmatrix}a'&b'\\a''& b''\end{smallmatrix}\right)$}}$ and let $\mathrm{H}_{\mathsf{C}}([a]_{\Phi^{\mathrm{ob}}},[b]_{\Phi^{\mathrm{ob}}})$ be $\bigcup_{(a',b')\in [a]_{\Phi^{\mathrm{ob}}}\times [b]_{\Phi^{\mathrm{ob}}}}\mathrm{Hom}_{\mathsf{C}}(a',b')$. Then 
$\Phi^{\mathrm{fl}}_{[a]_{\Phi^{\mathrm{ob}}},[b]_{\Phi^{\mathrm{ob}}}}$ is an equivalence relation on 
$\mathrm{H}_{\mathsf{C}}([a]_{\Phi^{\mathrm{ob}}},[b]_{\Phi^{\mathrm{ob}}})$.  For a morphism $f$ in 
$\mathrm{H}_{\mathsf{C}}([a]_{\Phi^{\mathrm{ob}}},[b]_{\Phi^{\mathrm{ob}}})$, we let 
$[f]_{\Phi^{\mathrm{fl}}_{[a]_{\Phi^{\mathrm{ob}}},[b]_{\Phi^{\mathrm{ob}}}}}$ stand for the 
$\Phi^{\mathrm{fl}}_{[a]_{\Phi^{\mathrm{ob}}},[b]_{\Phi^{\mathrm{ob}}}}$-equivalence class of $f$. If there is no risk of confusion, we will write $[f]_{\Phi^{\mathrm{fl}}}$ instead of $[f]_{\Phi^{\mathrm{fl}}_{[a]_{\Phi^{\mathrm{ob}}},[b]_{\Phi^{\mathrm{ob}}}}}$.
\end{proposition}

\begin{remark}\label{PRiguetclass}
Let $f \colon a \mor b$ and $f' \colon a' \mor b'$ be morphisms in $\mathsf{C}$. Then we have that
$[f']_{\Phi^{\mathrm{fl}}_{[a']_{\Phi^{\mathrm{ob}}},[b']_{\Phi^{\mathrm{ob}}}}} = [f]_{\Phi^{\mathrm{fl}}_{[a]_{\Phi^{\mathrm{ob}}},[b]_{\Phi^{\mathrm{ob}}}}}$ if, and only if, 
 $(a,a')$,  $(b,b') \in \Phi^{\mathrm{ob}}$ and $(f,f') \in \Phi^{\mathrm{fl}}_{\scalebox{0.7}{$\left(\begin{smallmatrix}a&b\\a'& b'\end{smallmatrix}\right)$}}$.
\end{remark}

In the following definition, for complete precision in the definition of the composition operation---a point requiring some care---we retain the non-abbreviated notation for equivalence classes of morphisms.

\begin{definition}\label{DRiguetQ}
Let $\mathsf{C}$ be a category and $\Phi$ a Riguet congruence on it. Then the \emph{Riguet quotient} of $\mathsf{C}$ by $\Phi$, denoted by $\mathsf{C}/{\Phi}$, is the category defined as follows:
\begin{enumerate}
\item $\mathrm{Ob}(\mathsf{C}/{\Phi})=\mathrm{Ob}(\mathsf{C})/{\Phi^{\mathrm{ob}}}$;
\item for every $a,b \in \mathrm{Ob}(\mathsf{C})$, the set of morphisms from $[a]_{\Phi^{\mathrm{ob}}}$ to  $[b]_{\Phi^{\mathrm{ob}}}$ is
\[
\mathrm{Hom}_{\mathsf{C}/{\Phi}} ([a]_{\Phi^{\mathrm{ob}}}, [b]_{\Phi^{\mathrm{ob}}}) = \mathrm{H}_{\mathsf{C}}([a]_{\Phi^{\mathrm{ob}}},[b]_{\Phi^{\mathrm{ob}}})/\Phi^{\mathrm{fl}}_{[a]_{\Phi^{\mathrm{ob}}},[b]_{\Phi^{\mathrm{ob}}}};
\]
\item for every $a$, $b$, $c\in \mathrm{Ob}(\mathsf{C})$, the composition operation 
$\circ_{[a]_{\Phi^{\mathrm{ob}}}, [b]_{\Phi^{\mathrm{ob}}},[c]_{\Phi^{\mathrm{ob}}}}$ associated to $([a]_{\Phi^{\mathrm{ob}}}, [b]_{\Phi^{\mathrm{ob}}},[c]_{\Phi^{\mathrm{ob}}})$ , also denoted by $\circ$ for short, 
\[
\xymatrix@C=20pt{
\mathrm{Hom}_{\mathsf{C}/{\Phi}} ([a]_{\Phi^{\mathrm{ob}}}, [b]_{\Phi^{\mathrm{ob}}})\times 
\mathrm{Hom}_{\mathsf{C}/{\Phi}} ([b]_{\Phi^{\mathrm{ob}}}, [c]_{\Phi^{\mathrm{ob}}})
\ar[r]^-{\circ} 
  & \mathrm{Hom}_{\mathsf{C}/{\Phi}} ([a]_{\Phi^{\mathrm{ob}}}, [c]_{\Phi^{\mathrm{ob}}}) }
\]
is defined as: 
for $\mathfrak{f}\in \mathrm{Hom}_{\mathsf{C}/{\Phi}}([a]_{\Phi^{\mathrm{ob}}}, [b]_{\Phi^{\mathrm{ob}}})$, i.e., if $\mathfrak{f} = [f]_{\Phi^{\mathrm{fl}}_{[a]_{\Phi^{\mathrm{ob}}},[b]_{\Phi^{\mathrm{ob}}}}}$, for some $a'\in [a]_{\Phi^{\mathrm{ob}}}$, $b'\in [b]_{\Phi^{\mathrm{ob}}}$ and $f\in \mathrm{Hom}_{\mathsf{C}}(a',b')$, and for
$\mathfrak{g}\in \mathrm{Hom}_{\mathsf{C}/{\Phi}} ([b]_{\Phi^{\mathrm{ob}}},[c]_{\Phi^{\mathrm{ob}}})$, i.e., if $\mathfrak{g} = [g]_{\Phi^{\mathrm{fl}}_{[b]_{\Phi^{\mathrm{ob}}},[c]_{\Phi^{\mathrm{ob}}}}}$, for some $b''\in [b]_{\Phi^{\mathrm{ob}}}$, $c'\in [c]_{\Phi^{\mathrm{ob}}}$ and $g\in \mathrm{Hom}_{\mathsf{C}}(b'',c')$, then the composite of $\mathfrak{f}$ and $\mathfrak{g}$ is given by 
$\mathfrak{g}\circ \mathfrak{f} = [g\circ h]_{\Phi^{\mathrm{fl}}_{[a]_{\Phi^{\mathrm{ob}}},[c]_{\Phi^{\mathrm{ob}}}}}$, where, by condition $(\mathrm{f})$ in Definition~\ref{DRiguet}, 
$h\in \mathrm{Hom}_{\mathsf{C}}(a',b'')$ is such that $(f,h)\in\Phi^{\mathrm{fl}}_{\scalebox{0.7}{$\left(\begin{smallmatrix}a'&b'\\a'& b''\end{smallmatrix}\right)$}}$;
\item for every $a\in \mathrm{Ob}(\mathsf{C})$, $\mathrm{id}_{[a]_{\Phi^{\mathrm{ob}}}}$, the identity at $[a]_{\Phi^{\mathrm{ob}}}$, is $[\mathrm{id}_{a'}]_{\Phi^{\mathrm{fl}}_{[a']_{\Phi^{\mathrm{ob}}},[a']_{\Phi^{\mathrm{ob}}}}}$, for some $a'\in [a]_{\Phi^{\mathrm{ob}}}$.
\end{enumerate}
\end{definition}

\begin{remark}\label{CompOp}
A comparison between Definition~\ref{DRiguetQ} and Riguet's definition of composition in~\cite{JR60}, p.~11, shows that his original formulation fails to define a valid composition operation. In contrast, our condition $(\mathrm{f})$ guarantees that the operation is well-defined, allowing for the construction of a quotient category.  Indeed, using the notation of Definition~\ref{DRiguetQ}, let $\mathfrak{f} = [f]_{\Phi^{\mathrm{fl}}_{[a]_{\Phi^{\mathrm{ob}}},[b]_{\Phi^{\mathrm{ob}}}}}$, for some $a'\in [a]_{\Phi^{\mathrm{ob}}}$, $b'\in [b]_{\Phi^{\mathrm{ob}}}$ and $f\in \mathrm{Hom}_{\mathsf{C}}(a',b')$, and $\mathfrak{g} = [g]_{\Phi^{\mathrm{fl}}_{[b]_{\Phi^{\mathrm{ob}}},[c]_{\Phi^{\mathrm{ob}}}}}$, for some $b''\in [b]_{\Phi^{\mathrm{ob}}}$, $c'\in [c]_{\Phi^{\mathrm{ob}}}$ and $g\in \mathrm{Hom}_{\mathsf{C}}(b'',c')$. In Riguet's formulation, the composition of $ \mathfrak{f}$ and $\mathfrak{g}$ is given by $\mathfrak{g}\circ \mathfrak{f} = [g\circ f]_{\Phi^{\mathrm{fl}}_{[a]_{\Phi^{\mathrm{ob}}},[c]_{\Phi^{\mathrm{ob}}}}}$. However, since $b'$ and $b''$ are only required to belong to the same equivalence class 
$[b]_{\Phi^{\mathrm{ob}}}$, there is no guarantee that the composite $g\circ f$ exists in $\mathsf{C}$. Thus, without additional assumptions, the composition of equivalence classes need not be defined. Condition~$(\mathrm{f})$ ensures the existence of compatible representatives whose composite is defined, and therefore guarantees that the composition in the quotient category is well defined. 
\end{remark}

We next define the notions of Riguet-full functor and of Riguet kernel of a functor.

\begin{definition}
Let $F$ be a functor from $\mathsf{C}$ to $\mathsf{D}$. We will call $F$ a \emph{Riguet-full functor} if, for every $a$, $a'$, $b$, $b'\in \mathrm{Ob}(\mathsf{C})$ and every $f\in \mathrm{Hom}_{\mathsf{C}}(a,b)$, if $F(a) = F(a')$ and $F(b) = F(b')$, then there exists an $f'\in \mathrm{Hom}_{\mathsf{C}}(a',b')$ such that $F(f) = F(f')$. We denote by $\mathrm{K}^{\mathrm{r}}(F)$ the ordered pair  
\[
\left(\mathrm{K}^{\mathrm{r}}(F)^{\mathrm{ob}}, \left(\mathrm{K}^{\mathrm{r}}(F)^{\mathrm{fl}}_{\scalebox{0.7}{$\left(\begin{smallmatrix}a&b\\a'& b'\end{smallmatrix}\right)$}}\right)_{\scalebox{0.7}{$\left(\begin{smallmatrix}a&b\\a'& b'\end{smallmatrix}\right)$} \in \mathrm{K}^{\mathrm{r}}(F)^{\mathrm{ob}}\times \mathrm{K}^{\mathrm{r}}(F)^{\mathrm{ob}}}\right), 
\]
where $\mathrm{K}^{\mathrm{r}}(F)^{\mathrm{ob}}$ is the kernel of the object mapping of $F$ and, for every $\scalebox{0.7}{$\left(\begin{smallmatrix}a&b\\a'& b'\end{smallmatrix}\right)$} \in \mathrm{K}^{\mathrm{r}}(F)^{\mathrm{ob}}\times \mathrm{K}^{\mathrm{r}}(F)^{\mathrm{ob}}$, $\mathrm{K}^{\mathrm{r}}(F)^{\mathrm{fl}}_{\scalebox{0.7}{$\left(\begin{smallmatrix}a&b\\a'& b'\end{smallmatrix}\right)$}}$ is the subset of $\mathrm{Hom}_{\mathsf{C}}(a,b) \times \mathrm{Hom}_{\mathsf{C}}(a',b')$ defined, for every $(f,f')\in \mathrm{Hom}_{\mathsf{C}}(a,b) \times \mathrm{Hom}_{\mathsf{C}}(a',b')$, as: $(f,f')\in \mathrm{K}^{\mathrm{r}}(F)^{\mathrm{fl}}_{\scalebox{0.7}{$\left(\begin{smallmatrix}a&b\\a'& b'\end{smallmatrix}\right)$}}$ if, and only if, $F(f) = F(f')$, and we will call it the \emph{Riguet kernel} of $F$ or, simply, the \emph{kernel} of $F$. Note that $\mathrm{K}^{\mathrm{r}}(F)^{\mathrm{ob}}$ is an equivalence relation on $\mathrm{Ob}(\mathsf{C})$ and that $\mathrm{K}^{\mathrm{r}}(F)^{\mathrm{fl}}$ satisfies the conditions $(\mathrm{a})$--$(\mathrm{e})$ of a Riguet congruence. 
\end{definition}

\begin{proposition} \label{RiguetFullKernel}
Let $F\colon \mathsf{C} \mor \mathsf{D}$ be a functor. Then $\mathrm{K}^{\mathrm{r}}(F)$ is a Riguet congruence on $\mathsf{C}$ if, and only if, $F$ is a Riguet-full functor.
\end{proposition}

\begin{proof}

$(\Rightarrow)$ Suppose that $\mathrm{K}^{\mathrm{r}}(F)$ is a Riguet congruence. Let $a, a', b, b' \in \mathrm{Ob}(\mathsf{C})$ such that $F(a)=F(a')$ and $F(b)=F(b')$. By the definition of $\mathrm{K}^{\mathrm{r}}(F)$, this implies $(a, a') \in \mathrm{K}^{\mathrm{r}}(F)^{\mathrm{ob}}$ and $(b, b') \in \mathrm{K}^{\mathrm{r}}(F)^{\mathrm{ob}}$. Let $f\colon a \mor b$ be a morphism. Then, by condition (f), there exists an $f'\colon a' \mor b'$ such that $(f, f') \in \mathrm{K}^{\mathrm{r}}(F)^{\mathrm{fl}}$. Hence, by the definition of $\mathrm{K}^{\mathrm{r}}(F)$, we have that $F(f) = F(f')$, which proves that $F$ is Riguet-full.

$(\Leftarrow)$ Suppose that $F$ is Riguet-full. To verify condition (f) for $\mathrm{K}^{\mathrm{r}}(F)$, let $(a, a')$, $(b, b') \in \mathrm{K}^{\mathrm{r}}(F)^{\mathrm{ob}}$, and $f\colon a \mor b$. These relations imply $F(a)=F(a')$ and $F(b)=F(b')$. The Riguet-fullness of $F$ guarantees the existence of a morphism $f'\colon a' \mor b'$ such that $F(f)=F(f')$, that is, $(f, f') \in \mathrm{K}^{\mathrm{r}}(F)^{\mathrm{fl}}$, thereby satisfying condition (f).
\end{proof}

In the following proposition, we define the canonical projection relative to a Riguet congruence on a category. The proof is immediate from the definitions and is omitted.

\begin{proposition}\label{PProj}
Let $\mathsf{C}$ be a category and $\Phi$ a Riguet congruence on it. Then we let 
$P_{\Phi}$ stand for the assignment from $\mathsf{C}$ to $\mathsf{C}/{\Phi}$ defined as follows:
\begin{enumerate}
\item for every object $a$ in $\mathsf{C}$, $P_{\Phi}(a)=[a]_{\Phi^{\mathrm{ob}}}$, and
\item for every mapping $f\colon a \mor b$ in $\mathsf{C}$, $P_{\Phi}(f)=[f]_{\Phi^{\mathrm{fl}}}$, as stated in Proposition~\ref{EquivDer}.
\end{enumerate}
Then $P_{\Phi}$ is a functor from $\mathsf{C}$ to $\mathsf{C}/{\Phi}$. We will call this functor the \emph{canonical projection relative to $\Phi$}. 
\end{proposition}

In the following proposition we state the universal property of the quotient category of a $\mathsf{C}$ by a Riguet congruence $\Phi$ on it. The proof is standard and hence omitted.

\begin{proposition}\label{UPQ}
Let $\mathsf{C}$ be a category, $\Phi$ a Riguet congruence on it and $F$ a functor from $\mathsf{C}$ to $\mathsf{D}$ such that $\Phi\subseteq \mathrm{K}^{\mathrm{r}}(F)$, i.e., (1) for every $a$, $a'\in \mathrm{Ob}(\mathsf{C})$, if $(a,a')\in \Phi^{\mathrm{ob}}$, then $F(a) = F(a')$ and (2) for every $f \colon a \mor b$ and $f' \colon a' \mor b'$ in $\mathsf{C}$, if $(f,f') \in \Phi^{\mathrm{fl}}_{\scalebox{0.7}{$\left(\begin{smallmatrix}a&b\\a'& b'\end{smallmatrix}\right)$}}$, then $F(f) = F(f')$. Then there exists a unique functor $F^{\sharp}$ from $\mathsf{C}/\Phi$ to $\mathsf{D}$ such that $F^{\sharp}\circ P_{\Phi} = F$. 
\end{proposition}

We next prove that the canonical projection $P_{\Phi}\colon \mathsf{C} \mor \mathsf{C}/\Phi$ is a regular epimorphism, i.e., the coequalizer of some parallel pair of functors.

\begin{proposition}\label{RCRegEpi}
Let $\Phi$ be a Riguet congruence on a category $\mathsf{C}$. Then the canonical projection $P_{\Phi}\colon \mathsf{C} \mor \mathsf{C}/\Phi$ is a regular epimorphism and a Riguet-full functor.
\end{proposition}

\begin{proof}
It suffices to consider the category $\mathsf{C}(\Phi)$ associated to $\Phi$, whose set of objects is $\Phi^{\mathrm{ob}} \subseteq \mathrm{Ob}(\mathsf{C}) \times \mathrm{Ob}(\mathsf{C})$ and, for any $(a,a'), (b,b') \in \Phi^{\mathrm{ob}}$, the hom-set is defined as
\[
\mathrm{Hom}_{\mathsf{C}(\Phi)}((a,a'), (b,b')) = \Phi^{\mathrm{fl}}_{\scalebox{0.7}{$\left(\begin{smallmatrix}a&b\\a'& b'\end{smallmatrix}\right)$}}.
\]
The identities and composition are inherited from $\mathsf{C} \times \mathsf{C}$, as ensured by conditions (a) and (e) of Definition~\ref{DRiguet}, respectively. This category comes equipped with two projection functors $P_{0}, P_{1} \colon \mathsf{C}(\Phi) \mor \mathsf{C}$ which map each morphism $(f, f')$ to $f$ and $f'$, respectively, and  $P_{\Phi}$ is the coequalizer of $(\mathsf{C}(\Phi), (P_{0}, P_{1}))$.
\end{proof}

\begin{proposition}\label{FSOReg}
Let $F\colon \mathsf{C} \mor \mathsf{D}$ be a functor. If $F$ is full and surjective on objects, then $F$ is a regular epimorphism. 
\end{proposition} 

\begin{proof}
The functor $F$ is the coequalizer of its kernel pair $(\mathsf{KP}(F),(P_{0},P_{1}))$, where $\mathsf{KP}(F)$ is the pullback of $F$ along itself and $P_{0}$, $P_{1}$ the projections from $\mathsf{KP}(F)$ to $\mathsf{C}$. We recall that
\begin{align*}
\mathrm{Ob}(\mathsf{KP}(F)) &= \{(x,y)\in \mathrm{Ob}(\mathsf{C})^{2}\mid F(x) = F(y)\} \text{ and} 
\\ 
\mathrm{Mor}(\mathsf{KP}(F)) &= \{(f,g)\in \mathrm{Mor}(\mathsf{C})^{2}\mid F(f) = F(g)\}.\qedhere
\end{align*}
\end{proof}

We next characterize functors that are full and surjective on objects in terms of extremal functors whose Riguet kernel is a Riguet congruence.

\begin{proposition}
A functor $F\colon  \mathsf{C} \mor \mathsf{D}$ is full and surjective on objects if, and only if, it is an extremal epimorphism and Riguet-full ($\equiv$ $\mathrm{K}^{\mathrm{r}}(F)$ is a Riguet congruence on $\mathsf{C}$).
\end{proposition}

\begin{proof}
$(\Rightarrow)$ If $F$ is full and surjective on objects, then, by Proposition~\ref{FSOReg}, it is a regular epimorphism. Hence it is extremal. Let $a$, $a'$, $b$, $b'\in \mathrm{Ob}(\mathsf{C})$ and $f\in \mathrm{Hom}_{\mathsf{C}}(a,b)$, such that $F(a) = F(a')$ and $F(b) = F(b')$. Then $F(f)\in 
\mathrm{Hom}_{\mathsf{D}}(F(a'),F(b'))$. Hence there exists an $f'\in \mathrm{Hom}_{\mathsf{C}}(a',b')$ such that $F(f) = F(f')$.

$(\Leftarrow)$ Let us suppose that $F$ is an extremal epimorphism and a Riguet-full functor. Then, because a functor $F\colon\mathsf{C}\mor\mathsf{D}$ is an extremal epimorphism if, and only if, $\mathsf{Sg}_{\mathsf{D}}(\mathrm{Im}(F)) = \mathsf{D}$, where $\mathsf{Sg}_{\mathsf{D}}(\mathrm{Im}(F))$ denotes the subcategory of $\mathsf{D}$ generated by $\mathrm{Im}(F)$, we have that $F$ is surjective on objects. We now prove that $F$ is full. Let $a$ and $b$ be objects of $\mathsf{C}$ and $u$ a morphism in $\mathsf{D}$ from $F(a)$ to $F(b)$. Since $\mathsf{Sg}_{\mathsf{D}}(\mathrm{Im}(F)) = \mathsf{D}$, there are a natural number $n\geq 1$, two families of objects $(a_{i})_{i\in n}$, $(b_{i})_{i\in n}$ in $\mathsf{C}$ and a family of morphisms $(f_{i})_{i\in n}$ in $\mathsf{C}$ such that, for every $i\in n$, $f_{i}\colon a_{i}\mor b_{i}$, for every $i\in n-1$, $F(b_{i}) = F(a_{i+1})$ and $u = F(f_{n-1})\circ \cdots\circ F(f_{0})$. We proceed by induction on $n$. For $n=1$, $u = F(f_{0})$ with $f_{0}\colon a_{0}\mor b_{0}$. By Riguet-fullness, since $F(a_{0})=F(a)$ and $F(b_{0})=F(b)$, there exists an $f\colon a\mor b$ such that  $F(f) = u = F(f_{0})$. Assume it holds for $n\geq 1$. For $n+1$, let $u = F(f_{n}) \circ (F(f_{n-1}) \circ \dots \circ F(f_{0}))$. By hypothesis, the tail is $F(h)$ for some $h\colon  a_{0} \mor b_{n-1}$. Since $F(b_{n-1}) = F(a_{n})$, Riguet-fullness provides $h'\colon  a_{0} \mor a_{n}$ such that $F(h')=F(h)$. Then $u = F(f_{n} \circ h')$. Therefore $F$ is full.
\end{proof}

\begin{remark}
Bednarczyk \textit{et al.}~\cite{BBP99} prove that $F\colon\mathsf{C}\mor\mathsf{D}$ is a regular  epimorphism if, and only if, $F$ is an extremal epimorphism and the generalized kernel $\mathrm{K}(F)$ of $F$ is regular.
\end{remark}

We next characterize functors that are full and surjective on objects in terms
of regular epimorphic functors whose Riguet kernel is a Riguet congruence. Before
proceeding, we recall the construction of the coequalizer of a pair of parallel functors
with common source and target and prove that, for a Riguet-full functor, every non-
identity morphism in the coequalizer of its kernel pair coincides with the equivalence
class determined by the sequence of length one associated to a morphism in the
domain of the functor.

\begin{lemma}\label{Coeq}
Let $F$ and $G$ be two functors from a category $\mathsf{A}$ to another $\mathsf{C}$. Then there exists the coequalizer of $F$ and $G$.
\end{lemma}

\begin{proof}
We denote by $\mathsf{Coeq}(F,G)$ the category defined as follows. The set of objects of $\mathsf{Coeq}(F,G)$, denoted by $\mathrm{Ob}(\mathsf{Coeq}(F,G))$, is the quotient of $\mathrm{Ob}(\mathsf{C})$ by the smallest equivalence relation $\sim^{\mathrm{ob}}$ on $\mathrm{Ob}(\mathsf{C})$ such that, for every $a \in \mathrm{Ob}(\mathsf{A})$, $F(a) \sim^{\mathrm{ob}} G(a)$. For an object $c$ in $\mathsf{C}$, we let $[c]_{\sim^{\mathrm{ob}}}$ stand for the $\sim^{\mathrm{ob}}$-equivalence class of $c$. Moreover, we denote by $\mathrm{coeq}(F,G)$ the canonical projection from $\mathrm{Ob}(\mathsf{C})$ to $\mathrm{Ob}(\mathsf{Coeq}(F,G))$. Note that $\mathrm{coeq}(F,G)$ will be the object mapping of the canonical functor from $\mathsf{C}$ to $\mathsf{Coeq}(F,G)$.

To define the morphisms of $\mathsf{Coeq}(F,G)$ we proceed as follows. Let $\mathsf{G}_{F,G}$ be the quiver defined as follows. The set of vertices of $\mathsf{G}_{F,G}$ is the set $\mathrm{Ob}(\mathsf{Coeq}(F,G)) = \mathrm{Ob}(\mathsf{C})/{\sim^{\mathrm{ob}}}$. For any pair of vertices $[c]_{\sim^{\mathrm{ob}}}$ and $[d]_{\sim^{\mathrm{ob}}}$ of $\mathsf{G}_{F,G}$, the set of edges from $[c]_{\sim^{\mathrm{ob}}}$ to $[d]_{\sim^{\mathrm{ob}}}$ is
\[
\textstyle
\mathrm{Hom}_{\mathsf{G}_{F,G}}
\big([c]_{\sim^{\mathrm{ob}}},[d]_{\sim^{\mathrm{ob}}}\bigr)
=
\bigcup_{\substack{a\in[c]_{\sim^{\mathrm{ob}}}\\ b\in[d]_{\sim^{\mathrm{ob}}}}}
\mathrm{Hom}_{\mathsf{C}}(a,b).
\]
Note that the edge-sets $\mathrm{Hom}_{\mathsf{G}_{F,G}}([c]_{\sim^{\mathrm{ob}}},[d]_{\sim^{\mathrm{ob}}})$ are pairwise disjoint. 
Now we consider the free category $\mathsf{L}(\mathsf{G}_{F,G})$ generated by this quiver. The objects in $\mathsf{L}(\mathsf{G}_{F,G})$ are the elements of $\mathrm{Ob}(\mathsf{Coeq}(F,G))$. The morphisms in $\mathsf{L}(\mathsf{G}_{F,G})$ are finite sequences of edges $(f_{i})_{i\in n}$ in the quiver $\mathsf{G}_{F,G}$, which we also denote by $(f_{0},\ldots,f_{n-1})$, such that $\mathrm{tg}(f_{i}) = [b_{i}]_{\sim^{\mathrm{ob}}}$ is $\mathrm{sc}(f_{i+1}) = [a_{i+1}]_{\sim^{\mathrm{ob}}}$. We write $\mathrm{id}_{[c]_{\sim^{\mathrm{ob}}}}$ for the empty sequence of edges at the object $[c]_{\sim^{\mathrm{ob}}}$, which is the identity morphism at $[c]_{\sim^{\mathrm{ob}}}$. Sometimes it is convenient to regard the morphisms in $\mathsf{L}(\mathsf{G}_{F,G})$ as ordered pairs $(([c_{i}]_{\sim^{\mathrm{ob}}})_{i\in n+1},(f_{i})_{i\in n})$ such that, for every $i\in n$, $\mathrm{tg}(f_{i}) = [c_{i+1}]_{\sim^{\mathrm{ob}}} = \mathrm{sc}(f_{i+1})$. Note that, in this case, $\mathrm{id}_{[c]_{\sim^{\mathrm{ob}}}}$ is the ordered pair 
$(([c]_{\sim^{\mathrm{ob}}}),\lambda)$, where $\lambda$ is the empty word of the free monoid on the set of edges of $\mathsf{G}_{F,G}$.  
The morphisms of $\mathsf{Coeq}(F,G)$ are the morphisms of the quotient category $\mathsf{L}(\mathsf{G}_{F,G})/{{\sim^{\mathrm{fl}}}}$, where ${\sim^{\mathrm{fl}}} = 
({\sim^{\mathrm{fl}}}_{[a]_{\sim^{\mathrm{ob}}},[b]_{\sim^{\mathrm{ob}}}})_{([a]_{\sim^{\mathrm{ob}}},[b]_{\sim^{\mathrm{ob}}})\in \mathrm{Ob}(\mathsf{Coeq}(F,G))^{2}}$ is the smallest congruence (in the sense of Mac Lane~\cite{sM98}) on $\mathsf{L}(\mathsf{G}_{F,G})$ that satisfies: (1) $(f_{0}, f_{1}) {\sim^{\mathrm{fl}}} (f_{1} \circ f_{0})$ whenever $f_{1} \circ f_{0}$ is defined in $\mathsf{C}$; (2) $(\mathrm{id}_{c}) {\sim^{\mathrm{fl}}} \mathrm{id}_{[c]_{\sim^{\mathrm{ob}}}}$; and (3) $(F(h)) {\sim^{\mathrm{fl}}} (G(h))$ for every morphism $h$ in $\mathsf{A}$. Thus, 
$\mathsf{Coeq}(F,G)$ is $\mathsf{L}(\mathsf{G}_{F,G})/{{\sim^{\mathrm{fl}}}}$. 
Moreover, with the customary abuse of notation, the same symbol is used for the mapping $\mathrm{coeq}(F,G)$ from $\mathrm{Mor}(\mathsf{C})$ to $\mathrm{Mor}(\mathsf{Coeq}(F,G))$ that sends a morphism $f$ in $\mathsf{C}$ to $\mathrm{coeq}(F,G)(f) = [(f)]_{{\sim^{\mathrm{fl}}}}$. The proof of the functoriality of $\mathrm{coeq}(F,G)$ and of the universal property are standard and hence omitted. We merely note the following. First, $\mathrm{coeq}(F,G)$ is the composite
\[
\xymatrix@C=30pt{ \mathsf{C}
\ar[r]^-{P_{\sim^{\mathrm{ob}}}} 
&
\mathsf{G}_{F,G} 
\ar[r]^-{\eta_{\mathsf{G}_{F,G}}} 
&
\mathsf{L}(\mathsf{G}_{F,G})
\ar[r]^-{P_{{\sim^{\mathrm{fl}}}}}
&
\mathsf{L}(\mathsf{G}_{F,G})/\!{\sim^{\mathrm{fl}}}, 
}
\]
where $P_{\sim^{\mathrm{ob}}}$, $\eta_{\mathsf{G}_{F,G}}$, and
$P_{{\sim^{\mathrm{fl}}}}$ denote, respectively, the canonical morphism
from the underlying quiver of $\mathsf{C}$ to the quiver
$\mathsf{G}_{F,G}$---that sends $c$ to $[c]_{\sim^{\mathrm{ob}}}$ and $f\colon c\mor d$ to 
$f\colon [c]_{\sim^{\mathrm{ob}}}\mor [d]_{\sim^{\mathrm{ob}}}$---, the canonical morphism from
$\mathsf{G}_{F,G}$ to the underlying quiver of
$\mathsf{L}(\mathsf{G}_{F,G})$, and the canonical functor from
$\mathsf{L}(\mathsf{G}_{F,G})$ to the quotient category
$\mathsf{L}(\mathsf{G}_{F,G})/\!{\sim^{\mathrm{fl}}}$. Second, for a functor $H$ from $\mathsf{C}$ to $\mathsf{E}$ such that $H\circ F = H\circ G$, the unique functor $H^{\sharp}$ from $\mathsf{Coeq}(F,G)$ to $\mathsf{E}$ such that $H^{\sharp}\circ \mathrm{coeq}(F,G) = H$, is obtained by first stating the existence of a unique morphism $H^{\natural}$ from  $\mathsf{G}_{F,G}$ to the underlying quiver of $\mathsf{E}$ such that $H^{\natural}\circ P_{\sim^{\mathrm{ob}}} = H$, then stating the existence of a unique functor $H^{\flat}$ from $\mathsf{L}(\mathsf{G}_{F,G})$ to $\mathsf{E}$ such that $H^{\flat}\circ \eta_{\mathsf{G}_{F,G}} = H^{\natural}$ and, finally, stating the existence of a unique functor $H^{\sharp}$ from $\mathsf{L}(\mathsf{G}_{F,G})/\!\sim^{\mathrm{fl}}$ to $\mathsf{E}$ such that $H^{\sharp}\circ P_{\sim^{\mathrm{fl}}} = H^{\flat}$.
\end{proof}

\begin{proposition}\label{red}
Let $F\colon \mathsf{C} \mor \mathsf{D}$ be a functor, $(\mathsf{KP}(F),(P_{0},P_{1}))$ its kernel pair and $\mathrm{coeq}(P_{0},P_{1})\colon \mathsf{C} \mor \mathsf{Coeq}(P_{0},P_{1})$ the coequalizer of $(P_{0},P_{1})$. If $F$ is Riguet-full, then for every non-identity morphism 
$\mathfrak{f} \in \mathrm{Mor}(\mathsf{Coeq}(P_{0},P_{1}))$, there exists a morphism $f \in \mathrm{Mor}(\mathsf{C})$ such that $\mathrm{coeq}(P_{0},P_{1})(f) = \mathfrak{f}$.
\end{proposition} 

\begin{proof}
By Lemma~\ref{Coeq}, a morphism $\mathfrak{f} \in \mathrm{Mor}(\mathsf{Coeq}(P_{0},P_{1}))$ is an equivalence class $\mathfrak{f} = [(f_{i})_{i\in n}]_{{\sim^{\mathrm{fl}}}}$ such that, for every $i\in n-1$, $\mathrm{tg}(f_{i}) = [b_{i}]_{\sim^{\mathrm{ob}}}$ is $\mathrm{sc}(f_{i+1}) = [a_{i+1}]_{\sim^{\mathrm{ob}}}$. 

The nontrivial base case is that of morphisms of length $2$. Let $[(f_{i})_{i\in 2}]_{{\sim^{\mathrm{fl}}}}$ be a morphism  of length $2$, with $f_{0}\colon a_{0} \mor b_{0}$ and $f_{1}\colon a_{1} \mor b_{1}$, such that $F(b_{0}) = F(a_{1})$. Then,  by the definition of Riguet-fullness, we have that for $f_{0}\colon a_{0} \mor b_{0}$, since $F(a_{0}) = F(a_{0})$ and $F(b_{0}) = F(a_{1})$, there exists an $f_{0}'\colon a_{0} \mor a_{1}$ in $\mathsf{C}$ such that $F(f_{0}') = F(f_{0})$. In $\mathsf{Coeq}(P_{0},P_{1})$ we have that:
\begin{enumerate}
\item Since $F(f_{0}') = F(f_{0})$, we have that $[(f_{0}')]_{{\sim^{\mathrm{fl}}}} = [(f_{0})]_{{\sim^{\mathrm{fl}}}}$.
\item The composition $f_{1} \circ f_{0}'$ is well-defined in $\mathsf{C}$.
\item By the properties of $\mathrm{coeq}(P_{0},P_{1})$:
\[
[(f_{i})_{i\in 2}]_{{\sim^{\mathrm{fl}}}} = 
[(f_{1})]_{{\sim^{\mathrm{fl}}}} \circ [(f_{0})]_{{\sim^{\mathrm{fl}}}} = 
[(f_{1})]_{{\sim^{\mathrm{fl}}}} \circ [(f_{0}')]_{{\sim^{\mathrm{fl}}}} = 
[(f_{1} \circ f_{0}')]_{{\sim^{\mathrm{fl}}}}.
\]
\end{enumerate}
Thus, any morphism of length $2$ is equivalent to a morphism of length $1$.

Assume that any morphism of length $n\geq 2$ is equivalent to a morphism of length $1$. Let 
$\mathfrak{f} = [(f_{i})_{i\in n+1}]_{{\sim^{\mathrm{fl}}}}$ be a morphism of length $n+1$.
By the inductive hypothesis, $[(f_{i})_{i\in n}]_{{\sim^{\mathrm{fl}}}} = [(h)]_{{\sim^{\mathrm{fl}}}}$ for some $h\colon a_{0} \mor b_{n-1}$. But, $F(b_{n-1}) = F(a_{n})$ or, what is equivalent, $[b_{n-1}]_{\sim^{\mathrm{ob}}} = [a_{n}]_{\sim^{\mathrm{ob}}}$.
Hence, for the pair $(f_{n}, h)$, there exists $h'\colon a_{0} \mor a_{n}$ such that $F(h') = F(h)$.
Then $\mathfrak{f} = [(f_{n})]_{{\sim^{\mathrm{fl}}}} \circ [(h')]_{{\sim^{\mathrm{fl}}}} = [(f_{n} \circ h')]_{{\sim^{\mathrm{fl}}}}$. Since $f_{n} \circ h'$ is a single morphism in $\mathsf{C}$, the induction is complete.
\end{proof}

\begin{proposition}
A functor $F\colon  \mathsf{C} \mor \mathsf{D}$ is full and surjective on objects if, and only if, it is a regular epimorphism and Riguet-full ($\equiv$ $\mathrm{K}^{\mathrm{r}}(F)$ is a Riguet congruence on $\mathsf{C}$). 
\end{proposition}

\begin{proof}
$(\Rightarrow)$ Suppose $F$ is full and surjective on objects. Then, by Proposition~\ref{FSOReg}, $F$ is a regular epimorphism. For Riguet-fullness, let $a$, $a'$, $b$, $b'\in \mathrm{Ob}(\mathsf{C})$ and $f\in \mathrm{Hom}_{\mathsf{C}}(a,b)$, such that $F(a) = F(a')$ and $F(b) = F(b')$. Then $F(f)\in 
\mathrm{Hom}_{\mathsf{D}}(F(a'),F(b'))$. Hence there exists an $f'\in \mathrm{Hom}_{\mathsf{C}}(a',b')$ such that $F(f) = F(f')$.

$(\Leftarrow)$ Suppose $F$ is the coequalizer of its kernel pair $(\mathsf{KP}(F),(P_{0},P_{1}))$ and is Riguet-full. The coequalizer of any pair of functors is surjective on objects. Thus, $F$ is surjective on objects. Fullness follows from Proposition~\ref{red}.
\end{proof}

\begin{corollary}
Let $\Phi$ be a Riguet congruence on a category $\mathsf{C}$. Then the projection functor $P_{\Phi}\colon \mathsf{C} \mor \mathsf{C}/\Phi$ is full and surjective on objects.
\end{corollary}

\subsection{Generalized congruences}\label{S:gcgr}
We next investigate the relationship between Riguet congruences and generalized congruences, as defined by Bruckner \cite{Br15} and Haucourt \cite{H06}. Generalized congruences were originally introduced in 1999 by Bednarczyk \textit{et al.} \cite{BBP99}; however, throughout this paper we adopt the formulation given in \cite{Br15} and \cite{H06}. The reason for this choice is that, for a category $\mathsf{C}$ and a generalized congruence $\left(\Phi^{\mathrm{ob}}, \Phi^{\mathrm{fl}}\right)$ on $\mathsf{C}$ in the sense of Bednarczyk \textit{et al.}, the domain of definition of the partial equivalence relation $\Phi^{\mathrm{fl}}$ is the set of all $\Phi^{\mathrm{ob}}$-composable sequences in $\mathsf{C}$ (see below for this notion). Thus, $\Phi^{\mathrm{fl}}$ restricted to it is an equivalence relation. 

We begin by introducing the notion of a generalized congruence on a category and by establishing its fundamental properties from a lattice-theoretic perspective. As a preliminary step, which will be essential for defining generalized congruences, we introduce, for a category $\mathsf{C}$ and an equivalence relation 
$\Phi^{\mathrm{ob}}$ on $\mathrm{Ob}(\mathsf{C})$, the set of all $\Phi^{\mathrm{ob}}$-composable sequences in 
$\mathsf{C}$.

\begin{definition}\label{DCongGen}
Let $\mathsf{C}$ be a category, $\Phi^{\mathrm{ob}}$ an equivalence relation on $\mathrm{Ob}(\mathsf{C})$ and 
$\mathrm{Mor}(\mathsf{C})^{+}$ the underlying set of $\mathbf{Mor}(\mathsf{C})^{+}$, the free semigroup on 
$\mathrm{Mor}(\mathsf{C})$. A $\Phi^{\mathrm{ob}}$-\emph{composable} sequence in $\mathsf{C}$ is a finite nonempty sequence $(f_{i})_{i\in m}$ of morphisms in $\mathsf{C}$, such that, for every $i\in m-1$, 
$(\mathrm{tg}(f_{i}),\mathrm{sc}(f_{i+1}))\in \Phi^{\mathrm{ob}}$. We denote by $\mathrm{Mor}(\mathsf{C})^{+}_{\Phi^{\mathrm{ob}}}$ the subset of $\mathrm{Mor}(\mathsf{C})^{+}$ consisting of all $\Phi^{\mathrm{ob}}$-composable sequences in $\mathsf{C}$. For every morphism $f$ of $\mathrm{Mor}(\mathsf{C})$, we will let $(f)$ stand for the corresponding  word of length one in $\mathrm{Mor}(\mathsf{C})^{+}$. Moreover, we denote by  
$\curlywedge$ the concatenation operation on the underlying set of $\mathbf{Mor}(\mathsf{C})^{+}$.
\end{definition}

\begin{remark}
Let $\mathsf{C}$ be a category. Then, for $\nabla_{\mathrm{Ob}(\mathsf{C})}$, the greatest equivalence relation on $\mathrm{Ob}(\mathsf{C})$, we have that   
$\mathrm{Mor}(\mathsf{C})^{+}_{\nabla_{\mathrm{Ob}(\mathsf{C})}} = \mathrm{Mor}(\mathsf{C})^{+}$. For 
$\Delta_{\mathrm{Ob}(\mathsf{C})}$, the smallest equivalence relation on $\mathrm{Ob}(\mathsf{C})$, we have that $\mathrm{Mor}(\mathsf{C})^{+}_{\Delta_{\mathrm{Ob}(\mathsf{C})}}$ is the set of all strictly composable sequences of morphisms in $\mathsf{C}$, i.e., the set of all finite nonempty sequences 
$(f_{i})_{i\in m}$ of morphisms in $\mathsf{C}$, such that, for every $i\in m-1$, 
$\mathrm{tg}(f_{i})=\mathrm{sc}(f_{i+1})$.
\end{remark}

\begin{definition}\label{DGCgr}
A \emph{generalized congruence} on a category $\mathsf{C}$ is an ordered pair
$\Psi=\left(\Psi^{\mathrm{ob}},\Psi^{\mathrm{fl}}\right)$ in which
\begin{enumerate}
\item
$\Psi^{\mathrm{ob}}$ is an equivalence relation on $\mathrm{Ob}(\mathsf{C})$ and

\item
$\Psi^{\mathrm{fl}}$ is an equivalence relation on $\mathrm{Mor}(\mathsf{C})^{+}_{\Psi^{\mathrm{ob}}}$ such that the following conditions are satisfied
\begin{enumerate}
\item
for every $a$, $b\in \mathrm{Ob}(\mathsf{C})$, if $(a,b) \in \Psi^{\mathrm{ob}}$, then $((\mathrm{id}_{a}),(\mathrm{id}_{b}))\in \Psi^{\mathrm{fl}}$;
\item
for every $(f_{i})_{i \in m}, (g_{j})_{j \in n} \in \mathrm{Mor}(\mathsf{C})^{+}_{\Psi^{\mathrm{ob}}}$, if $((f_{i})_{i \in m}, (g_{j})_{j\in n}) \in \Psi^{\mathrm{fl}}$, then $(\mathrm{sc}(f_{0}), \mathrm{sc}(g_{0})) \in \Psi^{\mathrm{ob}}$ and $(\mathrm{tg}(f_{m-1}), \mathrm{tg}(g_{n-1})) \in \Psi^{\mathrm{ob}}$;

\item
for every $f, g \in \mathrm{Mor}(\mathsf{C})$, if $\mathrm{sc}(g) = \mathrm{tg}(f)$, then $((f,g), (g \circ f)) \in \Psi^{\mathrm{fl}}$;

\item
for every $(f_{i})_{i \in m}, (f'_{i'})_{i' \in m'}, (g_{j})_{j \in n}, (g'_{j'})_{j' \in n'} \in \mathrm{Mor}(\mathsf{C})^{+}_{\Psi^{\mathrm{ob}}}$, \[
\begin{aligned}
\text{if }\;
&\left\{
\begin{array}{c}
((f_{i})_{i \in m}, (f'_{i'})_{i' \in m'}),
((g_{j})_{j \in n}, (g'_{j'})_{j' \in n'})
\in \Psi^{\mathrm{fl}}
\\[2pt]
\text{and}
\\[2pt]
(\mathrm{tg}(f_{m-1}), \mathrm{sc}(g_{0}))
\in \Psi^{\mathrm{ob}}
\end{array}
\right.,
\\[4pt]
\text{then }\;
&((f_{i})_{i \in m} \curlywedge (g_{j})_{j \in n},
(f'_{i'})_{i' \in m'} \curlywedge (g'_{j'})_{j' \in n'})
\in \Psi^{\mathrm{fl}}.
\end{aligned}
\]
\end{enumerate}
\end{enumerate}
We denote by $\mathrm{GCgr}(\mathsf{C})$ the set of all generalized congruences on $\mathsf{C}$.
\end{definition}

\begin{remark}
Let $\Psi=\left(\Psi^{\mathrm{ob}},\Psi^{\mathrm{fl}}\right)$ be a generalized congruence on $\mathsf{C}$. Then the objects of $\mathsf{C}/\Psi$, the quotient category of $\mathsf{C}$ by $\Psi$, are the equivalence classes of objects of $\mathsf{C}$ with respect to $\Psi^{\mathrm{ob}}$, whereas the morphisms are given by equivalence classes of $\Psi^{\mathrm{ob}}$-composable sequences in $\mathsf{C}$ with respect to $\Psi^{\mathrm{fl}}$. For details, including the universal property of the quotient, see~\cite{BBP99} or \cite{Br15}.
\end{remark}

In the following proposition we state that the set of all generalized congruences on a category, ordered by inclusion, is an algebraic lattice. The proof is immediate from the definition of a generalized congruence and is omitted.

\begin{proposition}\label{GCgrAlgLatt}
Let $\mathsf{C}$ be a category. Then 
\begin{enumerate}
\item $\nabla_{\mathsf{C}} = \left(\nabla_{\mathrm{Ob}(\mathsf{C})},\nabla_{\mathrm{Mor}(\mathsf{C})^{+}}\right)\in\mathrm{GCgr}(\mathsf{C})$ and it is the greatest generalized congruence on $\mathsf{C}$.
\item If $\mathcal{F}$ is a nonempty subset of $\mathrm{GCgr}(\mathsf{C})$, then $\bigcap \mathcal{F}\in\mathrm{GCgr}(\mathsf{C})$.
\item If $\mathcal{F}$ is a nonempty directed subset of $\mathrm{GCgr}(\mathsf{C})$, then $\bigcup \mathcal{F}\in\mathrm{GCgr}(\mathsf{C})$.
\end{enumerate}
Therefore, $\mathrm{GCgr}(\mathsf{C})$, ordered by inclusion, is an algebraic lattice. We denote by 
$\mathrm{GCg_{\mathsf{C}}}$ the corresponding algebraic closure operator that sends each relation $\Phi = (\Phi^{\mathrm{ob}},\Phi^{\mathrm{fl}})$ on $(\mathrm{Ob}(\mathsf{C}),\mathrm{Mor}(\mathsf{C})^{+})$ to 
$\mathrm{GCg_{\mathsf{C}}}(\Phi)$, the generalized congruence on $\mathsf{C}$  generated by $\Phi$, which is 
$\bigcap\{\Psi\in \mathrm{GCgr}(\mathsf{C})\mid \Phi\subseteq \Psi\}$.
\end{proposition}

\begin{remark}
For the generalized congruences on a category $\mathsf{C}$, the set operations $\bigcup$, $\bigcap$ and the inclusion relation are defined coordinatewise.
\end{remark}

\begin{remark}
The smallest generalized congruence on $\mathsf{C}$ is
\[
\Bigl(\Delta_{\mathrm{Ob}(\mathsf{C})},\Lambda_{\mathrm{Mor}(\mathsf{C})^{+}_{\Delta_{\mathrm{Ob}(\mathsf{C})}}}\Bigr),
\]
where 
$\Lambda_{\mathrm{Mor}(\mathsf{C})^{+}_{\Delta_{\mathrm{Ob}(\mathsf{C})}}}$ is the smallest equivalence relation on $\mathrm{Mor}(\mathsf{C})^{+}_{\Delta_{\mathrm{Ob}(\mathsf{C})}}$ generated by 
$\left\{ ((f,g), (g \circ f)) \in (\mathrm{Mor}(\mathsf{C})^{+}_{\Phi^{\mathrm{ob}}})^{2} \Bigr|\, \mathrm{sc}(g) = \mathrm{tg}(f)\right \}$  
and satisfying conditions $(\mathrm{a})$ and $(\mathrm{d})$ in Definition~\ref{DGCgr}.
Note that $\Lambda_{\mathrm{Mor}(\mathsf{C})^{+}_{\Delta_{\mathrm{Ob}(\mathsf{C})}}}$ is generally not $\Delta_{\mathrm{Mor}(\mathsf{C})^{+}_{\Delta_{\mathrm{Ob}(\mathsf{C})}}}$. Moreover, 
$
\textstyle
\Bigl(\Delta_{\mathrm{Ob}(\mathsf{C})},\Lambda_{\mathrm{Mor}(\mathsf{C})^{+}_{\Delta_{\mathrm{Ob}(\mathsf{C})}}}\Bigr) = \bigcap \mathrm{GCgr}(\mathsf{C}).
$ 

\end{remark} 

Our next aim is to elucidate, from a lattice-theoretic perspective, the relationship between Riguet congruences and generalized congruences on a category. As a preliminary step toward this goal, we introduce the notion of a strong generalized congruence, which will play a key role in achieving it, and clarify its relationship with regular generalized congruences and with the generalized kernels of functors---defined by Bednarczyk \textit{et al.} in~\cite{BBP99}---that are full and surjective on objects.

\begin{definition}
Let $\mathsf{C}$ be a category and $\Psi$ a generalized congruence on $\mathsf{C}$. We will say that $\Psi$ is \emph{strong} if, for every $(a,a'), (b,b') \in \Psi^{\mathrm{ob}}$ and every $f\in \mathrm{Mor}(\mathsf{C})$, if $f \in \mathrm{Hom}_{\mathsf{C}}(a,b)$, then there exists an $f' \in \mathrm{Hom}_{\mathsf{C}}(a',b')$ such that $((f),(f'))\in \Psi^{\mathrm{fl}}$. We denote by $\mathrm{SGCgr}(\mathsf{C})$ the set of all strong generalized congruences on $\mathsf{C}$ and by $(\mathrm{SGCgr}(\mathsf{C}),\subseteq)$ the corresponding ordered set.
\end{definition}

\begin{definition}\label{resregker}
Let $\Psi=\left(\Psi^{\mathrm{ob}},\Psi^{\mathrm{fl}}\right)$ be a generalized congruence on a category $\mathsf{C}$. Then the \emph{restriction of} $\Psi$ \emph{to singleton sequences of morphisms}, denoted by $\bb{\Psi}$, is the ordered pair $\left(\bb{\Psi}^{\mathrm{ob}},\bb{\Psi}^{\mathrm{fl}}\right)$ in which $\bb{\Psi}^{\mathrm{ob}} = \Psi^{\mathrm{ob}}$ and $\bb{\Psi}^{\mathrm{fl}}$ is the subrelation of $\Psi^{\mathrm{fl}}$ consisting of all pairs of singleton sequences, i.e., $((f_i)_{i\in m}, (g_j)_{j\in n}) \in \bb{\Psi}^{\mathrm{fl}}$ if, and only if, $((f_i)_{i\in m}, (g_j)_{j\in n}) \in \Psi^{\mathrm{fl}}$ and  $m=n=1$. The \emph{regular part} of $\Psi$, denoted by $\mathrm{Rg}(\Psi)$, is $\mathrm{GCg_{\mathsf{C}}}(\bb{\Psi})$, the generalized congruence on $\mathsf{C}$ generated by $\bb{\Phi}$. Note that $\mathrm{Rg}(\Psi)\subseteq \Psi$. We will say that $\Psi$ is \emph{regular} if $\mathrm{Rg}(\Psi) = \Psi$. We denote by $\mathrm{RgGCgr}(\mathsf{C})$ the set of all regular generalized congruences on $\mathsf{C}$  and by $(\mathrm{RgGCgr}(\mathsf{C}),\subseteq)$ the corresponding ordered set.

Let $F \colon \mathsf{C} \mor \mathsf{D}$ be a functor. The \emph{generalized kernel} of $F$ or, simply, the \emph{kernel} of $F$, denoted by $\mathrm{K}(F)$, is the generalized congruence on $\mathsf{C}$ whose object part $\mathrm{K}(F)^{\mathrm{ob}}$ is defined as: for every $a$, $b\in \mathrm{Ob}(\mathsf{C})$, $(a,b)\in \mathrm{K}(F)^{\mathrm{ob}}$ if, and only if, $F(a)= F(b)$; and whose morphism part $\mathrm{K}(F)^{\mathrm{fl}}$ is defined as: for every $(f_{i})_{i \in m}, (g_{j})_{j \in n} \in \mathrm{Mor}(\mathsf{C})^{+}_{\Psi^{\mathrm{ob}}}$, $((f_{i})_{i \in m}, (g_{j})_{j \in n}) \in \mathrm{K}(F)^{\mathrm{fl}}$ if, and only if, both  
$F(f_{m-1})\circ\cdots\circ F(f_{0})$ and $F(g_{n-1})\circ\cdots\circ F(g_{0})$ are defined in $\mathsf{D}$ and $F(f_{m-1})\circ\cdots\circ F(f_{0}) = F(g_{n-1})\circ\cdots\circ F(g_{0})$.  Note that $\mathrm{Rg}(\mathrm{K}(F))$, the regular part of $\mathrm{K}(F)$, is the generalized congruence on $\mathsf{C}$ generated by the kernel of the object mapping of $F$, i.e., by the equivalence relation $\{(a,b)\in \mathrm{Ob}(\mathsf{C})\mid F(a)= F(b)\}$ on $\mathrm{Ob}(\mathsf{C})$, together with the binary relation on $\mathrm{Hom}(1,\mathrm{Mor}(\mathsf{C}))\subseteq\mathrm{Mor}(\mathsf{C})^{+}$ consisting of all those ordered pairs $((f),(g))\in\mathrm{Hom}(1,\mathrm{Mor}(\mathsf{C}))^{2}$ such that $F(f) = F(g)$.
\end{definition}

In the following proposition, for a fixed category $\mathsf{C}$, we show that forming the regular part of a generalized congruence induces an algebraic interior operator, and we describe its adjoint relationship with the canonical embedding of the category of regular generalized congruences into the category of generalized congruences.

\begin{proposition}
The operator $\mathrm{Rg}\colon \mathrm{GCgr}(\mathsf{C}) \mor \mathrm{GCgr}(\mathsf{C})$ that sends a generalized congruence $\Psi$ on $\mathsf{C}$ to $\mathrm{Rg}(\Psi)$ is an algebraic interior operator. Note that $\mathrm{Rg}(\Psi)$ is the greatest regular generalized congruence on $\mathsf{C}$ contained in $\Psi$. Moreover, the canonical embedding of $\mathsf{RgGCgr}(\mathsf{C})$, the category associated to the ordered set $(\mathrm{RgGCgr}(\mathsf{C}),\subseteq)$ into $\mathsf{GCgr}(\mathsf{C})$, the category associated to the ordered set $(\mathrm{GCgr}(\mathsf{C}),\subseteq)$, has a right adjoint, precisely the functor $\mathrm{Rg}$ that sends $\Psi$ to $\mathrm{Rg}(\Psi)$.
\end{proposition}

\begin{proposition}\label{KerRg}
For every $\Psi \in \mathrm{GCgr}(\mathsf{C})$, $\mathrm{Rg}(\mathrm{K}(P_{\Psi})) = \mathrm{Rg}(\Psi)$, where $P_{\Psi}$ is the canonical projection from $\mathsf{C}$ to $\mathsf{C}/\Psi$.
\end{proposition}

\begin{proof}
For $P_{\Psi}$, by Definition~\ref{resregker}, we have that $\mathrm{Rg}(\mathrm{K}(P_{\Psi}))$ is the generalized congruence on $\mathsf{C}$ generated by $\{(a,b)\in \mathrm{Ob}(\mathsf{C})^{2}\mid [a]_{\Psi^{\mathrm{ob}}} = [b]_{\Psi^{\mathrm{ob}}}\} = \Psi^{\mathrm{ob}}$, together with the binary relation on $\mathrm{Hom}(1,\mathrm{Mor}(\mathsf{C}))\subseteq\mathrm{Mor}(\mathsf{C})^{+}$ consisting of all those ordered pairs $((f),(g))\in\mathrm{Hom}(1,\mathrm{Mor}(\mathsf{C}))^{2}$ such that $[f]_{\Psi^{\mathrm{fl}}} = [g]_{\Psi^{\mathrm{fl}}}$ or, what is equivalent, such that $((f),(g)) \in \Psi^{\mathrm{fl}}$, for all morphisms $f,g \in \mathrm{Mor}(\mathsf{C})$. Therefore $\mathrm{Rg}(\mathrm{K}(P_{\Psi})) = \mathrm{Rg}(\Psi)$.
\end{proof}

We next characterize functors that are full and surjective on objects in terms of extremal epimorphisms and generalized kernels that are both strong and regular generalized congruences.

\begin{proposition}
A functor $F \colon \mathsf{C} \mor \mathsf{D}$ is full and surjective on objects if, and only if, $F$ is an extremal epimorphism and 
$\mathrm{K}(F)\in \mathrm{RgGCgr}\cap \mathrm{SGCgr}(\mathsf{C})$.
\end{proposition}

\begin{proof}
$(\Rightarrow)$ Suppose $F$ is full and surjective on objects. Then, by Proposition~\ref{FSOReg}, $F$ is a regular epimorphism. Hence, by 
Bednarczyk \textit{et al.}~\cite{BBP99}, $F$ is extremal and $\mathrm{K}(F)$ is regular. 
Now we prove that $\mathrm{K}(F)$ is strong. Given $(a, a'), (b, b') \in \mathrm{K}(F)^{\mathrm{ob}}$, which means that $F(a) = F(a')$ and $F(b) = F(b')$, and $f \colon a \mor b$, then we have that  $F(f) \colon F(a') \mor F(b')$ is a morphism in $\mathsf{D}$. Hence, by fullness, there exists $f' \colon a' \mor b'$ such that $F(f') = F(f)$. Thus $((f), (f')) \in \mathrm{K}(F)^{\mathrm{fl}}$, satisfying the strong condition. 

$(\Leftarrow)$ Suppose $F$ is an extremal epimorphism and $\mathrm{K}(F)\in \mathrm{RgGCgr}\cap \mathrm{SGCgr}(\mathsf{C})$. Then the surjectivity on objects is immediate from $F$ being an extremal epimorphism. We next prove that $F$ is full. Let $a, b \in \mathrm{Ob}(\mathsf{C})$ and $u \colon F(a) \mor F(b)$ in $\mathsf{D}$. Since $F$ is an extremal epimorphism, there exists a family of morphisms $(f_{i})_{i \in n}$ in $\mathrm{Mor}(\mathsf{C})$ such that $u = F(f_{n-1}) \circ \dots \circ F(f_{0})$ and $(\mathrm{tg}(f_{i}), \mathrm{sc}(f_{i+1})) \in \mathrm{K}(F)^{\mathrm{ob}}$. We also have that $(a, \mathrm{sc}(f_{0})) \in \mathrm{K}(F)^{\mathrm{ob}}$ and $(b, \mathrm{tg}(f_{n-1})) \in \mathrm{K}(F)^{\mathrm{ob}}$. We now use the fact that $\mathrm{K}(F)$ is strong iteratively:
\begin{enumerate}
\item Since $(a, \mathrm{sc}(f_{0})) \in \mathrm{K}(F)^{\mathrm{ob}}$ and $f_{0} \colon \mathrm{sc}(f_{0}) \mor \mathrm{tg}(f_{0})$, there is a $g_{0} \colon a \mor b'_{0}$ such that $((f_{0}), (g_{0})) \in \mathrm{K}(F)^{\mathrm{fl}}$ and $(\mathrm{tg}(f_{0}), b'_{0}) \in \mathrm{K}(F)^{\mathrm{ob}}$.
\item For $f_{1} \colon \mathrm{sc}(f_{1}) \mor \mathrm{tg}(f_{1})$, we have $(\mathrm{tg}(f_{0}), \mathrm{sc}(f_{1})) \in \mathrm{K}(F)^{\mathrm{ob}}$ (by sequence composability) and $(\mathrm{tg}(f_{0}), b'_{0}) \in \mathrm{K}(F)^{\mathrm{ob}}$, so $(b'_{0}, \mathrm{sc}(f_{1})) \in \mathrm{K}(F)^{\mathrm{ob}}$. By strongness, there is a  $g_{1} \colon b'_{0} \mor b'_{1}$ such that $((f_{1}), (g_{1})) \in \mathrm{K}(F)^{\mathrm{fl}}$.
\item By repeating this process, we obtain a sequence $(g_{i})_{i \in n}$ such that $\mathrm{sc}(g_{0}) = a$, $\mathrm{tg}(g_{n-1}) = b$, and for every $i\in n$, $\mathrm{tg}(g_{i}) = \mathrm{sc}(g_{i+1})$. 
\end{enumerate}
This sequence $(g_{i})_{i \in n}$ is, by construction, $\mathrm{K}(F)^{\mathrm{ob}}$-composable in $\mathsf{C}$. Let $f = g_{n-1} \circ \dots \circ g_{0}$ be its composition in $\mathsf{C}$. By successive applications of Axiom 2(c) of Definition~\ref{DGCgr} and the transitivity of $\mathrm{K}(F)^{\mathrm{fl}}$, we have that $((g_{i})_{i \in n}, (f)) \in \mathrm{K}(F)^{\mathrm{fl}}$. Furthermore, since $((f_{i}), (g_{i})) \in \mathrm{K}(F)^{\mathrm{fl}}$ for each $i$ (by the strong property), Axiom 2(d) allows us to concatenate these pairs, yielding:
$((f_{i})_{i \in n}, (g_{i})_{i \in n}) \in \mathrm{K}(F)^{\mathrm{fl}}$. 
By the transitivity of $\mathrm{K}(F)^{\mathrm{fl}}$, we conclude that $((f), (f_{i})_{i \in n}) \in \mathrm{K}(F)^{\mathrm{fl}}$. Finally, by the definition of the generalized kernel we have that 
$F(f) = F(f_{n-1}) \circ \dots \circ F(f_{0}) = u$. The regularity of $\mathrm{K}(F)$ ensures that the generalized congruence is indeed generated by these singleton relations $\bb{\mathrm{K}(F)}$.
Thus, $F$ is full.
\end{proof}

We next prove one of the main results of this section, namely, the existence of a Scott continuous morphism from $(\mathrm{RCgr}(\mathsf{C}),\subseteq)$ to $(\mathrm{GCgr}(\mathsf{C}),\subseteq)$. As we will see in Theorem~\ref{AdjSGR}, this morphism induces a functor from the category $\mathsf{RCgr}(\mathsf{C})$ associated with $(\mathrm{RCgr}(\mathsf{C}),\subseteq)$ to the category $\mathsf{SGCgr}(\mathsf{C})$ associated with $(\mathrm{SGCgr}(\mathsf{C}),\subseteq)$, and this functor has a right adjoint.

\begin{theorem}\label{RCgrGCgr}
Let $\mathsf{C}$ be a category. Then there exists a Scott continuous morphism 
$(\bigcdot)^{\natural}$ from $(\mathrm{RCgr}(\mathsf{C}),\subseteq)$ to $(\mathrm{GCgr}(\mathsf{C}),\subseteq)$ such that, for every Riguet congruence $\Phi$, denoting by $\Phi^{\natural}$ the value of $(\bigcdot)^{\natural}$ at $\Phi$, the following properties hold:
\begin{enumerate}
\item $\Phi^{\natural\mathrm{ob}}$, the first coordinate of 
$\Phi^{\natural}$, is $\Phi^{\mathrm{ob}}$.
\item $\Phi^{\natural\mathrm{fl}}$, the second coordinate of $\Phi^{\natural}$, is the smallest equivalence relation on $\mathrm{Mor}(\mathsf{C})^{+}_{\Phi^{\mathrm{ob}}}$ generated by $\bigcup \Phi^{\mathrm{fl}}$ and satisfying the conditions (a)--(d).
\item $\Phi^{\natural}$ is a strong generalized congruence on $\mathsf{C}$.
\end{enumerate}
\end{theorem}

\begin{proof}
Let $\Phi$ be a Riguet congruence on $\mathsf{C}$. Then $\bigcup_{\scalebox{0.7}{$\left(\begin{smallmatrix}a&b\\a'& b'\end{smallmatrix}\right)$}\in \Phi^{\mathrm{ob}}\times \Phi^{\mathrm{ob}}}\Phi^{\mathrm{fl}}_{\scalebox{0.7}{$\left(\begin{smallmatrix}a&b\\a'& b'\end{smallmatrix}\right)$}}$, which is $\bigcup \Phi^{\mathrm{fl}}$, is a subset of $\mathrm{Mor}(\mathsf{C})\times \mathrm{Mor}(\mathsf{C})$. Henceforth, since there exists a canonical embedding $\eta_{\mathrm{Mor}(\mathsf{C})}$ from $\mathrm{Mor}(\mathsf{C})$ into $\mathrm{Mor}(\mathsf{C})^{+}$, that sends $f$ in $\mathrm{Mor}(\mathsf{C})$ to $(f)$ in $\mathrm{Mor}(\mathsf{C})^{+}$, we have that $\bigcup \Phi^{\mathrm{fl}}$ is isomorphic to $\eta_{\mathrm{Mor}(\mathsf{C})}^{2}[\bigcup \Phi^{\mathrm{fl}}]$, where 
\[
\textstyle
\eta_{\mathrm{Mor}(\mathsf{C})}^{2}[\bigcup \Phi^{\mathrm{fl}}] = 
\bigcup_{\scalebox{0.7}{$\left(\begin{smallmatrix}a&b\\a'& b'\end{smallmatrix}\right)$} \in \Phi^{\mathrm{ob}}\times \Phi^{\mathrm{ob}}}\left\{((f),(f'))\,\Bigr|\, (f,f')\in \Phi^{\mathrm{fl}}_{\scalebox{0.7}{$\left(\begin{smallmatrix}a&b\\a'& b'\end{smallmatrix}\right)$}}\right\},
\]
which is, by logical vacuity, a binary relation on $\mathrm{Mor}(\mathsf{C})^{+}_{\Phi^{\mathrm{ob}}}$.
From now on, by a harmless abuse of notation and for the sake of brevity, we shall write $\bigcup \Phi^{\mathrm{fl}}$ in place of $\eta_{\mathrm{Mor}(\mathsf{C})}^{2}[\bigcup \Phi^{\mathrm{fl}}]$. 
Then, by Proposition~\ref{GCgrAlgLatt}, there exists 
$\mathrm{GCg}_{\mathsf{C}}\left(\Phi^{\mathrm{ob}},\bigcup \Phi^{\mathrm{fl}}\right)$, the generalized congruence on $\mathsf{C}$ generated by $\left(\Phi^{\mathrm{ob}},\bigcup \Phi^{\mathrm{fl}}\right)$, which   
is the smallest generalized congruence on $\mathsf{C}$ that contains $\left(\Phi^{\mathrm{ob}},\bigcup \Phi^{\mathrm{fl}}\right)$. Bednarczyk \textit{et al.} in~\cite{BBP99} would call it the principal congruence generated by $\left(\Phi^{\mathrm{ob}},\bigcup \Phi^{\mathrm{fl}}\right)$. We denote by $\Phi^{\natural\mathrm{ob}}$ and $\Phi^{\natural\mathrm{fl}}$ the first and second coordinates, respectively, of $\mathrm{GCg}_{\mathsf{C}}\left(\Phi^{\mathrm{ob}},\bigcup \Phi^{\mathrm{fl}}\right)$. 

Then we let $(\bigcdot)^{\natural}$ stand for the mapping from $\mathrm{RCgr}(\mathsf{C})$ to $\mathrm{GCgr}(\mathsf{C})$ that sends a Riguet congruence $\Phi = \left(\Phi^{\mathrm{ob}},\Phi^{\mathrm{fl}}\right)$ on $\mathsf{C}$ to the generalized congruence 
$\Phi^{\natural} = \mathrm{GCg}_{\mathsf{C}}\left(\Phi^{\mathrm{ob}},\bigcup \Phi^{\mathrm{fl}}\right)$ on $\mathsf{C}$. We leave it to the reader to verify that the mapping $(\bigcdot)^{\natural}$ is 
Scott continuous, i.e., that, for every nonempty directed subset $\mathcal{F}$ of $\mathrm{RCgr}(\mathsf{C})$,  
$\left(\bigcup_{\Phi\in\mathcal{F}}\Phi\right)^{\natural} = \bigcup_{\Phi\in\mathcal{F}}\Phi^{\natural}$.

We next prove that $\Phi^{\natural\mathrm{ob}} = \Phi^{\mathrm{ob}}$ by explicitly exhibiting an upper bound for the morphism $(\bigcdot)^{\natural}$. Let $(\bigcdot)^{\dagger}$ be the mapping from $\mathrm{RCgr}(\mathsf{C})$ to $\mathrm{GCgr}(\mathsf{C})$ that sends
a Riguet congruence $\Phi$ on $\mathsf{C}$ to $\Phi^{\dagger} = \left(\Phi^{\mathrm{ob}},\Phi^{\dagger\mathrm{fl}}\right)$, where $\Phi^{\dagger\mathrm{fl}}$ is the binary relation on $\mathrm{Mor}(\mathsf{C})^{+}_{\Phi^{\mathrm{ob}}}$ defined, for every $((f_{i})_{i\in m},(g_{j})_{j\in n})$ in 
$\mathrm{Mor}(\mathsf{C})^{+}_{\Phi^{\mathrm{ob}}}\times \mathrm{Mor}(\mathsf{C})^{+}_{\Phi^{\mathrm{ob}}}$, as: $((f_{i})_{i\in m},(g_{j})_{j\in n})\in \Phi^{\dagger\mathrm{fl}}$ if, and only if, 
$[f_{m-1}]_{\Phi^{\mathrm{fl}}}\circ\cdots\circ [f_{0}]_{\Phi^{\mathrm{fl}}} = [g_{n-1}]_{\Phi^{\mathrm{fl}}}\circ\cdots\circ [g_{0}]_{\Phi^{\mathrm{fl}}}$. Then, by definition of $\Phi^{\dagger}$, we have that $\Phi^{\mathrm{ob}}$ is an equivalence relation on $\mathrm{Ob}(\mathsf{C})$ and that $\Phi^{\dagger\mathrm{fl}}$ is an equivalence relation on $\mathrm{Mor}(\mathsf{C})^{+}_{\Phi^{\mathrm{ob}}}$. Let us verify that $\Phi^{\dagger}$ satisfies the four conditions stated in Definition~\ref{DGCgr}.

$(\mathrm{a})$
Let $a, b \in \mathrm{Ob}(\mathsf{C})$ be such that $(a,b) \in \Phi^{\mathrm{ob}}$. 
Thus $[\mathrm{id}_{a}]_{\Phi^{\mathrm{fl}}} = [\mathrm{id}_{b}]_{\Phi^{\mathrm{fl}}}$ and so $((\mathrm{id}_{a}), (\mathrm{id}_{b})) \in \Phi^{\dagger\mathrm{fl}}$.

$(\mathrm{b})$
Let $(f_{i})_{i \in m}, (g_{j})_{j \in n} \in \mathrm{Mor}(\mathsf{C})^{+}$ be such that $((f_{i})_{i \in m}, (g_{j})_{j \in n}) \in \Phi^{\dagger\mathrm{fl}}$. 
Then we have that $[g_{n-1}]_{\Phi^{\mathrm{fl}}} \circ \cdots \circ [g_0]_{\Phi^{\mathrm{fl}}} = [f_{m-1}]_{\Phi^{\mathrm{fl}}} \circ \cdots \circ [f_0]_{\Phi^{\mathrm{fl}}}$ in $\mathsf{C}/\Phi$. 
Therefore, $(\mathrm{sc}(f_0), \mathrm{sc}(g_0)) \in \Phi^{\mathrm{ob}}$ and $(\mathrm{tg}(f_{m-1}), \mathrm{tg}(g_{n-1})) \in \Phi^{\mathrm{ob}}$. 

$(\mathrm{c})$
Let $f, g \in \mathrm{Mor}(\mathsf{C})$ be  such that $\mathrm{tg}(f) = \mathrm{sc}(g)$. Then,  
by Definition~\ref{DRiguetQ}, we have that $[g]_{\Phi^{\mathrm{fl}}} \circ [f]_{\Phi^{\mathrm{fl}}} = [g \circ f]_{\Phi^{\mathrm{fl}}}$. 
Thus, $((f,g), (g\circ f)) \in \Phi^{\dagger\mathrm{fl}}$.

$(\mathrm{d})$
Let $(f_{i})_{i \in m}, (f'_{i'})_{i' \in m'}, (g_{j})_{j \in n}, (g'_{j'})_{j' \in n'} \in \mathrm{Mor}(\mathsf{C})^{+}$ be such that 
\[
((f_{i})_{i \in m}, (f'_{i'})_{i' \in m'}), ((g_{j})_{j \in n}, (g'_{j'})_{j' \in n'}) \in \Phi^{\dagger\mathrm{fl}}
\]
and $(\mathrm{tg}(f_{m-1}), \mathrm{sc}(g_{n-1}) \in \Phi^{\mathrm{ob}}$. 
Since
\begin{align*}
[f_{m-1}]_{\Phi^{\mathrm{fl}}} \circ \cdots \circ [f_{0}]_{\Phi^{\mathrm{fl}}}
&=
[f'_{m'-1}]_{\Phi^{\mathrm{fl}}} \circ \cdots \circ [f'_{0}]_{\Phi^{\mathrm{fl}}}
\text{ and }
\\
[g_{n-1}]_{\Phi^{\mathrm{fl}}} \circ \cdots \circ [g_{0}]_{\Phi^{\mathrm{fl}}}
&=
[g'_{n'-1}]_{\Phi^{\mathrm{fl}}} \circ \cdots \circ [g'_{0}]_{\Phi^{\mathrm{fl}}},
\end{align*}
we have that 
\begin{multline*}
[g_{n-1}]_{\Phi^{\mathrm{fl}}} \circ \cdots \circ [g_{0}]_{\Phi^{\mathrm{fl}}}
\circ
[f_{m-1}]_{\Phi^{\mathrm{fl}}} \circ \cdots \circ [f_{0}]_{\Phi^{\mathrm{fl}}}
=\\
[g'_{n'-1}]_{\Phi^{\mathrm{fl}}} \circ \cdots \circ [g'_{0}]_{\Phi^{\mathrm{fl}}}
\circ
[f'_{m'-1}]_{\Phi^{\mathrm{fl}}} \circ \cdots \circ [f'_{0}]_{\Phi^{\mathrm{fl}}}.
\end{multline*}
Thus, by definition, $((f_{i})_{i \in m} \curlywedge (g_{j})_{ \in n}, (f'_{i'})_{i' \in m'} \curlywedge (g'_{j'})_{j' \in n'}) \in \Phi^{\dagger\mathrm{fl}}$.

We leave it to the reader to verify that the mapping $(\bigcdot)^{\dagger}$ is 
Scott continuous, i.e., that, for every nonempty directed subset $\mathcal{F}$ of $\mathrm{RCgr}(\mathsf{C})$,  
$\left(\bigcup_{\Phi\in\mathcal{F}}\Phi\right)^{\dagger} = \bigcup_{\Phi\in\mathcal{F}}\Phi^{\dagger}$. 

The morphism $(\bigcdot)^{\dagger}$ is an upper bound of $(\bigcdot)^{\natural}$, i.e., for every Riguet congruence $\Phi$ on $\mathsf{C}$, $\Phi^{\natural}\subseteq\Phi^{\dagger}$. This is so because 
$\left(\Phi^{\mathrm{ob}},\bigcup \Phi^{\mathrm{fl}}\right)\subseteq \Phi^{\dagger}$. Therefore 
$\Phi^{\natural\mathrm{ob}} = \Phi^{\mathrm{ob}}$.

We next prove that $\Phi^{\natural\mathrm{fl}}$ is the smallest equivalence relation on $\mathrm{Mor}(\mathsf{C})^{+}_{\Phi^{\mathrm{ob}}}$ generated by $\bigcup \Phi^{\mathrm{fl}}$ and satisfying the conditions $(\mathrm{a})$--$(\mathrm{d})$.
Since $\Phi^{\natural\mathrm{fl}}$ is the second coordinate of $\mathrm{GCg}_{\mathsf{C}}\left(\Phi^{\mathrm{ob}},\bigcup \Phi^{\mathrm{fl}}\right)$ and $\left(\Phi^{\mathrm{ob}},\bigcup \Phi^{\mathrm{fl}}\right)\subseteq \mathrm{GCg}_{\mathsf{C}}\left(\Phi^{\mathrm{ob}},\bigcup \Phi^{\mathrm{fl}}\right)$, it follows immediately that $\Phi^{\natural\mathrm{fl}}$ is an equivalence relation on $\mathrm{Mor}(\mathsf{C})^{+}_{\Phi^{\mathrm{ob}}}$ satisfying the conditions 
$(\mathrm{a})$--$(\mathrm{d})$ and containing $\bigcup \Phi^{\mathrm{fl}}$. Now, let $\Psi^{\mathrm{fl}}$ be another equivalence relation on $\mathrm{Mor}(\mathsf{C})^{+}_{\Phi^{\mathrm{ob}}}$ satisfying the conditions $(\mathrm{a})$--$(\mathrm{d})$ and containing $\bigcup \Phi^{\mathrm{fl}}$. Then $(\Phi^{\mathrm{ob}}, \Psi^{\mathrm{fl}}) \in \mathrm{GCgr}(\mathsf{C})$ and $(\Phi^{\mathrm{ob}}, \bigcup \Phi^{\mathrm{fl}}) \subseteq (\Phi^{\mathrm{ob}}, \Psi^{\mathrm{fl}})$. Hence $\Phi^{\natural}\subseteq (\Phi^{\mathrm{ob}}, \Psi^{\mathrm{fl}})$. Therefore, $\Phi^{\natural\mathrm{fl}}\subseteq \Psi^{\mathrm{fl}}$. From which it follows that $\Phi^{\natural\mathrm{fl}}$ is indeed the smallest equivalence relation on $\mathrm{Mor}(\mathsf{C})^{+}_{\Phi^{\mathrm{ob}}}$ generated by $\bigcup \Phi^{\mathrm{fl}}$ and satisfying the conditions $(\mathrm{a})$--$(\mathrm{d})$.

Verifying that $\Phi^{\natural}$ is a strong generalized congruence on $\mathsf{C}$ is left to the reader.
\end{proof}

\begin{remark}
For a detailed discussion of the relationship between congruences and generalized congruences, see~\cite{BBP99}, pp.~273--274. It can be summarized as follows: generalized congruences are a conservative extension of congruences.
\end{remark}

In the following proposition, we prove that the notion of Riguet congruence is strictly more restrictive than that of generalized congruence.

\begin{proposition}
There exists a category and a generalized congruence on it which is not in the image of $(\bigcdot)^{\natural}$. 
\end{proposition}

\begin{proof}
Let us consider $\mathsf{Dis}(2)$, the discrete category on the set $2 = \{0,1\}$. Then, for $\nabla_{\mathsf{Dis}(2)} = (\nabla_{\mathrm{Ob}(\mathsf{Dis}(2))}, \nabla_{\mathrm{Mor}(\mathsf{Dis}(2))^{+}})$, the greatest  generalized congruence on $\mathsf{Dis}(2)$, we have that $\mathsf{Dis}(2)/\nabla_{\mathsf{Dis}(2)}$ has only one object $[0]_{\nabla_{\mathrm{Ob}(\mathsf{Dis}(2))}}$ and only one arrow $[\mathrm{id}_{[0]_{\nabla_{\mathrm{Ob}(\mathsf{Dis}(2))}}}]_{\nabla_{\mathrm{Mor}(\mathsf{Dis}(2))^{+}}}$. However, there is no Riguet congruence $\Phi = (\Phi^{\mathrm{ob}}, \Phi^{\mathrm{fl}})$ on $\mathsf{Dis}(2)$ such that $\Phi^{\natural} = \nabla_{\mathsf{Dis}(2)}$. Otherwise, $(0,1) \in \Phi^{\mathrm{ob}}$. Thus, by condition $(\mathrm{f})$ in Definition~\ref{DRiguet}, there should exists a morphism $f \colon 0 \mor 1$ such that $(\mathrm{id}_{0}, f) \in \Phi^{\mathrm{fl}}_{\scalebox{.7}{$\left(\begin{smallmatrix}0&0\\0&1\end{smallmatrix}\right)$}}$ which is clearly not the case.
\end{proof}

In Theorem~\ref{RCgrGCgr}, for a Riguet congruence $\Phi$ on a category $\mathsf{C}$, we showed that  
$\Phi^{\natural\mathrm{fl}}$ is the smallest equivalence relation on $\mathrm{Mor}(\mathsf{C})^{+}_{\Phi^{\mathrm{ob}}}$ generated by $\bigcup \Phi^{\mathrm{fl}}$ and satisfying the conditions $(\mathrm{a})$--$(\mathrm{d})$. Our next objective is to give a constructive description of $\Phi^{\natural\mathrm{fl}}$, explicitly exhibiting the process by which it is obtained. This description will be especially useful in the proofs of Proposition~\ref{PNFuncMor}, Lemma~\ref{commutation} and Proposition~\ref{PComp2}.

\begin{definition}\label{CtvDef}
Let $\mathsf{C}$ be a category and $\Phi$ a Riguet congruence on $\mathsf{C}$. 
\begin{enumerate}
\item
We denote by $\mathrm{G}_{\Phi}$ the operator on $(\mathrm{Mor}(\mathsf{C})^{+}_{\Phi^{\mathrm{ob}}})^{2}$ that assigns to a relation $\Upsilon \subseteq (\mathrm{Mor}(\mathsf{C})^{+}_{\Phi^{\mathrm{ob}}})^{2}$, the relation
\[
\mathrm{G}_{\Phi}(\Upsilon) = (\Upsilon \circ \Upsilon) \cup \left( (\curlywedge \times \curlywedge) [\Upsilon \times \Upsilon] \cap (\mathrm{Mor}(\mathsf{C})^{+}_{\Phi^{\mathrm{ob}}})^{2} \right).
\]
\item
We define the family $(\Phi^{\mathrm{fl}}_{r})_{r \in \mathbb{N}}$ in $\mathrm{Sub}((\mathrm{Mor}(\mathsf{C})^{+}_{\Phi^{\mathrm{ob}}})^{2})$ recursively as follows:
\begin{align*}
\Phi^{\mathrm{fl}}_{0} 
&=
\textstyle
\left(\bigcup\Phi^{\mathrm{fl}}\right) \cup \mathrm{C} \cup \mathrm{C}^{-1} \cup \Delta_{\mathrm{Mor}(\mathsf{C})^{+}_{\Phi^{\mathrm{ob}}}},
\\
\Phi^{\mathrm{fl}}_{r+1} 
&=
\mathrm{G}_{\Phi}(\Phi^{\mathrm{fl}}_{r}), \text{ for } r \geq 0,
\end{align*}
where $\mathrm{C} = \left\{ ((f,g), (g \circ f)) \in (\mathrm{Mor}(\mathsf{C})^{+}_{\Phi^{\mathrm{ob}}})^{2} \Bigr|\, \mathrm{sc}(g) = \mathrm{tg}(f)\right \}$.
\item
We let $\Phi^{\mathrm{fl}}_{\omega}$ stand for $\bigcup_{r \in \mathbb{N}} \Phi^{\mathrm{fl}}_{r}$.
\end{enumerate}
\end{definition}

\begin{proposition}
\label{PGCongOpG}
Let $\mathsf{C}$ be a category and $\Phi$ a Riguet congruence on $\mathsf{C}$. Then, for every $r \in \mathbb{N}$, we have that
\begin{enumerate}
\item
$\Phi^{\mathrm{fl}}_{r}$ is a reflexive relation,
\item
$\Phi^{\mathrm{fl}}_{r} \subseteq \Phi^{\mathrm{fl}}_{r+1}$ \text{and}
\item
$\Phi^{\mathrm{fl}}_{r}$ is a symmetric relation.
\end{enumerate}
\end{proposition}

\begin{proof}
The proof of (1) is straightforward and follows immediately by induction on $r \in \mathbb{N}$, so we omit the details.

The proof of (2) follows readily from (1) and the definition of $ \Phi^{\mathrm{fl}}_{r+1}$, and the proof is omitted.  

We now prove the third statement by induction on $r \in \mathbb{N}$.

\textsf{Base step of the induction for the symmetry property.}

The inverse of $\Phi^{\mathrm{fl}}_{0}$ is 
$
\textstyle
\left( \Phi^{\mathrm{fl}}_{0} \right)^{-1}
=
\left(\bigcup\Phi^{\mathrm{fl}}\right)^{-1} \cup \mathrm{C}^{-1} \cup \mathrm{C} \cup \Delta_{\mathrm{Mor}(\mathsf{C})^{+}_{\Phi^{\mathrm{ob}}}}^{-1}.
$
But $\Delta_{\mathrm{Mor}(\mathsf{C})^{+}_{\Phi^{\mathrm{ob}}}}^{-1} = \Delta_{\mathrm{Mor}(\mathsf{C})^{+}_{\Phi^{\mathrm{ob}}}}$ and, since $\Phi$ is a Riguet congruence on $\mathsf{C}$ and  by the laws of the calculus of relations, we have that $\left(\bigcup\Phi^{\mathrm{fl}}\right)^{-1} = 
\bigcup \Phi^{\mathrm{fl}}$. Thus 
$
\textstyle
\left( \Phi^{\mathrm{fl}}_{0} \right)^{-1}
=
\bigcup\Phi^{\mathrm{fl}}.
$
Therefore, $\Phi^{\mathrm{fl}}_{0}$ is a symmetric relation.

\textsf{Inductive step of the induction for the symmetry property.}

Let us assume that the statement holds for $r \in \mathbb{N}$, i.e., that $\Phi^{\mathrm{fl}}_{r}$ is a symmetric relation, and let us prove it for $r+1$. Let $((f_{i})_{i \in m}, (g_{j})_{j \in n})$ be an element of $\Phi^{\mathrm{fl}}_{r+1}$. 
Then, by definition of $\Phi^{\mathrm{fl}}_{r+1}$, $((f_{i})_{i \in m}, (g_{j})_{j \in n})$ belongs either (1) to $\Phi^{\mathrm{fl}}_{r} \circ \Phi^{\mathrm{fl}}_{r}$ or (2) to $(\curlywedge \times \curlywedge) [\Phi^{\mathrm{fl}}_{r} \times \Phi^{\mathrm{fl}}_{r}] \cap (\mathrm{Mor}(\mathsf{C})^{+}_{\Phi^{\mathrm{ob}}})^{2}$.

If (1), then there exists an element $(h_{k})_{k \in p} \in \mathrm{Mor}(\mathsf{C})^{+}_{\Phi^{\mathrm{ob}}}$ such that the pairs $((f_{i})_{i \in m}, (h_{k})_{k \in p})$ and $((h_{k})_{k \in p}, (g_{j})_{j \in n})$ are in $\Phi^{\mathrm{fl}}_{r}$. By induction, $\Phi^{\mathrm{fl}}_{r}$ is a symmetric relation. Hence, the pairs $((h_{k})_{k \in p}, (f_{i})_{i \in m})$ and $((g_{j})_{j \in n}, (h_{k})_{k \in p})$ belong to $\Phi^{\mathrm{fl}}_{r}$. Thus, the pair $((g_{j})_{j \in n}, (f_{i})_{i \in m})$ also belongs to $\Phi^{\mathrm{fl}}_{r+1}$.

If (2), then there are two pairs $((f'_{i'})_{i' \in m'}, (g'_{j'})_{j' \in n'})$ and $((f''_{i''})_{i'' \in m''}, (g''_{j''})_{j'' \in n''})$ in $\Phi^{\mathrm{fl}}_{r}$ such that the pairs $(\mathrm{tg}(f'_{m'-1}), \mathrm{sc}(f''_{0}))$ and $(\mathrm{tg}(g'_{n'-1}), \mathrm{tg}(g''_{0}))$ belong to $\Phi^{\mathrm{ob}}$, $(f_{i})_{i \in m} = (f'_{i'})_{i' \in m'} \curlywedge (f''_{i''})_{i'' \in m''}$ and $(g_{j})_{j \in n} = (g'_{j'})_{j'\in n'} \curlywedge (g''_{j''})_{j'' \in n''}$. 
By induction, we have that $\Phi^{\mathrm{fl}}_{r}$ is a symmetric relation. 
Hence, the pairs $((g'_{j'})_{j' \in n'}, (f'_{i'})_{i' \in m'})$ and $((g''_{j''})_{j'' \in n''}, (f''_{i''})_{i'' \in m''})$ are in $\Phi^{\mathrm{fl}}_{r}$. 
Thus, $((g'_{j'})_{j' \in n'} \curlywedge (g''_{j''})_{j'' \in n''}, (f'_{i'})_{i' \in m'} \curlywedge (f''_{i''})_{i'' \in m''})$, which is equal to $((g_{j})_{j \in n}, (f_{i})_{i \in m})$, also belongs to $\Phi^{\mathrm{fl}}_{r+1}$.

In any case, the pair $((g_{j})_{j \in n}, (f_{i})_{i \in m})$ belongs to $\Phi^{\mathrm{fl}}_{r+1}$. Therefore, $\Phi^{\mathrm{fl}}_{r+1}$ is a symmetric relation. Hence, for every $r \in \mathbb{N}$, $\Phi^{\mathrm{fl}}_{r}$ is a symmetric relation.

This finishes the proof of the third item.
\end{proof}

As announced, we now give, for a Riguet congruence $\Phi$ on a category $\mathsf{C}$, a constructive description of $\Phi^{\natural\mathrm{fl}}$. More precisely, we prove that $\Phi^{\natural\mathrm{fl}} = \Phi^{\mathrm{fl}}_{\omega}$, where $\Phi^{\mathrm{fl}}_{\omega} = \bigcup_{r \in \mathbb{N}} \Phi^{\mathrm{fl}}_{r}$, as stated in Definition~\ref{CtvDef}.

\begin{proposition}
\label{PGCongIntG}
Let $\mathsf{C}$ be a category and $\Phi$ a Riguet congruence on $\mathsf{C}$. Then $\Phi^{\natural\mathrm{fl}} = \Phi^{\mathrm{fl}}_{\omega}$.
\end{proposition}

\begin{proof}
We first prove the inclusion from left to right. To do so, we prove that $(\Phi^{\mathrm{ob}}, \Phi^{\mathrm{fl}}_{\omega})$ is a generalized congruence on $\mathsf{C}$. 

Since $\Phi$ a Riguet congruence on $\mathsf{C}$, $\Phi^{\mathrm{ob}}$ is an equivalence relation on $\mathrm{Ob}(\mathsf{C})$. The proof that $\Phi^{\mathrm{fl}}_{\omega}$ is an equivalence relation on $\mathrm{Mor}(\mathsf{C})^{+}_{\Phi^{\mathrm{ob}}}$ follows immediately from Proposition~\ref{PGCongOpG} and is therefore omitted.

Finally, we show that $\Phi^{\mathrm{fl}}_{\omega}$ satisfies the four conditions stated in Definition~\ref{DGCgr}.

\textsf{(a)}
Let $(a,b)$ be an element of  $\Phi^{\mathrm{ob}}$.
Then $(\mathrm{id}_{a}, \mathrm{id}_{b}) \in \Phi^{\mathrm{fl}}_{\scalebox{.7}{$\left(\begin{smallmatrix}a&a\\b&b\end{smallmatrix}\right)$}} \subseteq \bigcup \Phi^{\mathrm{fl}} \subseteq \Phi^{\mathrm{fl}}_{0} \subseteq \Phi^{\mathrm{fl}}_{\omega}$. Recall the abuse of notation identifying $\eta_{\mathrm{Mor}(\mathsf{C})}^{2}[\bigcup \Phi^{\mathrm{fl}}]$ with $\bigcup \Phi^{\mathrm{fl}}$.

\textsf{(b)}
Let $((f_{i})_{i \in m}, (g_{j})_{j \in n})$ be an element of $\Phi^{\mathrm{fl}}_{\omega}$. 
We prove by induction on $r\in \mathbb{N}$ that $(\mathrm{sc}(f_{0}), \mathrm{sc}(g_0))$ and $(\mathrm{tg}(f_{m-1}), \mathrm{tg}(g_{n-1}))$ belong to $\Phi^{\mathrm{ob}}$.

\textsf{Base step of the induction.}

If $((f_{i})_{i \in m}, (g_{j})_{j \in n})$ belongs to $\Phi^{\mathrm{fl}}_{0}$, then, by definition of $\Phi^{\mathrm{fl}}_{0}$, 
this pair belongs either (1) to $\bigcup\Phi^{\mathrm{fl}}$, (2) to $\mathrm{C}$, (3) to $\mathrm{C}^{-1}$ or (4) to 
$\Delta_{\mathrm{Mor}(\mathsf{C})^{+}_{\Phi^{\mathrm{ob}}}}$.

If (1), then $m=n=1$ and $(f_{0}, g_{0})$ belongs to $\Phi^{\mathrm{fl}}_{\scalebox{.7}{$\left(\begin{smallmatrix}a&b\\a'&b'\end{smallmatrix}\right)$}}$, where $f_{0} \colon a \mor b$ and $g_{0} \colon a' \mor b'$.
Thus, by Definition~\ref{DRiguet}, $(\mathrm{sc}(f_{0}), \mathrm{sc}(g_{0}))$ and $(\mathrm{tg}(f_{m-1}), \mathrm{tg}(g_{n-1}))$ belong to $\Phi^{\mathrm{ob}}$. 

If (2), then $m = 2$, $n=1$ and $g_{0} = f_{1} \circ f_{0}$.
Thus, $\mathrm{sc}(g_{0}) = \mathrm{sc}(f_{1} \circ f_{0}) = \mathrm{sc}(f_{0})$ and $\mathrm{tg}(g_{0}) = \mathrm{tg}(f_{1} \circ f_{0}) = \mathrm{tg}(f_{1})$.
Therefore, by reflexivity of $\Phi^{\mathrm{ob}}$, $(\mathrm{sc}(f_{0}), \mathrm{sc}(g_{0}))$ and $(\mathrm{tg}(f_{m-1}), \mathrm{tg}(g_{n-1}))$ belong to $\Phi^{\mathrm{ob}}$. 

If (3), then $((g_{j})_{j \in n}, (f_{i})_{i \in m})$ belongs to $\mathrm{C}$. Hence $n=2$, $m=1$ and $f_{0} = g_{1} \circ g_{0}$. Thus, $\mathrm{sc}(f_{0}) = \mathrm{sc}(g_{1} \circ g_{0}) = \mathrm{sc}(g_{0})$ and $\mathrm{tg}(f_{0}) = \mathrm{tg}(g_{1} \circ g_{0}) = \mathrm{tg}(g_{1})$. Therefore, by reflexivity and symmetry of $\Phi^{\mathrm{ob}}$, $(\mathrm{sc}(f_{0}), \mathrm{sc}(g_{0}))$ and $(\mathrm{tg}(f_{m-1}), \mathrm{tg}(g_{n-1}))$ belong to $\Phi^{\mathrm{ob}}$. 

If (4), then $m=n$ and, for every $i \in m$, $f_{i} = g_{i}$.
Therefore, by reflexivity of $\Phi^{\mathrm{ob}}$, $(\mathrm{sc}(f_{0}), \mathrm{sc}(g_0))$ and $(\mathrm{tg}(f_{m-1}), \mathrm{tg}(g_{n-1}))$ belong to $\Phi^{\mathrm{ob}}$.

In any case, $(\mathrm{sc}(f_{0}), \mathrm{sc}(g_0))$ and $(\mathrm{tg}(f_{m-1}), \mathrm{tg}(g_{n-1}))$ belong to $\Phi^{\mathrm{ob}}$.

This completes the base step of the induction.

\textsf{Inductive step of the induction.}

Let us assume that the statement holds for $r \in \mathbb{N}$, i.e., that for every pair $((f_{i})_{i \in m}, (g_{j})_{j \in n})$ in $\Phi^{\mathrm{fl}}_{r}$, $(\mathrm{sc}(f_{0}), \mathrm{sc}(g_0))$ and $(\mathrm{tg}(f_{m-1}), \mathrm{tg}(g_{n-1}))$ belong to $\Phi^{\mathrm{ob}}$, and let us prove it for $r+1$. 
Let $((f_{i})_{i \in m}, (g_{j})_{j \in n})$ be an element of $\Phi^{\mathrm{fl}}_{r+1}$. 
Then, by definition of $\Phi^{\mathrm{fl}}_{r+1}$, $((f_{i})_{i \in m}, (g_{j})_{j \in n})$ belongs either (1) to $\Phi^{\mathrm{fl}}_{r} \circ \Phi^{\mathrm{fl}}_{r}$ or (2) to $(\curlywedge \times \curlywedge) [\Phi^{\mathrm{fl}}_{r} \times \Phi^{\mathrm{fl}}_{r}] \cap (\mathrm{Mor}(\mathsf{C})^{+}_{\Phi^{\mathrm{ob}}})^{2}$.

If (1), then there exists a family $(h_{k})_{k \in p} \in \mathrm{Mor}(\mathsf{C})^{+}_{\Phi^{\mathrm{ob}}}$ such that the pairs $((f_{i})_{i \in m}, (h_{k})_{k \in p})$ and $((h_{k})_{k \in p}, (g_{j})_{j \in n})$ belong to $\Phi^{\mathrm{fl}}_{r}$. 
By induction, we have that $(\mathrm{sc}(f_{0}), \mathrm{sc}(h_{0}))$, $(\mathrm{tg}(f_{m-1}), \mathrm{tg}(h_{k-1}))$, $(\mathrm{sc}(h_{0}), \mathrm{sc}(g_{0}))$ and $(\mathrm{tg}(h_{k-1}), \mathrm{tg}(g_{n-1}))$ belong to $\Phi^{\mathrm{ob}}$.
Therefore, since $\Phi^{\mathrm{ob}}$ is transitive, the pairs $(\mathrm{sc}(f_{0}), \mathrm{sc}(g_{0}))$ and $(\mathrm{tg}(f_{m-1}), \mathrm{tg}(g_{n-1}))$ belong to $\Phi^{\mathrm{ob}}$.

If (2), then there are two pairs $((f'_{i'})_{i' \in m'}, (g'_{j'})_{j' \in n'})$ and $((f''_{i''})_{i'' \in m''}, (g''_{j''})_{j'' \in n''})$ in $\Phi^{\mathrm{fl}}_{r}$ such that $(\mathrm{tg}(f'_{m'-1}), \mathrm{sc}(f''_{0}))$ and $(\mathrm{tg}(g'_{n'-1}), \mathrm{tg}(g''_{0}))$ belong to $\Phi^{\mathrm{ob}}$, $(f_{i})_{i \in m} = (f'_{i})_{i \in m'} \curlywedge (f''_{i})_{i \in m''}$ and $(g_{j})_{j \in n} = (g'_{j})_{j \in n'} \curlywedge (g''_{j})_{j \in n''}$. 
By induction, we have that $(\mathrm{sc}(f'_{0}), \mathrm{sc}(g'_{0}))$, $(\mathrm{tg}(f'_{m'-1}), \mathrm{tg}(g'_{n'-1}))$, $(\mathrm{sc}(f''_{0}), \mathrm{sc}(g''_{0}))$ and $(\mathrm{tg}(f''_{m''-1}), \mathrm{tg}(g''_{n''-1}))$ belong to $\Phi^{\mathrm{ob}}$.
Therefore, since
\begin{align*}
(\mathrm{sc}(f_{0}), \mathrm{sc}(g_{0})) &= (\mathrm{sc}(f'_{0}), \mathrm{sc}(g'_{0}))
\text{ and}
\\
(\mathrm{tg}(f_{m-1}), \mathrm{tg}(g_{n-1})) &= (\mathrm{tg}(f''_{m''-1}), \mathrm{tg}(g''_{n''-1})),
\end{align*}
$(\mathrm{sc}(f_{0}), \mathrm{sc}(g_{0}))$ and $(\mathrm{tg}(f_{m-1}), \mathrm{tg}(g_{n-1}))$ belong to $\Phi^{\mathrm{ob}}$.

In any case, $(\mathrm{sc}(f_{0}), \mathrm{sc}(g_0))$ and $(\mathrm{tg}(f_{m-1}), \mathrm{tg}(g_{n-1}))$ belong to $\Phi^{\mathrm{ob}}$.

This completes the inductive step of the induction.

This shows that $(\mathrm{sc}(f_{0}), \mathrm{sc}(g_0))$ and $(\mathrm{tg}(f_{m-1}), \mathrm{tg}(g_{n-1}))$ belong to $\Phi^{\mathrm{ob}}$.

\textsf{(c)}
Let $f$ and $g$ be morphisms of $\mathsf{C}$ such that $\mathrm{tg}(f) = \mathrm{sc}(g)$.
Then $((f,g), (g \circ f)) \in \mathrm{C} \subseteq \Phi^{\mathrm{fl}}_{0} \subseteq \Phi^{\mathrm{fl}}_{\omega}$.

\textsf{(d)}
Let $((f_{i})_{i \in m}, (f'_{i'})_{i' \in m'})$, $((g_{j})_{j \in n}, (g'_{j'})_{j' \in n'})$ be pairs in $\Phi^{\mathrm{fl}}_{\omega}$  satisfying that $(\mathrm{tg}(f_{m-1}), \mathrm{sc}(g_{0}))$ belongs to $\Phi^{\mathrm{ob}}$. Then there exist $r$ and $s \in \mathbb{N}$ such that the pair $((f_{i})_{i \in m}, (f'_{i'})_{i' \in m'})$ belongs to $\Phi^{\mathrm{fl}}_{r}$ and $((g_{j})_{j \in n}, (g'_{j'})_{j' \in n'})$ belongs to $\Phi^{\mathrm{fl}}_{s}$. Let $t$ be $\mathrm{max}(r, s)$. Then, by Proposition~\ref{PGCongOpG}, $((f_{i})_{i \in m}, (f'_{i'})_{i' \in m'})$ and $((g_{j})_{j \in n}, (g'_{j'})_{j' \in n'})$ belong to $\Phi^{\mathrm{fl}}_{t}$. Hence $((f_{i})_{i \in m} \curlywedge (g_{j})_{j \in n}, (f'_{i'})_{i' \in m'} \curlywedge (g'_{j'})_{j' \in n'})$ belongs to \[(\curlywedge \times \curlywedge) [\Phi^{\mathrm{fl}}_{t} \times\Phi^{\mathrm{fl}}_{t}] \cap (\mathrm{Mor}(\mathsf{C})^{+}_{\Phi^{\mathrm{ob}}})^{2} \subseteq \Phi^{\mathrm{fl}}_{t+1} \subseteq \Phi^{\mathrm{fl}}_{\omega}.\]

Therefore, we can affirm that $(\Phi^{\mathrm{ob}}, \Phi^{\mathrm{fl}}_{\omega})$ is a generalized congruence on $\mathsf{C}$.

Moreover, $(\Phi^{\mathrm{ob}}, \bigcup\Phi^{\mathrm{fl}}) \subseteq (\Phi^{\mathrm{ob}}, \Phi^{\mathrm{fl}}_{\omega})$. This is so because $\bigcup\Phi^{\mathrm{fl}} \subseteq \Phi^{\mathrm{fl}}_{0} \subseteq \Phi^{\mathrm{fl}}_{\omega}$.

From these facts, we can affirm that $\Phi^{\natural} \subseteq (\Phi^{\mathrm{ob}}, \Phi^{\mathrm{fl}}_{\omega})$, hence, in particular, that $\Phi^{\natural\mathrm{fl}} \subseteq \Phi^{\mathrm{fl}}_{\omega}$. 

To prove the other inclusion, we proceed by induction on $r \in \mathbb{N}$. 

\textsf{Base step of the induction}

Let us recall that $\Phi^{\natural}$ is the smallest generalized congruence on $\mathsf{C}$ generated by $(\Phi^{\mathrm{ob}}, \bigcup \Phi^{\mathrm{fl}})$ and $\Phi^{\natural\mathrm{ob}} = \Phi^{\mathrm{ob}}$.
Thus, we can affirm that $\bigcup \Phi^{\mathrm{fl}}$ is a subset of $\Phi^{\natural\mathrm{fl}}$.
Moreover, for every $f, g \in \mathrm{Mor}(\mathsf{C})$, if $\mathrm{sc}(g) = \mathrm{tg}(f)$, then $((f,g), (g\circ f)) \in \Phi^{\natural\mathrm{fl}}$, that is, $\mathrm{C}$ is a subset of $\Phi^{\natural\mathrm{fl}}$.
Therefore, since $\Phi^{\natural\mathrm{fl}}$ is symmetric, $\mathrm{C}^{-1}$ is also a subset of $\Phi^{\natural\mathrm{fl}}$.
Finally, since $\Phi^{\natural\mathrm{fl}}$ is reflexive, $\Delta_{\mathrm{Mor}(\mathsf{C})^{+}_{\Phi^{\mathrm{ob}}}}$ is a subset of $\Phi^{\natural\mathrm{fl}}$.

Therefore, $\Phi^{\mathrm{fl}}_{0}$ is included in $\Phi^{\natural\mathrm{fl}}$.

\textsf{Inductive step of the induction}

Let us assume that the statement holds for $r \in \mathbb{N}$, i.e., we have that $\Phi^{\mathrm{fl}}_{r} \subseteq \Phi^{\natural\mathrm{fl}}$, and let us prove it for $r+1$. Since $\Phi^{\mathrm{fl}}_{r}$ is included in $\Phi^{\natural\mathrm{fl}}$ and $\Phi^{\natural\mathrm{fl}}$ is transitive, we have that $\Phi^{\mathrm{fl}}_{r} \circ \Phi^{\mathrm{fl}}_{r}$ is included in $\Phi^{\natural\mathrm{fl}}$.

On the other hand, for every $((f_{i})_{i \in m}, (f'_{i'})_{i' \in m'})$ and $((g_{j})_{j \in n}, (g'_{j'})_{j' \in n'})$ in $\Phi^{\mathrm{fl}}_{r}$, if $(\mathrm{tg}(f_{m-1}), \mathrm{sc}(g_{0}))$ and $(\mathrm{tg}(f'_{m'-1}), \mathrm{sc}(g'_{0}))$ are in $\Phi^{\mathrm{ob}}$, then, since $\Phi^{\natural}$ is a generalized congruence, the inductive hypothesis implies that, $((f_{i})_{i \in m} \curlywedge (g_{j})_{j \in n}, (f'_{i'})_{i' \in m'} \curlywedge (g'_{j'})_{j' \in n'}) \in \Phi^{\natural\mathrm{fl}}$.
Hence, $(\curlywedge \times \curlywedge) [\Phi^{\mathrm{fl}}_{r} \times \Phi^{\mathrm{fl}}_{r}] \cap (\mathrm{Mor}(\mathsf{C})^{+}_{\Phi^{\mathrm{ob}}})^{2}$ is included in $\Phi^{\natural\mathrm{fl}}$.

Therefore, we can affirm that $\Phi^{\mathrm{fl}}_{r+1} \subseteq \Phi^{\natural\mathrm{fl}}$.

Hence, $\Phi^{\natural\mathrm{fl}} = \Phi^{\mathrm{fl}}_{\omega}$.

This finishes the proof.
\end{proof}

Our next aim is to clarify, from a category-theoretic perspective, the relationship between Riguet congruences and generalized congruences on a category. Concretely, we prove that for every category $\mathsf{C}$, there exists an adjunction from $\mathsf{RCgr}(\mathsf{C})$, the category canonically associated to the ordered set of Riguet congruences on $\mathsf{C}$, to $\mathsf{SGCgr}(\mathsf{C})$, the category arising from the ordered set of strong generalized congruences on $\mathsf{C}$. This falls under the notion of relative adjunction and constitutes another main result of this section. 

\begin{theorem}\label{AdjSGR}
Let $\mathsf{C}$ be a category. Then the corestriction of $(\bigcdot)^{\natural}$, as defined in Theorem~\ref{RCgrGCgr}, to $\mathrm{SGCgr}(\mathsf{C})$---which, by abuse of notation, we denote by the same symbol---is a functor from $\mathsf{RCgr}(\mathsf{C})$ to $\mathsf{SGCgr}(\mathsf{C})$ and has a right adjoint $(\bigcdot)^{\flat}$.
\end{theorem}

\begin{proof}
That the corestriction of $(\bigcdot)^{\natural}$ to $\mathrm{SGCgr}(\mathsf{C})$ is a functor from $\mathsf{RCgr}(\mathsf{C})$ to $\mathsf{SGCgr}(\mathsf{C})$ follows from the fact that $(\bigcdot)^{\natural}$ is isotone.

Let $(\bigcdot)^{\flat}$ be the mapping from $\mathrm{GCgr}(\mathsf{C})$ to $\mathrm{RCgr}(\mathsf{C})$ that sends a strong generalized congruence $\Psi$ on $\mathsf{C}$ to the ordered pair 
$\Psi^{\flat} = (\Psi^{\flat\mathrm{ob}},\Psi^{\flat\mathrm{fl}})$, where $\Psi^{\flat\mathrm{ob}}$ is $\Psi^{\mathrm{ob}}$ and, for every $\scalebox{0.7}{$\left(\begin{smallmatrix}a&b\\a'& b'\end{smallmatrix}\right)$} \in \Psi^{\mathrm{ob}}\times \Psi^{\mathrm{ob}}$, 
\[
\Psi^{\flat\mathrm{fl}}_{\scalebox{.7}{$\left(\begin{smallmatrix}a&b\\a'&b'\end{smallmatrix}\right)$}} = 
\left\{(f,f')\in \mathrm{Hom}_{\mathsf{C}}(a,b) \times \mathrm{Hom}_{\mathsf{C}}(a',b')\bigm | 
((f),(f'))\in \Psi^{\mathrm{fl}}\right\}.
\]
We leave it to the reader to verify that $\Psi^{\flat}$ is a Riguet congruence on $\mathsf{C}$. 
 
It is immediate that the mapping $(\bigcdot)^{\flat}$ is isotone, i.e., that, for every pair of strong generalized congruences $\Psi$ and $\Omega$ on $\mathsf{C}$, if $\Psi\subseteq \Omega$, then $\Psi^{\flat}\subseteq \Omega^{\flat}$. Moreover, if $\Phi$ is a Riguet congruence and $\Psi$ a strong generalized congruence, then $\Phi^{\natural}\subseteq \Psi$ if, and only if, $\Phi\subseteq \Psi^{\flat}$. This follows immediately from the definitions. Hence, the functor $(\bigcdot)^{\flat}$ from $\mathsf{SGCgr}(\mathsf{C})$ to $\mathsf{RCgr}(\mathsf{C})$ is such that $(\bigcdot)^{\natural}\dashv(\bigcdot)^{\flat}$.
\end{proof}

Since $\mathsf{RCgr}(\mathsf{C})$ and $\mathsf{SGCgr}(\mathsf{C})$ are categories canonically associated with ordered sets, the following corollary follows immediately.

\begin{corollary}
Let $\mathsf{C}$ be a category, and let $\eta$ and $\varepsilon$ denote the unit and counit, respectively, of the adjunction $(\bigcdot)^{\natural}\dashv(\bigcdot)^{\flat}$. Then  the full subcategory $\mathsf{RCgr}(\mathsf{C})_{\mathsf{fix}}$ of $\mathsf{RCgr}(\mathsf{C})$ determined by those Riguet congruences $\Phi$ on $\mathsf{C}$ for which $\eta_{\Phi}\colon \Phi\mor \Phi^{\natural\flat}$ is an isomorphism (equivalently, $\Phi = \Phi^{\natural\flat}$), is isomorphic to the full subcategory $\mathsf{SGCgr}(\mathsf{C})_{\mathsf{fix}}$ of $\mathsf{SGCgr}(\mathsf{C})$ determined by those strong generalized congruences $\Psi$ on $\mathsf{C}$ for which $\varepsilon_{\Psi}\colon \Psi^{\flat\natural}\mor \Psi$ is an isomorphism (equivalently, $\Psi = \Psi^{\flat\natural}$). 
\end{corollary}

In light of Theorem~\ref{AdjSGR} and taking into account the fundamental role of adjunctions in mathematics, we would not regard Riguet's notion of congruence as outdated or practically incorrect. In this sense, we hope that this work may serve to draw some attention to Riguet's importance in the history of category theory. 
Moreover, the adjunction established above substantiates the remark of Bednarczyk \textit{et al.} in~\cite{BBP99} (p.~273) that ``Often, important examples of generalized congruences are principal congruences generated by a relation on single morphisms only.'' In this formal sense, Riguet congruences serve as the principal generators for the category of strong generalized congruences ($\mathsf{SGCgr}(\mathsf{C})$). 

Next, after defining the categories $\mathsf{RCCat}$, of Riguet classified categories, and $\mathsf{GCCat}$, of  generalized classified categories, we prove that there exists a functor $(\bigcdot)^{\natural}$ from $\mathsf{RCCat}$ to $\mathsf{GCCat}$, and that there are functors $Q_{\mathrm{R}}$, $U_{\mathrm{R}}$ from $\mathsf{RCCat}$ to $\mathsf{Cat}$, and $Q_{\mathrm{G}}$, $U_{\mathrm{G}}$ from $\mathsf{GCCat}$ to $\mathsf{Cat}$, together with natural transformations $\mathrm{pr}_{\mathrm{R}}$ from $U_{\mathrm{R}}$ to $Q_{\mathrm{R}}$, $\mathrm{pr}_{\mathrm{G}}$ from $U_{\mathrm{G}}$ to $Q_{\mathrm{G}}$ and $\alpha$ from $Q_{\mathrm{R}}$ to $Q_{\mathrm{G}}\circ (\bigcdot)^{\natural}$.  
Furthermore, we prove that the corestriction of the functor $(\bigcdot)^{\natural}$ to $\mathsf{SGCCat}$ has a right adjoint $(\bigcdot)^{\flat}$ and that there are functors $Q_{\mathrm{S}}$, $U_{\mathrm{S}}$ from $\mathsf{SGCCat}$ to $\mathsf{Cat}$ together with natural transformations $\mathrm{pr}_{\mathrm{S}}$ from $U_{\mathrm{S}}$ to $Q_{\mathrm{S}}$ and $\beta$ from $Q_{\mathrm{R}}\circ (\bigcdot)^{\flat}$ to $Q_{\mathrm{S}}$. 

\begin{definition}
A \emph{Riguet classified category} is an ordered pair $(\mathsf{C},\Phi)$ in which $\mathsf{C}$ is a category and $\Phi$ a Riguet congruence on $\mathsf{C}$. A \emph{morphism} from a Riguet classified category $(\mathsf{C},\Phi)$ to another $(\mathsf{D},\Upsilon)$ is a functor $F$ from $\mathsf{C}$ to $\mathsf{D}$ which is $(\Phi,\Upsilon)$-compatible, i.e., 
\begin{enumerate}
\item for every $(a,b)\in \mathrm{Ob}(\mathsf{C})\times \mathrm{Ob}(\mathsf{C})$, 
if $(a,b)\in\Phi^{\mathrm{ob}}$, then $(F(a),F(b))\in\Upsilon^{\mathrm{ob}}$ and, 
\item for every $\left(\begin{smallmatrix}a&b\\a'& b'\end{smallmatrix}\right)\in\Phi^{\mathrm{ob}}\times \Phi^{\mathrm{ob}}$ and every $(f,f')\in \mathrm{Mor}(\mathsf{C})\times \mathrm{Mor}(\mathsf{C})$, 
if $(f,f')\in \Phi^{\mathrm{fl}}_{\scalebox{0.7}{$\left(\begin{smallmatrix}a&b\\a&b'\end{smallmatrix}\right)$}}$, then $(F(f),F(f'))\in \Upsilon^{\mathrm{fl}}_{\scalebox{0.7}{$\left(\begin{smallmatrix}F(a)&F(b)\\F(a')&F(b')\end{smallmatrix}\right)$}}$.
\end{enumerate}
We let $\mathsf{RCCat}$ stand for the corresponding category. We denote by $Q_{\mathrm{R}}$ the functor from $\mathsf{RCCat}$ to $\mathsf{Cat}$ that sends $(\mathsf{C},\Phi)$ to $\mathsf{C}/\Phi$ and a morphism $F$ from $(\mathsf{C},\Phi)$ to $(\mathsf{D},\Upsilon)$ to $Q_{\mathrm{R}}(F) = (P_{\Upsilon}\circ F)^{\sharp}$. Moreover, we denote by $U_{\mathrm{R}}$ the functor from $\mathsf{RCCat}$ to $\mathsf{Cat}$ that sends $(\mathsf{C},\Phi)$ to $\mathsf{C}$. 

A \emph{generalized classified category} is an ordered pair $(\mathsf{C},\Psi)$ in which $\mathsf{C}$ is a category and $\Psi$ a generalized congruence on $\mathsf{C}$. A \emph{morphism} from a generalized classified category $(\mathsf{C},\Psi)$ to another $(\mathsf{D},\Omega)$ is a functor $F$ from $\mathsf{C}$ to $\mathsf{D}$ which is $(\Psi,\Omega)$-compatible., i.e., 
\begin{enumerate}
\item for every $(a,b)\in \mathrm{Ob}(\mathsf{C})\times \mathrm{Ob}(\mathsf{C})$, 
if $(a,b)\in\Psi^{\mathrm{ob}}$, then $(F(a),F(b))\in\Omega^{\mathrm{ob}}$ and, 
\item for every $(f_{i})_{i \in m}, (g_{j})_{j \in n} \in \mathrm{Mor}(\mathsf{C})^{+}_{\Psi^{\mathrm{ob}}}$, if $((f_{i})_{i \in m}, (g_{j})_{j\in n}) \in \Psi^{\mathrm{fl}}$, then 
$((F(f_{i}))_{i \in m}, (F(g_{j}))_{j\in n}) \in \Omega^{\mathrm{fl}}$.
\end{enumerate}
We let $\mathsf{GCCat}$ stand for the corresponding category.
We denote by $Q_{\mathrm{G}}$ the functor from $\mathsf{GCCat}$ to $\mathsf{Cat}$ that sends $(\mathsf{C},\Psi)$ to $\mathsf{C}/\Psi$ and a morphism $F$ from $(\mathsf{C},\Psi)$ to $(\mathsf{D},\Omega)$ to $Q_{\mathrm{G}}(F) = (P_{\Omega}\circ F)^{\sharp}$.  Moreover, we denote by $U_{\mathrm{G}}$ the functor from $\mathsf{GCCat}$ to $\mathsf{Cat}$ that sends $(\mathsf{C},\Psi)$ to $\mathsf{C}$. 

A \emph{strong generalized classified category} is an ordered pair $(\mathsf{C},\Psi)$ in which $\mathsf{C}$ is a category and $\Psi$ a strong generalized congruence on $\mathsf{C}$. We let $\mathsf{SGCCat}$ stand for the full subcategory of $\mathsf{GCCat}$ determined by the strong generalized classified categories. We denote by 
$Q_{\mathrm{S}}$ and $U_{\mathrm{S}}$ the restrictions of  $Q_{\mathrm{G}}$ and $U_{\mathrm{G}}$ to $\mathsf{SGCCat}$.
\end{definition}

The following proposition, together with Theorem~\ref{RCgrGCgr}, will allow us to prove the existence of a functor from $\mathsf{RCCat}$ to $\mathsf{GCCat}$.

\begin{proposition}\label{PNFuncMor}
Let $(\mathsf{C}, \Phi)$ and $(\mathsf{D}, \Upsilon)$ be Riguet classified categories and let $F$ be a morphism of $\mathsf{RCCat}$ from $(\mathsf{C}, \Phi)$ to $(\mathsf{D}, \Upsilon)$. Then $F$ is a morphism of $\mathsf{GCCat}$ from $(\mathsf{C}, \Phi^{\natural})$ to $(\mathsf{D}, \Upsilon^{\natural})$.
\end{proposition}

\begin{proof}
The proof that if $(a,b)\in\Phi^{\natural\mathrm{ob}}$, then $(F(a),F(b))\in\Upsilon^{\natural\mathrm{ob}}$ follows from the fact that, by Theorem~\ref{RCgrGCgr}, we have that $\Phi^{\natural\mathrm{ob}}= \Phi^{\mathrm{ob}}$ and $\Upsilon^{\natural\mathrm{ob}} = \Upsilon^{\mathrm{ob}}$.

It remains to prove that if $(f_{i})_{i \in m}, (g_{j})_{j \in n} \in \mathrm{Mor}(\mathsf{C})^{+}_{\Psi^{\natural\mathrm{ob}}}$ and $((f_{i})_{i \in m}, (g_{j})_{j\in n}) \in \Psi^{\natural\mathrm{fl}}$, then 
$((F(f_{i}))_{i \in m}, (F(g_{j}))_{j\in n}) \in \Upsilon^{\natural\mathrm{fl}}$. The proof of this follows from Theorem~\ref{RCgrGCgr} and the fact that, by Proposition~\ref{PGCongIntG}, we have that $\Phi^{\natural\mathrm{fl}} = \Phi^{\mathrm{fl}}_{\omega}$ and $\Psi^{\natural\mathrm{fl}} = \Psi^{\mathrm{fl}}_{\omega}$, and it is carried out by induction. This induction proceeds as in Proposition~\ref{PGCongIntG}, distinguishing cases as in Definition~\ref{CtvDef}, and is therefore omitted.
\end{proof}

The following proposition, together with Theorem~\ref{AdjSGR}, will allow us to prove the existence of a functor from $\mathsf{SGCCat}$ to $\mathsf{RCCat}$. The proof follows directly from the definitions and is therefore omitted.

\begin{proposition}\label{PFFuncMor}
Let $(\mathsf{C}, \Psi)$ and $(\mathsf{D}, \Omega)$ be strong generalized classified categories and let $F$ be a morphism of $\mathsf{SGCCat}$ from $(\mathsf{C}, \Psi)$ to $(\mathsf{D}, \Omega)$. Then $F$ is a morphism of $\mathsf{RCCat}$ from $(\mathsf{C}, \Psi^{\flat})$ to $(\mathsf{D}, \Omega^{\flat})$.
\end{proposition}

The notation $(\bigcdot)^{\natural}$ and $(\bigcdot)^{\flat}$ used in the following theorem is consistent with its use in Theorems~\ref{RCgrGCgr} and~\ref{AdjSGR}, respectively.

\begin{theorem}\label{DNFunc}
There exists a functor $(\bigcdot)^{\natural}$ from $\mathsf{RCCat}$ to $\mathsf{GCCat}$ whose image is contained in $\mathsf{SGCCat}$. Moreover, the corestriction of $(\bigcdot)^{\natural}$ to $\mathsf{SGCCat}$---which, by abuse of notation, we denote by the same symbol---has a right adjoint $(\bigcdot)^{\flat}$.
\end{theorem}

\begin{proof}
Let $(\bigcdot)^{\natural}$ be the assignment defined as follows:
\begin{enumerate}
\item for every Riguet classified category $(\mathsf{C}, \Phi)$, $(\mathsf{C}, \Phi)^{\natural} = (\mathsf{C}, \Phi^{\natural})$, and
\item for every morphism $F$ from the Riguet classified category $(\mathsf{C}, \Phi)$ to $(\mathsf{D}, \Upsilon)$, $(F)^{\natural} = F$.
\end{enumerate}
Then $(\bigcdot)^{\natural}$ is a functor from $\mathsf{RCCat}$ to $\mathsf{GCCat}$, since Theorem~\ref{RCgrGCgr} ensures that $(\mathsf{C}, \Phi^{\natural})$ is a generalized classified category, and Proposition~\ref{PNFuncMor} ensures that $F$ is a morphism from $(\mathsf{C}, \Phi^{\natural})$ to $(\mathsf{D}, \Upsilon^{\natural})$. Moreover, by Theorem~\ref{RCgrGCgr}, $\Phi^{\natural}$ is a strong generalized congruence, and therefore the image of $(\bigcdot)^{\natural}$ is contained in $\mathsf{SGCCat}$.

Let $(\bigcdot)^{\flat}$ be the assignment defined as follows:
\begin{enumerate}
\item for every strong generalized classified category $(\mathsf{C}, \Psi)$, $(\mathsf{C}, \Psi)^{\flat} = (\mathsf{C}, \Psi^{\flat})$, and
\item for every morphism $F$ from the strong generalized classified category $(\mathsf{C}, \Psi)$ to $(\mathsf{D}, \Omega)$, $(F)^{\flat} = F$.
\end{enumerate} 
Then $(\bigcdot)^{\flat}$ is a functor from $\mathsf{SGCCat}$ to $\mathsf{RCCat}$, since Theorem~\ref{AdjSGR} ensures that $(\mathsf{C}, \Psi^{\flat})$ is a Riguet classified category, Proposition~\ref{PFFuncMor} ensures that $F$ is a morphism from $(\mathsf{C}, \Psi^{\flat})$ to $(\mathsf{D}, \Omega^{\flat})$. Moreover, taking into account the unit and counit of the adjunction $(\bigcdot)^{\natural}\dashv(\bigcdot)^{\flat}$ stated in  Theorem~\ref{AdjSGR}, $(\bigcdot)^{\flat}$ is right adjoint to $(\bigcdot)^{\natural}$.
\end{proof}

\begin{proposition}\label{nattrans}
There are natural transformations 
\begin{enumerate}
\item
$\mathrm{pr}_{\mathrm{R}}$ from $U_{\mathrm{R}}$ to $Q_{\mathrm{R}}$, 
\item
$\mathrm{pr}_{\mathrm{G}}$ from $U_{\mathrm{G}}$ to $Q_{\mathrm{G}}$,
\item
$\mathrm{pr}_{\mathrm{S}}$ from $U_{\mathrm{S}}$ to $Q_{\mathrm{S}}$, 
\item
$\alpha$ from $Q_{\mathrm{R}}$ to $Q_{\mathrm{G}}\circ (\bigcdot)^{\natural}$, \text{and}
\item
$\beta$ from $Q_{\mathrm{R}}\circ (\bigcdot)^{\flat}$ to $Q_{\mathrm{S}}$.
\end{enumerate}
\end{proposition}

\begin{proof}
\textsf{(1)}
By Proposition~\ref{PProj}, for every object $(\mathsf{C},\Phi)$ of $\mathsf{RCCat}$, we have a functor 
$P_{\Phi}$ from $\mathsf{C}$ to $\mathsf{C}/\Phi$. Then $\mathrm{pr}_{\mathrm{R}} = (P_{\Phi})_{(\mathsf{C}, \Phi) \in \mathrm{Ob}(\mathsf{RCCat})}$ is a natural transformation from $U_{\mathrm{R}}$ to $Q_{\mathrm{R}}$, i.e., $Q_{\mathrm{R}}(F) \circ P_{\Phi} = P_{\Upsilon} \circ F$. In fact, let $F$ be a morphism of $\mathsf{RCCat}$ from $(\mathsf{C}, \Phi)$ to $(\mathsf{D}, \Upsilon)$. Then, by Proposition~\ref{UPQ}, we have that 
\[
Q_{\mathrm{R}}(F) \circ P_{\Phi}
=
(P_{\Upsilon} \circ F)^{\sharp} \circ P_{\Phi}
=
P_{\Upsilon} \circ F.
\]

\textsf{(2)}
This follows by the same method as in (1).

\textsf{(3)}
This follows by the same method as in (1).

\textsf{(4)}
Let $(\mathsf{C},\Phi)$ be an object of $\mathsf{RCCat}$. Since $\Phi^{\natural} = \mathrm{GCg}_{\mathsf{C}}(\Phi^{\mathrm{ob}}, \bigcup \Phi^{\mathrm{fl}})$, we have that: 
\begin{enumerate}
\item for every $(a,a')\in \mathrm{Ob}(\mathsf{C})\times \mathrm{Ob}(\mathsf{C})$, if $(a,a')\in\Phi^{\mathrm{ob}}$, then $P_{\Phi^{\natural}}(a) = [a]_{\Phi^{\natural}} = [a']_{\Phi^{\natural}} = P_{\Phi^{\natural}}(a')$, and 
\item for every $\left(\begin{smallmatrix}a&b\\a'& b'\end{smallmatrix}\right)\in\Phi^{\mathrm{ob}}\times \Phi^{\mathrm{ob}}$ and every $(f,f')\in \mathrm{Mor}(\mathsf{C})\times \mathrm{Mor}(\mathsf{C})$, 
if $(f,f')\in \Phi^{\mathrm{fl}}_{\scalebox{0.7}{$\left(\begin{smallmatrix}a&b\\a&b'\end{smallmatrix}\right)$}}$, then $P_{\Phi^{\natural}}(f) = [f]_{\Phi^{\natural\mathrm{fl}}} = [f']_{\Phi^{\natural\mathrm{fl}}} = P_{\Phi^{\natural}}(f')$.
\end{enumerate}
Hence, by Proposition~\ref{UPQ}, there exists a unique functor $\alpha_{(\mathsf{C},\Phi)} = P_{\Phi^{\natural}}^{\sharp}$ from $\mathsf{C}/\Phi$ to $\mathsf{C}/\Phi^{\natural}$ such that $\alpha_{(\mathsf{C}, \Phi)} \circ P_{\Phi} = P_{\Phi^{\natural}}$.

We now prove that $\alpha = (\alpha_{(\mathsf{C},\Phi)})_{(\mathsf{C},\Phi)\in \mathrm{Ob}(\mathsf{RCCat})}$ is a natural transformation from $Q_{\mathrm{R}}$ to $Q_{\mathrm{G}}\circ (\bigcdot)^{\natural}$, i.e., that, for every morphism $F$ of $\mathsf{RCCat}$ from $(\mathsf{C}, \Phi)$ to $(\mathsf{D}, \Upsilon)$, $Q_{\mathrm{G}}(F) \circ \alpha_{(\mathsf{C}, \Phi)} = \alpha_{(\mathsf{D}, \Upsilon)} \circ Q_{\mathrm{R}}(F)$. To this aim, we start by noting that the functor $P_{\Upsilon^{\natural}} \circ F$ from $\mathsf{C}$ to $\mathsf{D}/\Upsilon^{\natural}$ is such that 
\begin{enumerate}
\item for every $a, a' \in \mathrm{Ob}(\mathsf{C})$, if $(a,a') \in \Phi^{\mathrm{ob}}$, then, since $F$ is $(\Phi, \Upsilon)$-compatible, we have $(F(a), F(a')) \in \Upsilon^{\mathrm{ob}}$. 
Thus $P_{\Upsilon^{\natural}}(F(a)) = P_{\Upsilon^{\natural}}(F(a'))$, and
\item for every $\left(\begin{smallmatrix}a&b\\a'& b'\end{smallmatrix}\right)\in\Phi^{\mathrm{ob}}\times \Phi^{\mathrm{ob}}$ and every $(f,f')\in \mathrm{Mor}(\mathsf{C})\times \mathrm{Mor}(\mathsf{C})$, 
if $(f,f')\in \Phi^{\mathrm{fl}}_{\scalebox{0.7}{$\left(\begin{smallmatrix}a&b\\a&b'\end{smallmatrix}\right)$}}$, then, since $F$ is $(\Phi, \Upsilon)$-compatible, we have $(F(f), F(f')) \in \Upsilon^{\mathrm{fl}}_{\scalebox{.7}{$\left(\begin{smallmatrix}F(a)&F(b)\\F(a')&F(b')\end{smallmatrix}\right)$}}$.
Thus $P_{\Upsilon^{\natural}}(F(f)) = P_{\Upsilon^{\natural}}(F(f'))$.
\end{enumerate}
Hence, by Proposition~\ref{UPQ}, there exists a unique functor $G = (P_{\Upsilon^{\natural}} \circ F)^{\sharp}$ from $\mathsf{C}/\Phi$ to $\mathsf{D}/\Upsilon^{\natural}$ such that $G \circ P_{\Phi} = P_{\Upsilon^{\natural}} \circ F$. On the other hand, the morphism $Q_{\mathrm{G}}(F) \circ \alpha_{(\mathsf{C}, \Phi)}$ is such that
\[
Q_{\mathrm{G}}(F) \circ \alpha_{(\mathsf{C}, \Phi)} \circ P_{\Phi}
=
Q_{\mathrm{G}}(F) \circ P_{\Phi^{\natural}}
= 
P_{\Upsilon^{\natural}} \circ F.
\]
In the just stated chain of equalities the first equality follows from the definition of the functor $\alpha_{(\mathsf{C}, \Phi)}$, and the second equality follows from the fact that $\mathrm{pr}_{\mathrm{G}}$ is a natural transformation from $U_{\mathrm{G}}$ to $Q_{\mathrm{G}}$.

Moreover, the morphism $\alpha_{(\mathsf{D}, \Upsilon)} \circ Q_{\mathrm{R}}(F)$ is such that
\[
\alpha_{(\mathsf{D}, \Upsilon)} \circ Q_{\mathrm{R}}(F) \circ P_{\Phi}
=
\alpha_{(\mathsf{D}, \Upsilon)} \circ P_{\Upsilon} \circ F
= 
P_{\Upsilon^{\natural}} \circ F.
\]
In the just stated chain of equalities the first equality follows from the fact that $\mathrm{pr}_{\mathrm{R}}$ is a natural transformation from $U_{\mathrm{R}}$ to $Q_{\mathrm{R}}$, and the second equality follows from the definition of the functor $\alpha_{(\mathsf{D}, \Upsilon)}$.

Therefore, by uniqueness, it follows that $Q_{\mathrm{G}}(F) \circ \alpha_{(\mathsf{C}, \Phi)} = \alpha_{(\mathsf{D}, \Upsilon)} \circ Q_{\mathrm{R}}(F)$.

In addition, the object mapping of the functor $\alpha_{(\mathsf{C},\Phi)}$ is, in particular, surjective (and thus essentially surjective), and their morphism mapping, restricted to the hom-sets, is full, yet the functor is not, in general, faithful. Consequently, whenever $\alpha_{(\mathsf{C},\Phi)}$ fails to be faithful, it cannot be an equivalence of categories.

\textsf{(5)}
This follows by the same method as in (4).
\end{proof}

Our next goal is to provide a sufficient condition under which, for a category $\mathsf{C}$ and a Riguet congruence $\Phi$ on it, the functor $\alpha_{(\mathsf{C},\Phi)}$ from $\mathsf{C}/\Phi$ to $\mathsf{C}/\Phi^{\natural}$ is an equivalence. To this end, we first define what it means for a choice function $F$ for the family $(\mathrm{Hom}_{\mathsf{C}}(a,b))_{(a,b) \in \Phi^{\mathrm{ob}}}$ to be $\Phi$-compatible, where $\Phi$ is a Riguet congruence on a category $\mathsf{C}$. 

\begin{definition}\label{ChfComp}
Let $\Phi$ be a Riguet congruence on a category $\mathsf{C}$ and $F$ a choice function for $(\mathrm{Hom}_{\mathsf{C}}(a,b))_{(a,b) \in \Phi^{\mathrm{ob}}}$. We will say that $F$ is \emph{$\Phi$-compatible} if
\begin{enumerate}
\item for every $a\in \mathrm{Ob}(\mathsf{C})$, $F(a,a) = \mathrm{id}_{a}$ and,  
\item for every $a$, $b$, $a'$, $b'\in \mathrm{Ob}(\mathsf{C})$, if $(a,a')$, $(b,b')$, $(a,b)$ and $(a', b') \in \Phi^{\mathrm{ob}}$, then $(F(a,b), F(a',b')) \in \Phi^{\mathrm{fl}}_{\scalebox{.7}{$\left(\begin{smallmatrix}a&b\\a'&b'\end{smallmatrix}\right)$}}$.
\end{enumerate}
Then, we let $\diamond_{F}$ stand for the mapping from $\mathrm{Mor}(\mathsf{C})^{+}_{ \Phi^{\mathrm{ob}}}$ to $\mathrm{Mor}(\mathsf{C})$ defined by recursion as follows:
\begin{enumerate}
\item $\diamond_{F}((f)) = f$, for every $f\in \mathrm{Mor}(\mathsf{C})$.
\item Suppose that the mapping $\diamond_{F}$ has been defined for all $\Phi^{\mathrm{ob}}$-composable sequence of length $m$. Let $(f_{i})_{i \in m+1}$ be a $\Phi^{\mathrm{ob}}$-composable sequences of length $m+1$. Then we define:
\[
\diamond_{F}((f_{i})_{i \in m+1}) = f_{m} \circ F(\mathrm{tg}(f_{m-1}), \mathrm{sc}(f_{m})) \circ \diamond_{F}((f_{i})_{i \in m}).
\]
\end{enumerate}
Note that, since $F$ is $\Phi$-compatible, if a sequence $(f_{i})_{i \in m}$, is strictly composable, i.e., if, for every $i \in m-1$, $\mathrm{tg}(f_{i}) = \mathrm{sc}(f_{i+1})$, then $\diamond_{F}((f_{i})_{i \in m}) = f_{m-1} \circ \cdots \circ f_{0}$.
\end{definition}

In the following two propositions, under the hypothesis that a choice function $F$ for $(\mathrm{Hom}_{\mathsf{C}}(a,b))_{(a,b) \in \Phi^{\mathrm{ob}}}$ is $\Phi$-compatible, we establish relationships between $\diamond_{F}$, $\Phi^{\natural\mathrm{fl}}$ and $\bigcup\Phi^{\mathrm{fl}}$.

\begin{proposition}
\label{PComp1}
Let $\mathsf{C}$ be a category, $\Phi$ a Riguet congruence on $\mathsf{C}$ and $F$ a choice function for $(\mathrm{Hom}_{\mathsf{C}}(a,b))_{(a,b) \in \Phi^{\mathrm{ob}}}$. If $F$ is $\Phi$-compatible, then, for every $(f_{i})_{i \in m} \in \mathrm{Mor}(\mathsf{C})^{+}_{ \Phi^{\mathrm{ob}}}$, we have that $((f_{i})_{i \in m}, (\diamond_{F}((f_{i})_{i \in m}))) \in \Phi^{\natural\mathrm{fl}}$.
\end{proposition}

\begin{proof}
We show the statement by induction on $m \in \mathbb{N}-1$.

\textsf{Base step of the induction}

If $m = 1$ and $f\in \mathrm{Mor}(\mathsf{C})$, then, for $(f_{i})_{i\in 1}$ where $f_{0} = f$, $\diamond_{F}((f)) = f$. Therefore, by reflexivity of $\Phi^{\natural\mathrm{fl}}$, it follows that $((f_{i})_{i \in 1}, (\diamond_{F}((f_{i})_{i \in 1}))) = ((f),(f)) \in \Phi^{\natural\mathrm{fl}}$.

\textsf{Inductive step of the induction}

Let us assume that the statement holds for $m \in \mathbb{N}-1$, i.e., that for every sequence $(f_{i})_{i \in m}$ of length $m$, $((f_{i})_{i \in m}, (\diamond_{F}((f_{i})_{i \in m}))) \in \Phi^{\natural\mathrm{fl}}$, and let us prove it for $m+1$.
Let $(f_{i})_{i \in m+1}$ be an element of $\mathrm{Mor}(\mathsf{C})^{+}_{\Phi^{\mathrm{ob}}}$. Thus, in particular, we have that $(\mathrm{tg}(f_{m-1}), \mathrm{sc}(f_{m})) \in \Phi^{\mathrm{ob}}$.

We first show that $((f_{i})_{i \in m+1}, (\diamond_{F}((f_{i})_{i \in m}), f_{m})) \in \Phi^{\natural\mathrm{fl}}$. By induction, we have that $((f_{i})_{i \in m}, (\diamond_{F}((f_{i})_{i \in m}))) \in \Phi^{\natural\mathrm{fl}}$. Thus, by item (d) in Definition~\ref{DGCgr}, we have that $((f_{i})_{i \in m+1}, (\diamond_{F}((f_{i})_{i \in m}), f_{m})) \in \Phi^{\natural\mathrm{fl}}$.

We now prove that $((\diamond_{F}((f_{i})_{i \in m}), f_{m}), (\diamond_{F}((f_{i})_{i \in m}), F(\mathrm{tg}(f_{m-1}), \mathrm{sc}(f_{m})), f_{m})) \in \Phi^{\natural\mathrm{fl}}$. Since $\Phi^{\natural\mathrm{fl}}$ is reflexive, we have that $((\diamond_{F}((f_{i})_{i \in m})), (\mathrm{id}_{\mathrm{tg}(f_{m-1})} \circ \diamond_{F}((f_{i})_{i \in m})))$ belongs to $\Phi^{\natural\mathrm{fl}}$.
Thus, by item (c) in Definition~\ref{DGCgr} and since $\Phi^{\natural\mathrm{fl}}$ is symmetric, we also have  that $((\mathrm{id}_{\mathrm{tg}(f_{m-1})} \circ \diamond_{F}((f_{i})_{i \in m})), (\diamond_{F}((f_{i})_{i \in m}), \mathrm{id}_{\mathrm{tg}(f_{m-1})}))$ belongs to $\Phi^{\natural\mathrm{fl}}$.
Moreover, since $F$ is $\Phi$-compatible and $(\mathrm{tg}(f_{m-1}), \mathrm{tg}(f_{m-1}))$ and $(\mathrm{tg}(f_{m-1}), \mathrm{sc}(f_{m}))$ belongs to $\Phi^{\mathrm{ob}}$, we have that $((\mathrm{id}_{\mathrm{tg}(f_{m-1})}), (F(\mathrm{tg}(f_{m-1}), \mathrm{sc}(f_{m}))))$ belongs to $\Phi^{\natural\mathrm{fl}}$.
Therefore, by item (d) in Definition~\ref{DGCgr}, we have that 
\[
( (\diamond_{F}((f_{i})_{i \in m}), \mathrm{id}_{\mathrm{tg}(f_{m-1})}), (\diamond_{F}((f_{i})_{i \in m}), F(\mathrm{tg}(f_{m-1}), \mathrm{sc}(f_{m}))) ) \in \Phi^{\natural\mathrm{fl}}.
\]
Hence, since $\Phi^{\natural\mathrm{fl}}$ is transitive, we have that
\[
((\diamond_{F}((f_{i})_{i \in m}), f_{m}), (\diamond_{F}((f_{i})_{i \in m}), F(\mathrm{tg}(f_{m-1}), \mathrm{sc}(f_{m})), f_{m})) \in \Phi^{\natural\mathrm{fl}}.
\]

Finally, by applying items (c) and (d) in Definition~\ref{DGCgr}, we have that
\begin{align*}
\left( (\diamond_{F}((f_{i})_{i \in m}), f_{m}), (F(\mathrm{tg}(f_{m-1}), \mathrm{sc}(f_{m})) \circ \diamond_{F}((f_{i})_{i \in m}), f_{m}) \right) 
&\in 
\Phi^{\natural\mathrm{fl}}, \text{ and}
\\
\left((F(\mathrm{tg}(f_{m-1}), \mathrm{sc}(f_{m})) \circ \diamond_{F}((f_{i})_{i \in m}), f_{m}), (\diamond_{F}((f_{i})_{i \in m+1}))\right) 
&\in 
\Phi^{\natural\mathrm{fl}}
\end{align*}

Thus, since $\Phi^{\natural\mathrm{fl}}$ is transitive, 
$
((f_{i})_{i \in m+1}, (\diamond_{F}((f_{i})_{i \in m+1}))) \in \Phi^{\natural\mathrm{fl}}.
$
\end{proof}

\begin{proposition}\label{PComp2}
Let $\mathsf{C}$ be a category, $\Phi$ a Riguet congruence on $\mathsf{C}$, $F$ a choice function for $(\mathrm{Hom}_{\mathsf{C}}(a,b))_{(a,b) \in \Phi^{\mathrm{ob}}}$ and $(f_{i})_{i \in m}, (g_{j})_{j \in n} \in \mathrm{Mor}(\mathsf{C})^{+}_{ \Phi^{\mathrm{ob}}}$. If $F$ is $\Phi$-compatible, then the following statements are equivalent:
\begin{enumerate}
\item
$((f_{i})_{i \in m}, (g_{j})_{j \in n}) \in \Phi^{\natural\mathrm{fl}}$.
\item
$(\diamond_{F}((f_{i})_{i \in m}), \diamond_{F}((g_{j})_{j \in n})) \in \bigcup\Phi^{\mathrm{fl}}$.
\end{enumerate}
In particular, $\Phi^{\natural\mathrm{fl}} \cap \mathrm{Im}\left(\eta^{2}_{\mathrm{Mor}(\mathsf{C})}\right) = \eta^{2}_{\mathrm{Mor}(\mathsf{C})}[\bigcup\Phi^{\mathrm{fl}}]$.
\end{proposition}

\begin{proof}
We begin by proving that (1) implies (2).

Let us recall that, by Proposition~\ref{PGCongIntG}, $\Phi^{\natural\mathrm{fl}} = \Phi^{\mathrm{fl}}_{\omega}$. Thus, we prove by induction on $r \in \mathbb{N}$ that, for every $(f_{i})_{i \in m}, (g_{j})_{j \in n} \in \mathrm{Mor}(\mathsf{C})^{+}_{ \Phi^{\mathrm{ob}}}$, if $((f_{i})_{i \in m}, (g_{j})_{j \in n}) \in \Phi^{\mathrm{fl}}_{r}$, then $(\diamond_{F}((f_{i})_{i \in m}), \diamond_{F}((g_{j})_{j \in n})) \in \bigcup\Phi^{\mathrm{fl}}$.

\textsf{Base step of the induction}

If $((f_{i})_{i \in m}, (g_{j})_{j \in n})$ belongs to $\Phi^{\mathrm{fl}}_{0}$, then, by definition of $\Phi^{\mathrm{fl}}_{0}$, $((f_{i})_{i \in m}, (g_{j})_{j \in n})$ belongs either (1) to $\bigcup\Phi^{\mathrm{fl}}$, (2) to $\mathrm{C}$, (3) to $\mathrm{C}^{-1}$ or (4) to $\Delta_{\mathrm{Mor}(\mathsf{C})^{+}_{\Phi^{\mathrm{ob}}}}$.

If (1), then $m=n=1$ and $(f_{0}, g_{0})$ belongs to $\Phi^{\mathrm{fl}}_{\scalebox{.7}{$\left(\begin{smallmatrix}a&b\\a'&b'\end{smallmatrix}\right)$}} \subseteq \bigcup\Phi^{\mathrm{fl}}$.  

If (2), then $m = 2$, $n=1$ and $g_{0} = f_{1} \circ f_{0}$. 
In particular, the sequence $(f_{i})_{i \in 2}$ is composable and $\diamond_{F}((f_{i})_{i \in 2}) = f_{1} \circ f_{0} = g_{0}$.
Thus $(\diamond_{F}((f_{i})_{i \in 2}), g_{0}) \in \Delta_{\mathrm{Hom}_{\mathsf{C}}(a,b)} \subseteq \Phi^{\mathrm{fl}}_{\scalebox{.7}{$\left(\begin{smallmatrix}a&b\\a&b\end{smallmatrix}\right)$}} \subseteq \bigcup \Phi^{\mathrm{fl}}$.

If (3), then $((g_{j})_{j \in n}, (f_{i})_{i \in m})$ belongs to $\mathrm{C}$. 
Hence $n=2$, $m=1$ and $f_{0} = g_{1} \circ g_{0}$.
In particular, the sequence $(g_{j})_{j \in 2}$ is composable and $\diamond_{F}((g_{j})_{j \in 2}) = g_{1} \circ g_{0} = f_{0}$.
Thus $(\diamond_{F}((g_{j})_{j \in 2}), f_{0}) \in \Delta_{\mathrm{Hom}_{\mathsf{C}}(a,b)} \subseteq \Phi^{\mathrm{fl}}_{\scalebox{.7}{$\left(\begin{smallmatrix}a&b\\a&b\end{smallmatrix}\right)$}} \subseteq \bigcup \Phi^{\mathrm{fl}}$.

If (4), then $m=n$ and, for every $i \in m$, $f_{i} = g_{i}$.
Thus, $\diamond_{F}((f_{i})_{i \in m}) = \diamond_{F}((g_{j})_{j \in n})$ and $(\diamond_{F}((f_{i})_{i \in m}), \diamond_{F}((g_{j})_{j \in n})) \in \Delta_{\mathrm{Hom}_{\mathsf{C}}(a,b)} \subseteq \Phi^{\mathrm{fl}}_{\scalebox{.7}{$\left(\begin{smallmatrix}a&b\\a&b\end{smallmatrix}\right)$}} \subseteq \bigcup \Phi^{\mathrm{fl}}$.

In any case $(\diamond_{F}((f_{i})_{i \in m}), \diamond_{F}((g_{j})_{j \in n}))$ belongs to $\bigcup\Phi^{\mathrm{fl}}$.

This completes the base step of the induction.

\textsf{Inductive step of the induction}

Let us assume that the statement holds for $r \in \mathbb{N}$, i.e., that for every pair $((f_{i})_{i \in m}, (g_{j})_{j \in n})$ in $\Phi^{\mathrm{fl}}_{r}$, it follows that $(\diamond_{F}((f_{i})_{i \in m}), \diamond_{F}((g_{j})_{j \in n}))$ belongs to $\bigcup\Phi^{\mathrm{fl}}$, and let us prove it for $r+1$. Let $((f_{i})_{i \in m}, (g_{j})_{j \in n})$ be an element of $\Phi^{\mathrm{fl}}_{r+1}$. Then, by definition of $\Phi^{\mathrm{fl}}_{r+1}$, $((f_{i})_{i \in m}, (g_{j})_{j \in n})$ belongs either (1) to $\Phi^{\mathrm{fl}}_{r} \circ \Phi^{\mathrm{fl}}_{r}$ or (2) to $(\curlywedge \times \curlywedge) [\Phi^{\mathrm{fl}}_{r} \times \Phi^{\mathrm{fl}}_{r}] \cap (\mathrm{Mor}(\mathsf{C})^{+}_{\Phi^{\mathrm{ob}}})^{2}$.

If (1), then there exists a family $(h_{k})_{k \in p} \in \mathrm{Mor}(\mathsf{C})^{+}_{\Phi^{\mathrm{ob}}}$ such that the pairs $((f_{i})_{i \in m}, (h_{k})_{k \in p})$ and $((h_{k})_{k \in p}, (g_{j})_{j \in n})$ belong to $\Phi^{\mathrm{fl}}_{r}$. 
By induction, we have that $(\diamond_{F}((f_{i})_{i \in m}), \diamond_{F}((h_{k})_{k \in p}))$ and $(\diamond_{F}((h_{k})_{k \in p}), \diamond_{F}((g_{j})_{j \in n}))$ belong to $\bigcup\Phi^{\mathrm{fl}}$.
Thus, since $\bigcup\Phi^{\mathrm{fl}}$ is transitive, we have that $(\diamond_{F}((f_{i})_{i \in m}), \diamond_{F}((g_{j})_{j \in n}))$ belongs to $\bigcup\Phi^{\mathrm{fl}}$.

If (2), then there are two pairs $((f'_{i'})_{i' \in m'}, (g'_{j'})_{j' \in n'})$ and $((f''_{i''})_{i'' \in m''}, (g''_{j''})_{j'' \in n''})$ in $\Phi^{\mathrm{fl}}_{r}$ such that $(\mathrm{tg}(f'_{m'-1}), \mathrm{sc}(f''_{0}))$ and $(\mathrm{tg}(g'_{n'-1}), \mathrm{tg}(g''_{0}))$ belong to $\Phi^{\mathrm{ob}}$, $(f_{i})_{i \in m} = (f'_{i'})_{i' \in m'} \curlywedge (f''_{i''})_{i'' \in m''}$ and $(g_{j})_{j \in n} = (g'_{j'})_{j' \in n'} \curlywedge (g''_{j''})_{j'' \in n''}$. 
By induction, we have that $(\diamond_{F}((f'_{i'})_{i' \in m'}), \diamond_{F}((g'_{j'})_{j' \in n'}))$ and $(\diamond_{F}((f''_{i''})_{i'' \in m''}), \diamond_{F}((g''_{j''})_{j'' \in n''}))$ belong to $\bigcup\Phi^{\mathrm{fl}}$.
Note that since $(f_{i})_{i \in m}$ and $(g_{j})_{j \in m}$ are $\Phi$-composable, we have that $(\mathrm{tg}(f'_{m'-1}), \mathrm{sc}(f''_{0}))$ and $(\mathrm{tg}(g'_{n'-1}), \mathrm{sc}(g''_{0}))$ belong to $\Phi^{\mathrm{ob}}$. 
Moreover, by induction, we also have that $(\mathrm{sc}(f''_{0}), \mathrm{sc}(g''_{0}))$ and $(\mathrm{tg}(f'_{m'-1}), \mathrm{tg}(g'_{n'-1}))$ belong to $\Phi^{\mathrm{ob}}$.
Therefore, since $F$ is $\Phi$-compatible, we have that 
\[
\textstyle
\left(F\left(\mathrm{tg}\left(f'_{m'-1}\right), \mathrm{sc}\left(f''_{0}\right)\right), F\left(\mathrm{tg}\left(g'_{n'-1}\right), \mathrm{tg}\left(g''_{0}\right)\right)\right) \in \bigcup \Phi^{\mathrm{fl}}.
\]
Thus, by item (e) in the Definition~\ref{DRiguet}, we have that
\begin{align*}
&(\diamond_{F}((f_{i})_{i \in m}), \diamond_{F}((g_{j})_{j \in n}))
=
\\
&\hspace{1cm}
\left(\diamond_{F}((f''_{i''})_{i'' \in m''}) \circ F\left(\mathrm{tg}\left(f'_{m'-1}\right), \mathrm{sc}\left(f''_{0}\right)\right) \circ \diamond_{F}((f'_{i'})_{i' \in m'})\right., 
\\
&\hspace{2.5cm}
\textstyle
\left.\diamond_{F}((g''_{j''})_{j'' \in n''}) \circ F\left(\mathrm{tg}\left(g'_{n'-1}\right), \mathrm{sc}\left(g''_{0}\right)\right) \circ \diamond_{F}((g'_{j'})_{j' \in n'})\right)
\in
\bigcup\Phi^{\mathrm{fl}}.
\end{align*}

In any case $(\diamond_{F}((f_{i})_{i \in m}), \diamond_{F}((g_{j})_{j \in n})$ belongs to $\bigcup\Phi^{\mathrm{fl}}$.

This completes the inductive step of the induction.

This proves that $(\diamond_{F}((f_{i})_{i \in m}), \diamond_{F}((g_{j})_{j \in n})$ belongs to $\bigcup\Phi^{\mathrm{fl}}$.

This completes the proof.

Now we prove that (2) implies (1).

According to Proposition~\ref{PComp1}, we have that the pairs $((f_{i})_{i \in m}, \diamond_{F}((f_{i})_{i \in m}))$ and $((g_{j})_{j \in n}, \diamond_{F}((g_{j})_{j \in n}))$ belong to $\Phi^{\natural\mathrm{fl}}$. Thus, since the relation $\Phi^{\natural\mathrm{fl}}$ is symmetric and transitive, it follows that $(\diamond_{F}((f_{i})_{i \in m}), \diamond_{F}((g_{j})_{j \in n})) \in \bigcup\Phi^{\mathrm{fl}}$.
\end{proof}

We are now in a position to state a sufficient condition under which the functor $\alpha_{(\mathsf{C},\Phi)}$ from $\mathsf{C}/\Phi$ to $\mathsf{C}/\Phi^{\natural}$ is an equivalence.

\begin{corollary}\label{EqvQRQN}
Let $(\mathsf{C}, \Phi)$ be a Riguet classified category. If there exists a choice function for $(\mathrm{Hom}_{\mathsf{C}}(a,b))_{(a,b) \in \Phi^{\mathrm{ob}}}$ that is $\Phi$-compatible, then $\alpha_{(\mathsf{C}, \Phi)}$ is an equivalence from $\mathsf{C}/\Phi$ to $\mathsf{C}/\Phi^{\natural}$. 
\end{corollary}

Next, within Manes' framework of categories of $\mathsf{K}$-objects with structure~\cite{em76}, we study the wide subcategory $\mathsf{RCCat}_{\mathrm{full}}$ of $\mathsf{RCCat}$, whose morphisms are the full morphisms of $\mathsf{RCCat}$, and the category $\mathsf{GCCat}$. For $\mathsf{RCCat}_{\mathrm{full}}$ we take 
$\mathsf{K} = \mathsf{Cat}_{\mathrm{full}}$, the wide subcategory of $\mathsf{Cat}$ whose morphisms are the full functors, while for $\mathsf{GCCat}$, we take $\mathsf{K} = \mathsf{Cat}$. In both cases, these categories are further related to the Grothendieck theory of fibrations~\cite{Gro71}, by considering the restriction of $U_{\mathrm{R}}$ to $\mathsf{RCCat}_{\mathrm{full}}$, the functor $U_{\mathrm{G}}$, the restriction of $(\bigcdot)^{\natural}$ to $\mathsf{RCCat}_{\mathrm{full}}$ and the canonical embedding of $\mathsf{Cat}_{\mathrm{full}}$ into $\mathsf{Cat}$.

We begin by proving the following result for $\mathsf{RCCat}_{\mathrm{full}}$: if $F\colon \mathsf{C}\mor (\mathsf{D},\Upsilon)$, where $F$ is a full functor from $\mathsf{C}$ to $\mathsf{D}$ and $\Upsilon$ a Riguet congruence on $\mathsf{D}$, then $F$ has an optimal lift.

\begin{proposition}\label{ROpLift}
Let $\mathsf{C}$ and $\mathsf{D}$ be categories, $F$ a full functor from $\mathsf{C}$ to $\mathsf{D}$ and 
$\Upsilon$ a Riguet congruence on $\mathsf{D}$. Then there exists a Riguet congruence $F^{\ast}_{\mathrm{R}}[\Upsilon]$ on $\mathsf{C}$ such that $F$ is a morphism from $(\mathsf{C},F^{\ast}_{\mathrm{R}}[\Upsilon])$ to  
$(\mathsf{D},\Upsilon)$ and, for every Riguet classified category $(\mathsf{E},\Xi)$ and every full functor $H$ from $\mathsf{E}$ to $\mathsf{C}$, if $F\circ H$ is a morphism from $(\mathsf{E},\Xi)$ to $(\mathsf{D},\Upsilon)$, then $H$ is a morphism from $(\mathsf{E},\Xi)$ to $(\mathsf{C},F^{\ast}_{\mathrm{R}}[\Upsilon])$. We will call $F^{\ast}_{\mathrm{R}}[\Upsilon]$ the optimal lift of $\Upsilon$ through $F$. 
\end{proposition}

\begin{proof}
It suffices to take as $F^{\ast}_{\mathrm{R}}[\Upsilon]$ the ordered pair 
$(F^{\ast}_{\mathrm{R}}[\Upsilon]^{\mathrm{ob}},F^{\ast}_{\mathrm{R}}[\Upsilon]^{\mathrm{fl}})$ defined as follows:
\begin{enumerate}
\item $F^{\ast}_{\mathrm{R}}[\Upsilon]^{\mathrm{ob}} = (F\times F)^{-1}[\Upsilon^{\mathrm{ob}}]$.
\item For every $\scalebox{0.7}{$\left(\begin{smallmatrix}a&b\\a'& b'\end{smallmatrix}\right)$} \in F^{\ast}_{\mathrm{R}}[\Upsilon]^{\mathrm{ob}}\times F^{\ast}_{\mathrm{R}}[\Upsilon]^{\mathrm{ob}}$, 
$F^{\ast}_{\mathrm{R}}[\Upsilon]^{\mathrm{fl}}_{\scalebox{0.7}{$\left(\begin{smallmatrix}a&b\\a'& b'\end{smallmatrix}\right)$}} = (F\times F)^{-1}\left[\Upsilon^{\mathrm{fl}}_{\scalebox{0.7}{$\left(\begin{smallmatrix}F(a)&F(b)\\F(a')& F(b')\end{smallmatrix}\right)$}}\right]$.
\end{enumerate}

The details are left to the reader.
\end{proof}

Before proving that for $\mathsf{GCCat}$: if $F\colon  \mathsf{C}\mor(\mathsf{D},\Omega)$, where $F$ is a functor from $\mathsf{C}$ to $\mathsf{D}$ and $\Omega$ a generalized congruence on $\mathsf{D}$, then $F$ has an optimal lift, we first present a preliminary lemma that will be needed to carry this out. Since the proof of the lemma is straightforward we omit the proof.

\begin{lemma}\label{ExtCompPaths}
Let $\mathsf{C}$ and $\mathsf{D}$ be categories, $F$ a functor from $\mathsf{C}$ to $\mathsf{D}$, 
$\Psi^{\mathrm{ob}}$ an equivalence relation on $\mathrm{Ob}(\mathsf{C})$ and $\Omega^{\mathrm{ob}}$ an equivalence relation on $\mathrm{Ob}(\mathsf{D})$. If $F$ is $(\Psi^{\mathrm{ob}},\Omega^{\mathrm{ob}})$-compatible, then there exists a unique mapping $(\eta_{\mathrm{Mor}(\mathsf{D})}\circ F)^{\sharp}$ from 
$\mathrm{Mor}(\mathsf{C})^{+}$ to $\mathrm{Mor}(\mathsf{D})^{+}$ such that 
$\eta_{\mathrm{Mor}(\mathsf{D})}\circ F = (\eta_{\mathrm{Mor}(\mathsf{D})}\circ F)^{\sharp}\circ \eta_{\mathrm{Mor}(\mathsf{C})}$ and 
$(\eta_{\mathrm{Mor}(\mathsf{D})}\circ F)^{\sharp}\left[\mathrm{Mor}(\mathsf{C})^{+}_{\Psi^{\mathrm{ob}}}\right]\subseteq \mathrm{Mor}(\mathsf{D})^{+}_{\Psi^{\mathrm{ob}}}$. We denote by $F^{+}$ the birestriction of $(\eta_{\mathrm{Mor}(\mathsf{D})}\circ F)^{\sharp}$ to $\mathrm{Mor}(\mathsf{C})^{+}_{\Psi^{\mathrm{ob}}}$ and $\mathrm{Mor}(\mathsf{D})^{+}_{\Omega^{\mathrm{ob}}}$.
\end{lemma}

\begin{proposition}\label{GOpLift}
Let $\mathsf{C}$ and $\mathsf{D}$ be categories, $F$ a functor from $\mathsf{C}$ to $\mathsf{D}$ and $\Omega$ a generalized congruence on $\mathsf{D}$, then there exists a generalized congruence $F^{\ast}_{\mathrm{G}}[\Omega]$ on $\mathsf{C}$ such that $F$ is a morphism from $(\mathsf{C},F^{\ast}_{\mathrm{G}}[\Omega])$ to $(\mathsf{D},\Omega)$ and, for every generalized classified category $(\mathsf{E},\Xi)$ and every functor $H$ from $\mathsf{E}$ to $\mathsf{C}$, if $F\circ H$ is a morphism from $(\mathsf{E},\Xi)$ to $(\mathsf{D},\Omega)$, then $H$ is a morphism from $(\mathsf{E},\Xi)$ to $(\mathsf{C},F^{\ast}_{\mathrm{G}}[\Omega])$. We will call $F^{\ast}_{\mathrm{G}}[\Omega]$ the optimal lift of $\Omega$ through $F$.
\end{proposition}

\begin{proof}
It suffices to take as $F^{\ast}_{\mathrm{G}}[\Omega]$ the ordered pair 
$(F^{\ast}_{\mathrm{G}}[\Omega]^{\mathrm{ob}},F^{\ast}_{\mathrm{G}}[\Omega]^{\mathrm{fl}})$ defined as follows:
\begin{enumerate}
\item $F^{\ast}_{\mathrm{G}}[\Omega]^{\mathrm{ob}} = (F\times F)^{-1}[\Omega^{\mathrm{ob}}]$.
\item $F^{\ast}_{\mathrm{G}}[\Omega]^{\mathrm{fl}} = (F^{+}\times F^{+})^{-1}[\Omega^{\mathrm{fl}}]$, where $F^{+}$ is the mapping from $\mathrm{Mor}(\mathsf{C})^{+}_{F^{\ast}_{\mathrm{G}}[\Omega]^{\mathrm{ob}}}$ to 
$\mathrm{Mor}(\mathsf{D})^{+}_{\Omega^{\mathrm{ob}}}$ obtained as a special case of Lemma~\ref{ExtCompPaths}.
\end{enumerate}

The details are left to the reader.
\end{proof}

Let $\mathsf{C}$ be an arbitrary category. By Proposition~\ref{GCgrAlgLatt}, the set $\mathrm{GCgr}(\mathsf{C})$ of all generalized congruences on $\mathsf{C}$, ordered by inclusion, is an algebraic lattice and therefore has arbitrary infima. Combined with Proposition~\ref{GOpLift}, this fact allows us to establish the following stronger result for generalized congruences.

\begin{proposition}\label{OpLift}
Let $L$ be a set, $\mathsf{C}$ a category, $(\mathsf{D}_{\lambda},\Omega_{\lambda})_{\lambda\in L}$ an $L$-indexed family of generalized classified categories and $(F_{\lambda})_{\lambda\in L}$ an $L$-indexed family of functors such that, for every $\lambda\in L$, $F_{\lambda}\colon \mathsf{C}\mor \mathsf{D}_{\lambda}$. Then the generalized congruence $\bigcap_{\lambda\in L}F^{\ast}_{\lambda,\mathrm{G}}[\Omega_{\lambda}]$ on $\mathsf{C}$ is such that, for every $\lambda\in L$, $F_{\lambda}$ is a morphism from $(\mathsf{C},\bigcap_{\lambda\in L}F^{\ast}_{\lambda,\mathrm{G}}[\Omega_{\lambda}])$ to $(\mathsf{D}_{\lambda},\Omega_{\lambda})$ and, for every generalized classified category $(\mathsf{E},\Xi)$ and every functor $H$ from $\mathsf{E}$ to $\mathsf{C}$, if, for every $\lambda\in L$, $F_{\lambda}\circ H$ is a morphism from $(\mathsf{E},\Xi)$ to $(\mathsf{D}_{\lambda},\Omega_{\lambda})$, then $H$ is a morphism from $(\mathsf{E},\Xi)$ to $(\mathsf{C},\bigcap_{\lambda\in L}F^{\ast}_{\lambda,\mathrm{G}}[\Omega_{\lambda}])$. We summarize the just stated result by saying that, for every set $L$, every $L$-indexed family $F_{\lambda}\colon \mathsf{C}\mor (\mathsf{D}_{\lambda},\Omega_{\lambda})$ has an optimal lift. 

We will call $\bigcap_{\lambda\in L}F^{\ast}_{\lambda,\mathrm{G}}[\Omega_{\lambda}]$ the optimal lift of $(\Omega_{\lambda})_{\lambda\in L}$ through $(F_{\lambda})_{\lambda\in L}$. In particular, for $(F_{\lambda})_{\lambda\in\varnothing}$, the optimal lift is $\nabla_{\mathsf{C}}$, the greatest generalized congruence on $\mathsf{C}$.
\end{proposition}

Moreover, since Proposition~\ref{OpLift} satisfies the hypothesis of Proposition 3.31 of Manes~\cite{em76} (pp.~151--152), we may assert, as stated in the following corollary, that, for every set $L$, every $L$-indexed family $F_{\lambda}\colon (\mathsf{D}_{\lambda},\Omega_{\lambda})\mor\mathsf{C}$ has a co-optimal lift. 

\begin{corollary}
Let $L$ be a set, $\mathsf{C}$ a category, $(\mathsf{D}_{\lambda},\Omega_{\lambda})_{\lambda\in L}$ an $L$-indexed family of generalized classified categories and $(F_{\lambda})_{\lambda\in L}$ an $L$-indexed family of functors such that, for every $\lambda\in L$, $F_{\lambda}\colon \mathsf{D}_{\lambda}\mor \mathsf{C}$. Then 
$F_{\lambda}\colon (\mathsf{D}_{\lambda},\Omega_{\lambda})\mor\mathsf{C}$ has a co-optimal lift. In particular, for $(F_{\lambda})_{\lambda\in\varnothing}$, the co-optimal lift is $\bigcap \mathrm{GCgr}(\mathsf{C})$, the smallest generalized congruence on $\mathsf{C}$. 
\end{corollary}

\begin{proof}
Indeed, it suffices to consider the optimal lift of the family 
$\mathrm{Id}_{\mathsf{C}}\colon \mathsf{C}\mor(\mathsf{C},\Psi)$ indexed by the set $\mathcal{C}_{(F_{\lambda})_{\lambda\in L}}$ of all generalized congruences $\Psi$ on $\mathsf{C}$ such that, for every $\lambda\in L$, $F_{\lambda}$ is a morphism from $(\mathsf{D}_{\lambda},\Omega_{\lambda})$ to $(\mathsf{C},\Psi)$ (note that $\nabla_{\mathsf{C}}$ belongs to $\mathcal{C}_{(F_{\lambda})_{\lambda\in L}}$). Thus, the co-optimal lift is $\bigcap \mathcal{C}_{(F_{\lambda})_{\lambda\in L}}$.  We call $\bigcap \mathcal{C}_{(F_{\lambda})_{\lambda\in L}}$ the co-optimal lift of $(\Omega_{\lambda})_{\lambda\in L}$ through $(F_{\lambda})_{\lambda\in L}$.
\end{proof}

\begin{remark}
Since $\mathsf{GCCat}$ is such that every family $F_{\lambda}\colon \mathsf{C}\mor (\mathsf{D}_{\lambda},\Omega_{\lambda})$ has an optimal lift and every family $F_{\lambda}\colon (\mathsf{D}_{\lambda},\Omega_{\lambda})\mor\mathsf{C}$ has a co-optimal lift, it follows from~\cite{em76} that the functor $U_{\mathrm{G}}$ from $\mathsf{GCCat}$ to $\mathsf{Cat}$ has both left and right adjoints and constructs projective and inductive limits.
\end{remark}

We next prove that the birestriction of $U_{\mathrm{R}}$ to $\mathsf{RCCat}_{\mathrm{full}}$ and $\mathsf{Cat}_{\mathrm{full}}$, which we continue to denote by $U_{\mathrm{R}}$, is a fibration, that $U_{\mathrm{G}}$ is a bifibration and that there exists a morphism of fibrations from $U_{\mathrm{R}}$ to $U_{\mathrm{G}}$. Before establishing this result, we state the following (commutation) lemma.

\begin{lemma}\label{commutation}
Let $\mathsf{C}$ and $\mathsf{D}$ be categories, $F$ a full functor from $\mathsf{C}$ to $\mathsf{D}$ and 
$\Upsilon$ a Riguet congruence on $\mathsf{D}$. Then $F^{\ast}_{\mathrm{R}}[\Upsilon]^{\natural} = F^{\ast}_{\mathrm{G}}[\Upsilon^{\natural}]$.
\end{lemma}

\begin{proof}
We begin by proving that $F^{\ast}_{\mathrm{R}}[\Upsilon]^{\natural}\subseteq F^{\ast}_{\mathrm{G}}[\Upsilon^{\natural}]$. By Proposition~\ref{ROpLift}, $F$ is a morphism from 
$(\mathsf{C},F^{\ast}_{\mathrm{R}}[\Upsilon])$ to $(\mathsf{D},\Upsilon)$ and, by Proposition~\ref{GOpLift}, $F$ is also a morphism from $(\mathsf{C},F^{\ast}_{\mathrm{G}}[\Upsilon^{\natural}])$ to $(\mathsf{D},\Upsilon^{\natural})$. Hence, by Proposition~\ref{PNFuncMor}, $F = F\circ \mathrm{Id}_{\mathsf{C}}$ is a morphism from $(\mathsf{C},F^{\ast}_{\mathrm{R}}[\Upsilon]^{\natural})$ to $(\mathsf{D},\Upsilon^{\natural})$. Therefore, $\mathrm{Id}_{\mathsf{C}}$ is a morphism from $(\mathsf{C},F^{\ast}_{\mathrm{R}}[\Upsilon]^{\natural})$ to 
$(\mathsf{C}, F^{\ast}_{\mathrm{G}}[\Upsilon^{\natural}])$, i.e., $F^{\ast}_{\mathrm{R}}[\Upsilon]^{\natural}\subseteq F^{\ast}_{\mathrm{G}}[\Upsilon^{\natural}]$.

We next prove that $F^{\ast}_{\mathrm{G}}[\Upsilon^{\natural}]\subseteq F^{\ast}_{\mathrm{R}}[\Upsilon]^{\natural}$. 
We first show that $F^{\ast}_{\mathrm{G}}[\Upsilon^{\natural}]^{\mathrm{ob}} \subseteq F^{\ast}_{\mathrm{R}}[\Upsilon]^{\natural \mathrm{ob}}$. The following chain of equalities holds
\[
F_{\mathrm{G}}^{\ast}[\Upsilon^{\natural}]^{\mathrm{ob}}
=
(F \times F)^{-1}[\Upsilon^{\natural \mathrm{ob}}]
=
(F \times F)^{-1}[\Upsilon^{\mathrm{ob}}]
=
F_{\mathrm{R}}^{\ast}[\Upsilon]^{\mathrm{ob}}
=
F_{\mathrm{R}}^{\ast}[\Upsilon]^{\natural\mathrm{ob}}.
\]

The first equality follows from the Definition of $F_{\mathrm{G}}^{\ast}[\Upsilon^{\natural}]$, introduced in Proposition~\ref{GOpLift};
the second equality follows from Theorem~\ref{RCgrGCgr};
the third equality follows from the Definition of $F_{\mathrm{R}}^{\ast}[\Upsilon]$, introduced in Proposition~\ref{ROpLift};
finally, the last equality follows from Theorem~\ref{RCgrGCgr}.

Now we show that $F^{\ast}_{\mathrm{G}}[\Upsilon^{\natural}]^{\mathrm{fl}} \subseteq F^{\ast}_{\mathrm{R}}[\Upsilon]^{\natural \mathrm{fl}}$. Let $((f_{i})_{i \in m}, (g_{j})_{j \in n})$ be a pair in $F^{\ast}_{\mathrm{G}}[\Upsilon^{\natural}]^{\mathrm{fl}}$, that is, according to Proposition~\ref{GOpLift}, a pair in $(F^{+} \times F^{+})^{-1}[\Upsilon^{\natural \mathrm{fl}}]$. Thus, $((F(f_{i}))_{i \in m}, (F(g_{j}))_{j \in n})$ is a pair in $\Upsilon^{\natural \mathrm{fl}}$. According to Proposition~\ref{PGCongIntG}, it follows that $((F(f_{i}))_{i \in m}, (F(g_{j}))_{j \in n})$ is a pair in $\Upsilon^{\mathrm{fl}}_{\omega}$, i.e., there exists an $r \in \mathbb{N}$ such that $((F(f_{i}))_{i \in m}, (F(g_{j}))_{j \in n}) \in \Upsilon^{\mathrm{fl}}_{r}$.

We show by induction on $r$ that $((f_{i})_{i \in m}, (g_{j})_{j \in n}) \in F_{\mathrm{R}}^{\ast}[\Upsilon]^{\natural \mathrm{fl}}$.

\textsf{Base step of the induction}

If $((F(f_{i}))_{i \in m}, (F(g_{j}))_{j \in n})$ belongs to $\Upsilon^{\mathrm{fl}}_{0}$, then, by definition of $\Upsilon^{\mathrm{fl}}_{0}$, the pair $((F(f_{i}))_{i \in m}, (F(g_{j}))_{j \in n})$ belongs either (1) to $\bigcup\Upsilon^{\mathrm{fl}}$, (2) to $\mathrm{C}$, (3) to $\mathrm{C}^{-1}$ or (4) to $\Delta_{\mathrm{Mor}(\mathsf{C})^{+}_{F_{\mathrm{G}}^{\ast}[\Upsilon^{\natural}]^{\mathrm{ob}}}}$.

If (1), then $m=n=1$ and $(F(f_{0}), F(g_{0}))$ belongs to $\Upsilon^{\mathrm{fl}}_{\scalebox{.7}{$\left(\begin{smallmatrix}F(a)&F(b)\\F(a')&F(b')\end{smallmatrix}\right)$}}$.  
Therefore,
\[
\textstyle
(f_{0}, g_{0})
\in 
(F \times F)^{-1}\left[
\Upsilon^{\mathrm{fl}}_{\scalebox{.7}{$\left(\begin{smallmatrix}F(a)&F(b)\\F(a')&F(b')\end{smallmatrix}\right)$}}
\right]
\subseteq 
\bigcup F_{\mathrm{R}}^{\ast}[\Upsilon]^{\mathrm{fl}}
\subseteq 
F_{\mathrm{R}}^{\ast}[\Upsilon]^{\mathrm{fl}}_{0}
\subseteq 
F_{\mathrm{R}}^{\ast}[\Upsilon]^{\mathrm{fl}}_{\omega}
=
F_{\mathrm{R}}^{\ast}[\Upsilon]^{\natural\mathrm{fl}}.
\]

If (2), then $m = 2$, $n=1$ and $F(g_{0}) = F(f_{1}) \circ F(f_{0}) = F(f_{1} \circ f_{0})$. In particular, $(F(f_{1} \circ f_{0}), F(g_{0}))$ belongs to $\Upsilon^{\mathrm{fl}}_{\scalebox{.7}{$\left(\begin{smallmatrix}F(a)&F(b)\\F(a)&F(b)\end{smallmatrix}\right)$}}$. Therefore,
\[
(f_{1} \circ f_{0}, g_{0})
\in 
(F \times F)^{-1}\left[
\Upsilon^{\mathrm{fl}}_{\scalebox{.7}{$\left(\begin{smallmatrix}F(a)&F(b)\\F(a)&F(b)\end{smallmatrix}\right)$}}
\right]
\subseteq
F_{\mathrm{R}}^{\ast}[\Upsilon]^{\natural\mathrm{fl}}.
\]
Moreover, since $F_{\mathrm{R}}^{\ast}[\Upsilon]^{\natural\mathrm{fl}}$ is a generalized congruence, it follows that $((f_{0}, f_{1}), (f_{1} \circ f_{0}))$ belongs to $F_{\mathrm{R}}^{\ast}[\Upsilon]^{\natural\mathrm{fl}}$. Finally, by transitivity, $((f_{0}, f_{1}), (g_{0}))$ belongs to $F_{\mathrm{R}}^{\ast}[\Upsilon]^{\natural\mathrm{fl}}$.

If (3), then $((F(g_{j}))_{j \in n}, (F(f_{i}))_{i \in m})$ belongs to $\mathrm{C}$. Thus, by the second case, $((g_{j})_{j \in n}, (f_{i})_{i \in m})$ belongs to $F_{\mathrm{R}}^{\ast}[\Upsilon]^{\natural\mathrm{fl}}$. Finally, by symmetry of $F_{\mathrm{R}}^{\ast}[\Upsilon]^{\natural\mathrm{fl}}$, it follows that $((f_{i})_{i \in m}, (g_{j})_{j \in n})$ also belongs to $F_{\mathrm{R}}^{\ast}[\Upsilon]^{\natural\mathrm{fl}}$.

If (4), then $m=n$ and, for every $i \in m$, $F(f_{i}) = F(g_{i})$. In particular, for every $i \in m$, $(F(f_{i}), F(g_{i}))$ belongs to $\Upsilon^{\mathrm{fl}}_{\scalebox{.7}{$\left(\begin{smallmatrix}F(a)&F(b)\\F(a)&F(b)\end{smallmatrix}\right)$}}$. Therefore, for every $i \in m$,
\[
(f_{i}, g_{i})
\in 
(F \times F)^{-1}\left[
\Upsilon^{\mathrm{fl}}_{\scalebox{.7}{$\left(\begin{smallmatrix}F(a)&F(b)\\F(a)&F(b)\end{smallmatrix}\right)$}}
\right]
\subseteq
F_{\mathrm{R}}^{\ast}[\Upsilon]^{\natural\mathrm{fl}}.
\]
Thus, since $F_{\mathrm{R}}^{\ast}[\Upsilon]^{\natural\mathrm{fl}}$ is a generalized congruence, it follows that $((f_{i})_{i \in m}, (g_{j})_{j \in n})$ also belongs to $F_{\mathrm{R}}^{\ast}[\Upsilon]^{\natural\mathrm{fl}}$.

In any case $((f_{i})_{i \in m}, (g_{j})_{j\in n})$ belongs to $F_{\mathrm{R}}^{\ast}[\Upsilon]^{\natural\mathrm{fl}}$.

This completes the base step of the induction.

\textsf{Inductive step of the induction.}

Let us assume that the statement holds for $r \in \mathbb{N}$, i.e., that for every pair $((F(f_{i}))_{i \in m}, (F(g_{j}))_{j \in n})$ in $\Upsilon^{\mathrm{fl}}_{r}$, the pair $((f_{i})_{i \in m}, (g_{j})_{j\in n})$ is in $F_{\mathrm{R}}^{\ast}[\Upsilon]^{\natural\mathrm{fl}}$, and let us prove it for $r+1$. Let $((F(f_{i}))_{i \in m}, (F(g_{j}))_{j \in n})$ be an element of $\Upsilon^{\mathrm{fl}}_{r+1}$. Then, by definition of $\Upsilon^{\mathrm{fl}}_{r+1}$, $((F(f_{i}))_{i \in m}, (F(g_{j}))_{j \in n})$ belongs either (1) to $\Upsilon^{\mathrm{fl}}_{r} \circ \Upsilon^{\mathrm{fl}}_{r}$ or (2) to $(\curlywedge \times \curlywedge) [\Upsilon^{\mathrm{fl}}_{r} \times \Upsilon^{\mathrm{fl}}_{r}] \cap (\mathrm{Mor}(\mathsf{C})^{+}_{F^{\ast}_{\mathrm{G}}[\Upsilon^{\natural}]^{\mathrm{ob}}})^{2}$.

If (1), then there exists a family $(F(h_{k}))_{k \in p} \in \mathrm{Mor}(\mathsf{C})^{+}_{F^{\ast}_{\mathrm{G}}[\Upsilon^{\natural}]^{\mathrm{ob}}}$ such that the pairs $((F(f_{i}))_{i \in m}, (F(h_{k}))_{k \in p})$ and $((F(h_{k}))_{k \in p}, (F(g_{j}))_{j \in n})$ belong to $\Phi^{\mathrm{fl}}_{r}$. Note that the functor $F$ is full.
By induction, it follows that $((f_{i})_{i \in m}, (h_{k})_{k\in p})$ and $((h_{k})_{k \in p}, (g_{j})_{j\in n})$ belong to $F_{\mathrm{R}}^{\ast}[\Upsilon]^{\natural\mathrm{fl}}$.
Since $F_{\mathrm{R}}^{\ast}[\Upsilon]^{\natural\mathrm{fl}}$ is a generalized congruence, by transitivity, it follows that $((f_{i})_{i \in m}, (g_{j})_{j\in n})$ is in $F_{\mathrm{R}}^{\ast}[\Upsilon]^{\natural\mathrm{fl}}$.

If (2), then there are two pairs 
\[
((F(f'_{i'}))_{i' \in m'}, (F(g'_{j'}))_{j' \in n'})
\text{ and }
((F(f''_{i''}))_{i'' \in m''}, (F(g''_{j''}))_{j'' \in n''})
\] 
in $\Upsilon^{\mathrm{fl}}_{r}$ such that the pairs $(\mathrm{tg}(F(f'_{m'-1})), \mathrm{sc}(F(f''_{0})))$ and $(\mathrm{tg}(F(g'_{n'-1})), \mathrm{tg}(F(g''_{0})))$  belong to $F^{\ast}_{\mathrm{G}}[\Upsilon^{\natural}]^{\mathrm{ob}}$, $(F(f_{i}))_{i \in m} = (F(f'_{i'}))_{i' \in m'} \curlywedge (F(f''_{i''}))_{i'' \in m''}$ and $(F(g_{j}))_{j \in n} = (F(g'_{j'}))_{j' \in n'} \curlywedge (F(g''_{j''}))_{j'' \in n''}$. Note that the functor $F$ is full. 

By induction, it follows that $((f'_{i'})_{i' \in m'}, (g'_{j'})_{j' \in n'})$ and $((f''_{i''})_{i'' \in m''}, (g''_{j''})_{j'' \in n''})$ belong to $F_{\mathrm{R}}^{\ast}[\Upsilon]^{\natural\mathrm{fl}}$. Therefore, since $F_{\mathrm{R}}^{\ast}[\Upsilon]^{\natural}$ is a generalized congruence, 
\begin{equation}
\label{Eq1}
\left(
(f'_{i'})_{i' \in m'} \curlywedge (f''_{i''})_{i'' \in m''}, (g'_{j'})_{j' \in n'} \curlywedge (g''_{j''})_{j'' \in n''}
\right)
\in
F_{\mathrm{R}}^{\ast}[\Upsilon]^{\natural\mathrm{fl}}.
\end{equation}

Moreover, note that $(F(f_{i'}))_{i' \in m'} = (F(f'_{i'}))_{i' \in m'}$, that is, for every $i' \in m'$, $F(f_{i'}) = F(f'_{i'})$. In particular, for every $i' \in m'$, $(F(f_{i'}), F(f'_{i'}))$ belongs to $\Upsilon^{\mathrm{fl}}_{\scalebox{.7}{$\left(\begin{smallmatrix}F(a_{i'})&F(b_{i'})\\F(a_{i'})&F(b_{i'})\end{smallmatrix}\right)$}}$. Therefore, for every $i' \in m'$,
\[
(f_{i'}, f'_{i'})
\in 
(F \times F)^{-1}\left[
\Upsilon^{\mathrm{fl}}_{\scalebox{.7}{$\left(\begin{smallmatrix}F(a_{i'})&F(b_{i'})\\F(a_{i'})&F(b_{i'})\end{smallmatrix}\right)$}}
\right]
\subseteq
F_{\mathrm{R}}^{\ast}[\Upsilon]^{\natural\mathrm{fl}}.
\]
Thus, since $F_{\mathrm{R}}^{\ast}[\Upsilon]^{\natural\mathrm{fl}}$ is a generalized congruence, $((f_{i'})_{i' \in m'}, (f'_{i'})_{i' \in m'})$ is in $F_{\mathrm{R}}^{\ast}[\Upsilon]^{\natural\mathrm{fl}}$.
By a similar argument, it also follows that the pairs
\begin{gather*}
\left( (f_{m'+i''})_{i'' \in m''}, (f''_{i''})_{i'' \in m''} \right), \quad
\left( (g_{j'})_{j' \in n'}, (g'_{j'})_{j' \in n'} \right)\\
\text{and}\quad
\left( (g_{n'+j''})_{j'' \in n''}, (g''_{j''})_{j'' \in n''} \right)
\end{gather*}
also belong to $F_{\mathrm{R}}^{\ast}[\Upsilon]^{\natural\mathrm{fl}}$. Thus, taking into account that $(f_i)_{i \in m} = (f_{i'})_{i' \in m'} \curlywedge (f_{m'+i''})_{i'' \in m''}$ and $(g_j)_{j \in n} = (g_{j'})_{j' \in n'} \curlywedge (g_{n'+j''})_{j'' \in n''}$, it follows that
\begin{equation}
\label{Eq2}
\left(
(f_i)_{i \in m},
(f'_{i'})_{i' \in m'} \curlywedge (f''_{i''})_{i'' \in m''}
\right), 
\left(
(g_j)_{j \in n},
(g'_{j'})_{j' \in n'} \curlywedge (g''_{j''})_{j'' \in n''}
\right)
\in 
F_{\mathrm{R}}^{\ast}[\Upsilon]^{\natural\mathrm{fl}}.
\end{equation}
All in all, from Equations~\eqref{Eq1} and~\eqref{Eq2}, it follows that the pair $((f_i)_{i \in m}, (g_j)_{j \in n})$ belongs to $F_{\mathrm{R}}^{\ast}[\Upsilon]^{\natural\mathrm{fl}}$.

In any case $((f_i)_{i \in m}, (g_j)_{j \in n})$ is in  $F_{\mathrm{R}}^{\ast}[\Upsilon]^{\natural\mathrm{fl}}$.

This completes the inductive step of the induction.

This completes the proof.
\end{proof}

\begin{proposition}
The birestriction of $U_{\mathrm{R}}$ to $\mathsf{RCCat}_{\mathrm{full}}$ and $\mathsf{Cat}_{\mathrm{full}}$, which we continue to denote by $U_{\mathrm{R}}$, is a fibration and $U_{\mathrm{G}}$ a bifibration. Moreover, the restriction of $(\bigcdot)^{\natural}$ to $\mathsf{RCCat}_{\mathrm{full}}$, again denoted by  
$(\bigcdot)^{\natural}$, together with $\mathrm{In}_{\mathsf{Cat}_{\mathrm{full}},\mathsf{Cat}}$, the canonical embedding of $\mathsf{Cat}_{\mathrm{full}}$ into $\mathsf{Cat}$, is a morphism of fibrations from $U_{\mathrm{R}}$ to $U_{\mathrm{G}}$.
\end{proposition}

\begin{proof}
It is straightforward to verify that $U_{\mathrm{R}}$ is a fibration and that $U_{\mathrm{G}}$ is a bifibration. 
Moreover, by Lemma~\ref{commutation}, the pair $((\bigcdot)^{\natural},\mathrm{In}_{\mathsf{Cat}_{\mathrm{full}},\mathsf{Cat}})$ is a morphism of fibrations from $U_{\mathrm{R}}$ to $U_{\mathrm{G}}$.
\end{proof}

\section{Examples of Riguet congruences}\label{S:examples}
We now present several applications where Riguet congruences are used to systematically recover---up to isomorphism or equivalence---fundamental structures, including: the category of finite cardinals $\mathsf{Card}^{A}_{\mathrm{f}}$ from the category associated with the free monoid on $A$; the fundamental group of a path-connected space from its fundamental groupoid; the Serre-Gabriel localization of a Grothendieck category with respect to a localizing subcategory; the action category of a deterministic finite automaton from the category associated to the automaton; and a category of cardinal numbers $\mathsf{Card}_{\mathrm{orb}}$ whose morphisms are defined by orbits under the action of product symmetric groups from the category of relations $\mathsf{Rel}$.

\subsection{A Riguet Congruence on the Category associated to the Free Monoid on a Set}\label{S:freemonoid}

\begin{example}
After introducing the notions and constructions required for this example, we show that there exists an equivalence $\vs^{A}{(\,\bigcdot\,)}$ from $\mathsf{C}(\mathbf{A}^{\star})$, the category associated to the free monoid on a set $A$, to $\mathsf{Set}^{A}_{\mathrm{f}}$, the category of finite $A$-sorted sets. 
Next, considering a quasi-inverse $J$ of the canonical inclusion of the standard skeleton $\mathsf{Card}^{A}_{\mathrm{f}}$ of $\mathsf{Set}^{A}_{\mathrm{f}}$, and defining a suitable Riguet congruence $\Phi$ on $\mathsf{C}(\mathbf{A}^{\star})$, we verify that the functor $F = J \circ \vs^{A}{(\,\bigcdot\,)}$ satisfies conditions~(1) and~(2) of Proposition~\ref{UPQ}. Hence, Proposition~\ref{UPQ} yields a unique functor $F^{\sharp}$ from the skeletal category $\mathsf{C}(\mathbf{A}^{\star})/\Phi$ to $\mathsf{Card}^{A}_{\mathrm{f}}$ such that $F = F^{\sharp} \circ P_{\Phi}$. Moreover, $F^{\sharp}$ is an isomorphism, and $\mathsf{C}(\mathbf{A}^{\star})/\Phi$ is equivalent to the non-skeletal category $\mathsf{C}(\mathbf{A}^{\star})$.
 
Let $A$ be a set. The \emph{free monoid on $A$}, denoted by $\mathbf{A}^{\star}$, is $(A^{\star}, \curlywedge, \lambda)$, where $A^{\star}$, the set of all \emph{words} on $A$, is $\bigcup_{n \in \mathbb{N}} \mathrm{Hom}(n,A)$. For a word $\mathbf{a} \in A^{\star}$ its domain, also called its \emph{length}, is denoted by $\bb{\mathbf{a}}$. Moreover, a word $\mathbf{a}$ is usually denoted as a sequence $(a_{i})_{i \in \bb{\mathbf{a}}}$ in which, for every $i \in \bb{\mathbf{a}}$, $a_{i}$ is $\mathbf{a}(i)$. Furthermore, $\curlywedge$, the \emph{concatenation} of words on $A$, is the binary operation on $A^{\star}$ which sends a pair of words $(\mathbf{a},\mathbf{b})\in A^{\star}$ to the mapping $\mathbf{a} \curlywedge \mathbf{b}$ from $\bb{\mathbf{a}}+\bb{\mathbf{b}}$ to $A$  and $\lambda$, the \emph{empty word} on $A$, is the unique mapping from $\varnothing$ to $A$. For a finite family of words $(\mathbf{a}_{j})_{j \in n}$ in $A^{\star}$, we let  $\bigcurlywedge_{j \in n} \mathbf{a}_{j}$  stand for the concatenation of the words in the family. We will denote by $\eta_{A}$ the mapping from $A$ to $A^{\star}$ that sends $a \in A$ to $(a) \in A^{\star}$. 
The ordered pair $(\mathbf{A}^{\star}, \eta_{A})$ is a universal morphism from $A$ to the forgetful functor from the category $\mathsf{Mon}$, of monoids, to $\mathsf{Set}$. For a monoid $\mathbf{M}=(M,\cdot, 1)$ and a mapping $f \colon A \mor M$ we will denote by $f^{\sharp}$ the unique monoid homomorphism from $\mathbf{A}^{\star}$ to $\mathbf{M}$ such that $f^{\sharp} \circ \eta_{A}=f$.
Let $\mathbf{a}$ be a word in $A^{\star}$ and $a \in A$. Then we let $\bb{\mathbf{a}}_{a}$ stand for the number of occurrences of $a$ in $\mathbf{a}$, i.e., for the cardinal of $\mathbf{a}^{-1}[\{a\}]$. 
For $\mathbf{a} \in A^{\star}$ and $l \in \bb{\mathbf{a}}$, we will denote by $\mathbf{a}^{0,l}$ the \emph{initial segment} of $\mathbf{a}$ beginning at position $0$ and ending at position $l$.

We let $\mathsf{C}(\mathbf{A}^{\star})$ stand for the category whose objects are the words on $A$ and whose morphisms from $\mathbf{a}$ to $\mathbf{b}$ are the ordered triples 
$(\mathbf{a},\varphi,\mathbf{b})$, written simply as $\varphi$ for short and denoted by $\varphi\colon \mathbf{a}\mor \mathbf{b}$, in which $\varphi$ is a mapping from $\bb{\mathbf{a}}$ to $\bb{\mathbf{b}}$ such that $\mathbf{a}=\mathbf{b} \circ \varphi$. For a word $\mathbf{a}$, the identity morphism at $\mathbf{a}$, denoted by  $\mathrm{id}_{\mathbf{a}}$, is $\mathrm{id}_{\bb{\mathbf{a}}}$. Finally, for two morphisms $\varphi \colon \mathbf{a} \mor \mathbf{b}$ and $\psi \colon \mathbf{b} \mor \mathbf{c}$ in $\mathsf{C}(\mathbf{A}^{\star})$, its composition in $\mathsf{C}(\mathbf{A}^{\star})$ is $\psi \circ \varphi$, the standard composition of $\varphi$ and $\psi$. The category $\mathsf{C}(\mathbf{A}^{\star})$ has finite coproducts. In fact $\lambda$ is initial and, for every $\mathbf{a}$, $\mathbf{b}\in A^{\star}$, 
$(\mathbf{a}\curlywedge \mathbf{b},(\iota_{\mathbf{a}},\iota_{\mathbf{b}}))$, where 
$\iota_{\mathbf{a}}$ and $\iota_{\mathbf{b}}$ are the canonical embedding of 
$\bb{\mathbf{a}}$ and $\bb{\mathbf{b}}$ into $\bb{\mathbf{a}\curlywedge \mathbf{b}}$, respectively, is the coproduct of $\mathbf{a}$ and $\mathbf{b}$. Moreover, $\mathsf{C}(\mathbf{A}^{\star})$ is a concrete realization of the free category with finite coproducts on the discrete category on the set $A$---a particular case of the well-known free coproduct completion of a category---and it is a strict monoidal category. Note that for a set $A$ such that $\mathrm{card}(A)\geq 2$, the category $\mathsf{C}(\mathbf{A}^{\star})$ is not skeletal.

Let $\mathsf{Set}^{A}_{\mathrm{f}}$ be the category whose objects are the finite $A$-sorted sets, i.e., the $A$-sorted sets $X = (X_{a})_{a \in A}$ such that $\mathrm{supp}_{A}(X)$ is finite and, for every $a \in \mathrm{supp}_{A}(X)$, $X_{a}$ is finite, where $\mathrm{supp}_{A}(X)$ is $\{a\in A\mid X_{a}\neq \varnothing\}$. And whose morphisms are the $A$-sorted  mappings $f = (f_{a})_{a \in A}$ from $X$ to $Y$, where, for every $a \in A$, $f_{a}$ is a mapping from $X_{a}$ to $Y_{a}$. Note that $\mathsf{Set}^{A}_{\mathrm{f}}$ is equivalent to the comma category $(\mathsf{Set}_{\mathrm{f}}{\downarrow}A)$. Then the categories $\mathsf{C}(\mathbf{A}^{\star})$ and  $\mathsf{Set}^{A}_{\mathrm{f}}$, are equivalent. In fact, let $\vs^{A}{(\,\bigcdot\,)}$ be the assignment from $\mathsf{C}(\mathbf{A}^{\star})$ to $\mathsf{Set}^{A}_{\mathrm{f}}$ defined as follows:
\begin{enumerate}
\item for every object $\mathbf{a}$ of $\mathsf{C}(\mathbf{A}^{\star})$, $\vs^{A}\mathbf{a}$ is the $A$-sorted set defined, for every $a \in A$, as $(\vs^{A}{\mathbf{a}})_{a}=\mathbf{a}^{-1}[\{a\}]$, and
\item for every morphism $\varphi \colon \mathbf{a} \mor \mathbf{b}$ of $\mathsf{C}(\mathbf{A}^{\star})$, $\vs^{A}{\varphi}$ is the $A$-sorted mapping from $\vs^{A}{\mathbf{a}}$ to $\vs^{A}{\mathbf{b}}$ defined, for every $a \in A$ and every $i \in (\vs^{A}{\mathbf{a}})_{a}$, as $(\vs^{A}{\varphi})_{a}(i)= \varphi(i)$.
\end{enumerate}
Thus defined,  $\vs^{A}{(\,\bigcdot\,)}$ is a functor from $\mathsf{C}(\mathbf{A}^{\star})$ to $\mathsf{Set}^{A}_{\mathrm{f}}$. It is straightforward to verify functoriality, and we leave the details to the reader. 
Hence, it suffices to prove that  $\vs^{A}{(\,\bigcdot\,)}$ is an equivalence from $\mathsf{C}(\mathbf{A}^{\star})$ to $\mathsf{Set}^{A}_{\mathrm{f}}$.

\textsf{Fullness.}
Let $\mathbf{a}$ and $\mathbf{b}$ be objects of $\mathsf{C}(\mathbf{A}^{\star})$ and $f \colon {\vs^{A}{\mathbf{a}}} \mor {\vs^{A}{\mathbf{b}}}$ an $A$-sorted mapping. We want to prove that there is a morphism $\varphi$ from $\mathbf{a}$ to $\mathbf{b}$ such that $\vs^{A}{\varphi}=f$, i.e., a mapping $\varphi \colon \bb{\mathbf{a}} \mor \bb{\mathbf{b}}$ such that $\mathbf{b} \circ \varphi = \mathbf{a}$ and $\vs^{A}{\varphi}=f$. It suffices to take as $\varphi$ precisely $[\mathrm{in}_{a} \circ f_{a}]_{a \in A}$, i.e., the unique mapping from $\bb{\mathbf{a}}$ to $\bb{\mathbf{b}}$ such that $[\mathrm{in}_{a} \circ f_{a}]_{a \in A} \circ \mathrm{in}_{A}=\mathrm{in}_{A} \circ f_{a}$, where, for every $a \in A$, $f_{a}$ is the $a$-th component of $f$ and where, by abuse of notation, we use $\mathrm{in}_{a}$ both for the canonical inclusion of $(\vs^{A}{\mathbf{a}})_{a}$ into $\bb{\mathbf{a}}$ and for the canonical inclusion of $(\vs^{A}{\mathbf{b}})_{a}$ into $\bb{\mathbf{b}}$.

\textsf{Faithfulness.}
Let $\mathbf{a}$ and $\mathbf{b}$ be objects of $\mathsf{C}(\mathbf{A}^{\star})$ and let $\varphi, \psi \colon \mathbf{a} \mor \mathbf{b}$ be morphisms of $\mathsf{C}(\mathbf{A}^{\star})$ such that $\vs^{A}{\varphi} = \vs^{A}{\psi}$. We want to prove that $\varphi = \psi$. 
Let $i \in \bb{\mathbf{a}}$. By definition of $\vs^{A}{\varphi}$ and $\vs^{A}{\psi}$, we have that
\[\varphi (i) = (\vs^{A}{\varphi})_{a_{i}}(i) = (\vs^{A}{\psi})_{a_{i}}(i) = \psi (i).\]
Thus, $\varphi = \psi$.

\textsf{Essential surjectivity.}
Let $X$ be a finite $A$-sorted set. Hence $\mathrm{supp}_{A}(X)$ is finite and, for every $a \in \mathrm{supp}_{A}(X)$, $X_{a}$ is finite. Let $\mathrm{card}(\mathrm{supp}_{A}(X))=m$, let $a_{0}, \ldots, a_{m-1}$ be the distinct elements of $\mathrm{supp}_{A}(X)$ numbered in an arbitrary manner by the natural numbers in $m$ and, for every $j \in m$, let $n_{j}=\mathrm{card}(X_{a_{j}})$. Then, for the word $\bigcurlywedge_{j \in m} (a_{j})^{n_{j}}$, where, for every $j \in m$, $(a_{j})^{n_{j}}$ is the word of length $n_{j}$ whose image is $\{a_{j}\}$, i.e., the constant mapping from $n_{j}$ to $A$ with value $a_{j}$, we have that $X$ and $\vs^{A}{(\bigcurlywedge_{j \in m} (a_{j})^{n_{j}})}$ are isomorphic. 

Therefore $\vs^{A}{(\,\bigcdot\,)} \colon \mathsf{C}(\mathbf{A}^{\star}) \mor \mathsf{Set}^{A}_{\mathrm{f}}$ is an equivalence.
 
Let $\mathsf{Card}^{A}_{\mathrm{f}}$ be the full subcategory of $\mathsf{Set}^{A}_{\mathrm{f}}$ determined by the finite $A$-sorted cardinals, i.e., by the set \[\mathbb{N}^{(A)}=\{ (n_{a})_{a \in A} \in \mathrm{Hom}(A, \mathbb{N}) \mid \mathrm{card}(\{a \in A \mid n_{a} \neq 0\}) < \aleph_{0}\},\]
and $J$ the quasi-inverse of the canonical inclusion of $\mathsf{Card}^{A}_{\mathrm{f}}$ into  
$\mathsf{Set}^{A}_{\mathrm{f}}$. We recall that, for every finite $A$-sorted set $X$, $J(X) = (\mathrm{card}(X_{a}))_{a\in A}$ and that, for every $A$-sorted mapping $f = (f_{a})_{a \in A}$ from $X$ to $Y$, by choosing, for every finite $A$-sorted set $X$, a bijection $\theta_{X}$ from $X$ to $J(X)$, we have that $J(f) = \theta_{Y}\circ f\circ \theta_{X}^{-1}$.

We next define a congruence on $\mathbf{A}^{\star}$ that will serve as the object component of a Riguet congruence on $\mathsf{C}(\mathbf{A}^{\star})$, such that the quotient category of $\mathsf{C}(\mathbf{A}^{\star})$ by it is isomorphic to $\mathsf{Card}^{A}_{\mathrm{f}}$. Let $\equiv^{A}$ be the binary relation on $A^{\star}$ defined as follows: two words $\mathbf{a}$ and $\mathbf{b}$ in $A^{\star}$ satisfy $\mathbf{a} \equiv^{A} \mathbf{b}$ if, and only if, $\bb{\mathbf{a}} = \bb{\mathbf{b}}$ and there exists a permutation $\sigma$ of $\bb{\mathbf{a}}$ such that $\mathbf{a} = \mathbf{b} \circ \sigma$ (i.e., $a_{i} = b_{\sigma(i)}$ for every $i \in \bb{\mathbf{a}}$). Then $\equiv^{A}$ is a congruence on $\mathbf{A}^{\star}$, and $\mathbf{A}^{\star}/{\equiv^{A}}$ is the free abelian monoid on $A$. Note that $\mathbf{a} \equiv^{A} \mathbf{b}$ if, and only if, for every $a \in A$, $\bb{\mathbf{a}}_{a} = \bb{\mathbf{b}}_{a}$.

Next, we assign to every pair $(\mathbf{a},\mathbf{b})$ of $\equiv^{A}$-equivalent words the so-called \emph{canonical permutation} for $(\mathbf{a},\mathbf{b})$. These canonical permutations will be a key tool in our definition of the Riguet congruence on $\mathbf{A}^{\star}$.

We begin by defining the notions of the occurrence of an index in a word and of the position of a letter within a word. 
Let $\mathbf{a}$ be an element of $A^{\star}$. Then we let $\mathrm{occ}_{\mathbf{a}}$ stand for the mapping from $\bb{\mathbf{a}}$ to $\coprod_{i\in \bb{\mathbf{a}}}\bb{\mathbf{a}}_{\mathbf{a}(i)}$ defined, for every $i\in \bb{\mathbf{a}}$, as
\[
\textstyle
\mathrm{occ}_{\mathbf{a}}(i) = \left(\bigcup ((i_{j})_{j\in\bb{\mathbf{a}}_{\mathbf{a}(i)}})^{-1}[\{i\}]\right)\times\{i\},
\]
where $(i_{j})_{j\in\bb{\mathbf{a}}_{\mathbf{a}(i)}}$ is the enumeration in ascending order of the occurrences of $\mathbf{a}(i)$ in $\mathbf{a}$. Thus $\mathrm{occ}_{\mathbf{a}}(i)$ is the pair $(j,i)$, where $j$ is the unique element of $\bb{\mathbf{a}}_{\mathbf{a}(i)}$ such that $i_{j} = i$. We call $\mathrm{occ}_{\mathbf{a}}(i)$ the \emph{occurrence} of the index $i$ in the word $\mathbf{a}$. By abuse of notation we will write $j$ instead of $(j,i)$.
Let $(\mathbf{a},a)$ be an element of $\coprod_{\mathbf{a}\in A^{\star}}\mathrm{Im}(\mathbf{a}) = \bigcup_{\mathbf{a}\in A^{\star}}(\{\mathbf{a\}}\times \mathrm{Im}(\mathbf{a}))$. Then we let $\mathrm{pos}_{\mathbf{a},a}$ stand for the mapping from 
$\bb{\mathbf{a}}_{a}$ to $\bb{\mathbf{a}}$ that sends $j\in \bb{\mathbf{a}}_{a}$ to the unique $i\in \bb{\mathbf{a}}$ such that (1) $\mathbf{a}(i) = a$ and 
(2) $\bb{\mathbf{a}^{0,i}}_{a} = j+1$. We call $\mathrm{pos}_{\mathbf{a},a}(j)$ the \emph{position} of the $j$-th occurrence of $a$ in $\mathbf{a}$. 

We define $\overline{\mathrm{occ}}_{\mathbf{a}, a}$ as the birestriction of the global mapping 
$\mathrm{occ}_{\mathbf{a}}$ to $\mathbf{a}^{-1}[\{a\}]$ and $\bb{\mathbf{a}}_{a}$, its corresponding cardinal in the disjoint union:
\[
\overline{\mathrm{occ}}_{\mathbf{a}, a} \colon \mathbf{a}^{-1}[\{a\}] \mor \bb{\mathbf{a}}_{a}.
\]
Note that $\overline{\mathrm{occ}}_{\mathbf{a}, a}(i) = j$, where $j$ is the first coordinate of the pair $\mathrm{occ}_{\mathbf{a}}(i)$. Moreover, we define $\overline{\mathrm{pos}}_{\mathbf{a}, a}$ as the corestriction of the position mapping to the fiber of $a$:
\[
\overline{\mathrm{pos}}_{\mathbf{a}, a} \colon \bb{\mathbf{a}}_{a} \mor \mathbf{a}^{-1}[\{a\}].
\]
Note that $\overline{\mathrm{pos}}_{\mathbf{a}, a}(j) = i$, where $i$ is the unique index in the word such that $\mathbf{a}(i) = a$ and $\overline{\mathrm{occ}}_{\mathbf{a}, a}(i) = j$.

By construction, since we are using the same ordering (the natural order of $\mathbb{N}$) for both the indices $i$ in the word and the count $j$ in the cardinal, these mappings are strictly monotonic and cover the same number of elements. Thus, for every $i \in  \mathbf{a}^{-1}[\{a\}]$, $\overline{\mathrm{pos}}_{\mathbf{a}, a}(\overline{\mathrm{occ}}_{\mathbf{a}, a}(i)) = i$ and, for every $j \in \bb{\mathbf{a}}_{a}$, $\overline{\mathrm{occ}}_{\mathbf{a}, a}(\overline{\mathrm{pos}}_{\mathbf{a}, a}(j)) = j$.

Note that the mappings $\mathrm{occ}_{\mathbf{a}}$ and $\mathrm{pos}_{\mathbf{a},a}$ are constructive; this will ensure that the canonical permutations defined below are effectively computable.

We next assign to every pair $(\mathbf{a},\mathbf{b})$ of $\equiv^{A}$-equivalent words the so-called canonical permutation for $(\mathbf{a},\mathbf{b})$. 
Let $\mathbf{a}$ and $\mathbf{b}$ be words in $A^{\star}$ such that $\mathbf{a} \equiv^{A} \mathbf{b}$, hence, in particular, we have that $\bb{\mathbf{a}}=\bb{\mathbf{b}}$. Then the permutation $\sigma_{\mathbf{a},\mathbf{b}}$ of 
$\bb{\mathbf{a}}$, defined as:  
\[
\sigma_{\mathbf{a},\mathbf{b}}
\left\lbrace
\begin{array}{ccl}
\bb{\mathbf{a}} & \mor & \bb{\mathbf{b}}\\
i & \longmapsto & \mathrm{pos}_{\mathbf{b}, \mathbf{a}(i)}(\mathrm{occ}_{\mathbf{a}}(i))
\end{array}
\right.
\]
and called the \emph{canonical permutation} for the pair $(\mathbf{a},\mathbf{b})$, is such that $\mathbf{a}=\mathbf{b} \circ \sigma_{\mathbf{a},\mathbf{b}}$. These canonical permutations have the following properties:  Let $\mathbf{a}, \mathbf{b}, \mathbf{c}$ be words in $A^{\star}$ such that $\mathbf{a} \equiv^{A} \mathbf{b}$ and $\mathbf{b} \equiv^{A} \mathbf{c}$. Then
(1) $\sigma_{\mathbf{a},\mathbf{a}}=\mathrm{id}_{\bb{\mathbf{a}}}$,
(2) $\sigma_{\mathbf{b},\mathbf{c}} \circ \sigma_{\mathbf{a},\mathbf{b}}= \sigma_{\mathbf{a},\mathbf{c}}$, and
(3) $\sigma_{\mathbf{a},\mathbf{b}}^{-1}=\sigma_{\mathbf{b},\mathbf{a}}$.

We define a Riguet congruence $\Phi$ on $\mathsf{C}(\mathbf{A}^{\star})$ as follows: 
\begin{enumerate}
\item \textsf{On Objects:} We set $\Phi^{\mathrm{ob}} = {\equiv^{A}}$. 
\item \textsf{On Arrows:} For every $\mathbf{a}, \mathbf{a}', \mathbf{b}, \mathbf{b}' \in A^{\star}$ such that $\mathbf{a}\equiv^{A}\mathbf{a}'$ and $\mathbf{b} \equiv^{A} \mathbf{b}'$, we let ${\equiv}^{A}_{\scalebox{0.7}{$\left(\begin{smallmatrix}\mathbf{a}&\mathbf{b}\\ \mathbf{a}'&\mathbf{b}'\end{smallmatrix}\right)$}}$ stand for the subset of $\mathrm{Hom}_{\mathsf{C}(\mathbf{A}^{\star})}(\mathbf{a},\mathbf{b}) \times \mathrm{Hom}_{\mathsf{C}(\mathbf{A}^{\star})}(\mathbf{a}',\mathbf{b}')$ defined as follows: For every $\varphi \colon \mathbf{a} \mor \mathbf{b}$ and $\varphi' \colon \mathbf{a}' \mor \mathbf{b}'$, we have that $(\varphi, \varphi') \in {\equiv}^{A}_{\scalebox{0.7}{$\left(\begin{smallmatrix}\mathbf{a}&\mathbf{b}\\ \mathbf{a}'&\mathbf{b}'\end{smallmatrix}\right)$}}$ if, and only if, for the canonical permutations $\sigma_{\mathbf{a}, \mathbf{a}'}$ for $(\mathbf{a}, \mathbf{a}')$ and $\sigma_{\mathbf{b},\mathbf{b}'}$ for $(\mathbf{b}, \mathbf{b}')$ (which are such that $\mathbf{a}=\mathbf{a}'\circ \sigma_{\mathbf{a},\mathbf{a}'}$ and $\mathbf{b}=\mathbf{b}'\circ \sigma_{\mathbf{b},\mathbf{b}'}$, respectively), we have that 
\[
\varphi' \circ \sigma_{\mathbf{a},\mathbf{a}'} = \sigma_{\mathbf{b},\mathbf{b}'} \circ \varphi.
\]
\end{enumerate}
 
By definition, $\Phi^{\mathrm{ob}} = {\equiv^{A}}$ and $\equiv^{A}$ is a congruence on $\mathbf{A}^{\star}$, hence $\Phi^{\mathrm{ob}}$ is an equivalence relation on $\mathrm{Ob}(\mathsf{C}(\mathbf{A}^{\star}))$. The proof of conditions $(\mathrm{a})$ through $(\mathrm{f})$ is immediate and is omitted. 

The quotient category $\mathsf{C}(\mathbf{A}^{\star})/{\equiv^{A}}$ is a skeletal strict symmetric monoidal category.

To prove that $\mathsf{C}(\mathbf{A}^{\star})/{\equiv^{A}}$ is isomorphic to $\mathsf{Card}^{A}_{\mathrm{f}}$, 
which will yield the equivalence of the former with both $\mathsf{Set}^{A}_{\mathrm{f}}$ and $\mathsf{C}(\mathbf{A}^{\star})$, we begin by proving that the functor $F = J\circ\vs^{A}{(\,\bigcdot\,)}$ from $\mathsf{C}(\mathbf{A}^{\star})$ to $\mathsf{Card}^{A}_{\mathrm{f}}$ verifies the conditions in Proposition~\ref{UPQ}. Before doing it we agree that for the finite $A$-sorted sets of type $(\mathbf{a}^{-1}[\{a\}])_{a\in A}$, the bijection $\theta_{(\mathbf{a}^{-1}[\{a\}])_{a\in A}}$, denoted by 
$\theta_{\mathbf{a}}$, for short, from $(\mathbf{a}^{-1}[\{a\}])_{a\in A}$ to 
$J((\mathbf{a}^{-1}[\{a\}])_{a\in A}) = (\bb{\mathbf{a}}_{a})_{a|\in A}$, is such that, for every $a\in A$,  $\theta_{\mathbf{a},a} = \overline{\mathrm{occ}}_{\mathbf{a}, a}$. Its inverse $\theta^{-1}_{\mathbf{a}}$ is $(\overline{\mathrm{pos}}_{\mathbf{a}, a})_{a\in A}$.
Then, given a morphism $\varphi \colon \mathbf{a} \mor \mathbf{b}$ in $\mathsf{C}(\mathbf{A}^{\star})$, the action of $F$ on the $a$-component is defined by:
\[
F(\varphi)_{a} = \theta_{\mathbf{b}, a} \circ \varphi|_{\mathbf{a}^{-1}[\{a\}]} \circ \theta_{\mathbf{a}, a}^{-1}.
\]
Thus, we have that:
\[
F(\varphi)_{a} = \overline{\mathrm{occ}}_{\mathbf{b}, a} \circ \varphi|_{\mathbf{a}^{-1}[\{a\}]} \circ \overline{\mathrm{pos}}_{\mathbf{a}, a}.
\]
Note that we let $\varphi|_{\mathbf{a}^{-1}[\{a\}]}$ stand for the birestriction of $\varphi$ to $\mathbf{a}^{-1}[\{a\}]$ and $\mathbf{b}^{-1}[\{a\}]$.

We are now in position to verify the conditions in Proposition~\ref{UPQ}.

(1) Let $\mathbf{a}, \mathbf{a}' \in \mathrm{Ob}(\mathsf{C}(\mathbf{A}^{\star}))$ be such that $(\mathbf{a}, \mathbf{a}') \in \Phi^{\mathrm{ob}}$. By definition, this means that $\mathbf{a} \equiv^{A} \mathbf{a}'$, i.e., that $\bb{\mathbf{a}} = \bb{\mathbf{a}'}$ and there exists a permutation $\sigma$ of $\bb{\mathbf{a}}$ such that $\mathbf{a} = \mathbf{a}' \circ \sigma$, which in turn is equivalent to: for every $a \in A$, $\bb{\mathbf{a}}_{a} = \bb{\mathbf{a}'}_{a}$. 
Consequently, the $A$-sorted sets $\vs^{A}\mathbf{a}$ and $\vs^{A}\mathbf{a}'$ have the same cardinalities for each component. Thus $F(\mathbf{a}) = F(\mathbf{a}')$.

(2) Let $\varphi \colon \mathbf{a} \mor \mathbf{b}$ and $\varphi' \colon \mathbf{a}' \mor \mathbf{b}'$ be morphisms in $\mathsf{C}(\mathbf{A}^{\star})$ such that $(\varphi, \varphi') \in {\equiv}^{A}_{\scalebox{0.7}{$\left(\begin{smallmatrix}\mathbf{a}&\mathbf{b}\\ \mathbf{a}'&\mathbf{b}'\end{smallmatrix}\right)$}}$. This condition is equivalent to: 
\[
\varphi' \circ \sigma_{\mathbf{a}, \mathbf{a}'} = \sigma_{\mathbf{b}, \mathbf{b}'} \circ \varphi.
\]
We wish to show that $F(\varphi)_a = F(\varphi')_a$ for every $a \in A$. First, we recall the relationship between the canonical permutations and our fiber-wise bijections. For any $a \in A$, the restriction of $\sigma_{\mathbf{a}, \mathbf{a}'}$ to the fiber $\mathbf{a}^{-1}[\{a\}]$, denoted by $\sigma_{\mathbf{a}, \mathbf{a}', a}$, is:
\[
\sigma_{\mathbf{a}, \mathbf{a}', a} = \overline{\mathrm{pos}}_{\mathbf{a}', a} \circ \overline{\mathrm{occ}}_{\mathbf{a}, a}
\]
By applying $\overline{\mathrm{pos}}_{\mathbf{a}, a}$ to both sides and using the identity 
$\overline{\mathrm{occ}}_{\mathbf{a}, a} \circ \overline{\mathrm{pos}}_{\mathbf{a}, a} = \mathrm{id}_{\bb{\mathbf{a}}_a}$, we obtain:
\[
\sigma_{\mathbf{a}, \mathbf{a}', a} \circ \overline{\mathrm{pos}}_{\mathbf{a}, a} = \overline{\mathrm{pos}}_{\mathbf{a}', a}.
\]
Similarly, for the codomain words $\mathbf{b}$ and $\mathbf{b}'$:
\[
\overline{\mathrm{occ}}_{\mathbf{b}', a} \circ \sigma_{\mathbf{b}, \mathbf{b}', a} = \overline{\mathrm{occ}}_{\mathbf{b}, a}.
\]
Starting with the definition of the functor $F$ on the morphism $\varphi'$ for the component $a$:
\[
F(\varphi')_a = \overline{\mathrm{occ}}_{\mathbf{b}', a} \circ \varphi'|_{\mathbf{a}'^{-1}[\{a\}]} \circ \overline{\mathrm{pos}}_{\mathbf{a}', a}.
\]
Substitute the identity $\overline{\mathrm{pos}}_{\mathbf{a}', a} = \sigma_{\mathbf{a}, \mathbf{a}', a} \circ \overline{\mathrm{pos}}_{\mathbf{a}, a}$:
\[
F(\varphi')_a = \overline{\mathrm{occ}}_{\mathbf{b}', a} \circ \varphi'|_{\mathbf{a}'^{-1}[\{a\}]} \circ \sigma_{\mathbf{a}, \mathbf{a}', a} \circ \overline{\mathrm{pos}}_{\mathbf{a}, a}.
\]
Since $\varphi'|_{\mathbf{a}'^{-1}[\{a\}]} \circ \sigma_{\mathbf{a}, \mathbf{a}', a} = (\varphi' \circ \sigma_{\mathbf{a}, \mathbf{a}'})|_{\mathbf{a}^{-1}[\{a\}]}$, we can apply the congruence condition $\varphi' \circ \sigma_{\mathbf{a}, \mathbf{a}'} = \sigma_{\mathbf{b}, \mathbf{b}'} \circ \varphi$:
\[
F(\varphi')_a = \overline{\mathrm{occ}}_{\mathbf{b}', a} \circ (\sigma_{\mathbf{b}, \mathbf{b}'} \circ \varphi)|_{\mathbf{a}^{-1}[\{a\}]} \circ \overline{\mathrm{pos}}_{\mathbf{a}, a}.
\]
Decomposing the restricted composition:
\[
F(\varphi')_a = \overline{\mathrm{occ}}_{\mathbf{b}', a} \circ \sigma_{\mathbf{b}, \mathbf{b}', a} \circ \varphi|_{\mathbf{a}^{-1}[\{a\}]} \circ \overline{\mathrm{pos}}_{\mathbf{a}, a}.
\]
Finally, applying the second identity $\overline{\mathrm{occ}}_{\mathbf{b}', a} \circ \sigma_{\mathbf{b}, \mathbf{b}', a} = \overline{\mathrm{occ}}_{\mathbf{b}, a}$:
\[
F(\varphi')_a = \overline{\mathrm{occ}}_{\mathbf{b}, a} \circ \varphi|_{\mathbf{a}^{-1}[\{a\}]} \circ \overline{\mathrm{pos}}_{\mathbf{a}, a}
\]
Comparing this to the definition of $F(\varphi)_a$, we see that:
\[
F(\varphi')_a = F(\varphi)_a.
\]
Thus, $F(\varphi) = F(\varphi')$, which completes the verification.

Hence, by Proposition~\ref{UPQ}, there exists a unique functor $F^{\sharp}$ from $\mathsf{C}(\mathbf{A}^{\star})/{\equiv^{A}}$ to $\mathsf{Card}^{A}_{\mathrm{f}}$ such that $F^{\sharp}\circ P_{\equiv^{A}} = F$, where $P_{\equiv^{A}}$ is the canonical projection from $ \mathsf{C}(\mathbf{A}^{\star})$ to $\mathsf{C}(\mathbf{A}^{\star})/{\equiv^{A}}$.

To show that $F^{\sharp}$ is an isomorphism, we verify that it is bijective on objects, faithful, and full.

\textsf{Bijectivity on objects.}
The bijectivity on objects of $F^{\sharp}$ follows from the essential surjectivity of $\vs^{A}$ and the fact that $\mathbf{a} \equiv^{A} \mathbf{b}$ if, and only if, they have the same $A$-sorted cardinalities, i.e., for every $a\in A$, $\bb{\mathbf{a}}_{a} = \bb{\mathbf{b}}_{a}$. 

\textsf{Faithfulness.}
We want to prove that, for every objects $[\mathbf{a}]_{\Phi^{\mathrm{ob}}}$, $[\mathbf{b}]_{\Phi^{\mathrm{ob}}}$ in $\mathsf{C}(\mathbf{A}^{\star})/{\equiv^{A}}$ and every $[\varphi]_{\Phi^{\mathrm{fl}}}, [\psi]_{\Phi^{\mathrm{fl}}} \in \mathrm{Hom}_{\mathsf{C}(\mathbf{A}^{\star})/\equiv^{A}}([\mathbf{a}]_{\Phi^{\mathrm{ob}}}, [\mathbf{b}]_{\Phi^{\mathrm{ob}}})$, if $F^{\sharp}([\varphi]_{\Phi^{\mathrm{fl}}}) = F^{\sharp}([\psi]_{\Phi^{\mathrm{fl}}})
$, then $[\varphi]_{\Phi^{\mathrm{fl}}} = [\psi]_{\Phi^{\mathrm{fl}}}$. Suppose:
\[
F^{\sharp}([\varphi]_{\Phi^{\mathrm{fl}}}) = F^{\sharp}([\psi]_{\Phi^{\mathrm{fl}}}).
\]
By the definition of $F^{\sharp}$, this means that $F(\varphi) = F(\psi)$ for any representatives $\varphi \colon \mathbf{a}_{0} \mor \mathbf{b}_{0}$ and $\psi \colon \mathbf{a}_{1} \mor \mathbf{b}_{1}$ where $\mathbf{a}_{0} \equiv^{A} \mathbf{a}_{1}$ and $\mathbf{b}_{0} \equiv^{A} \mathbf{b}_{1}$. For every $a \in A$, the components are equal:
\[
\overline{\mathrm{occ}}_{\mathbf{b}_{0}, a} \circ \varphi|_{\mathbf{a}_{0}^{-1}[\{a\}]} \circ \overline{\mathrm{pos}}_{\mathbf{a}_{0}, a} = \overline{\mathrm{occ}}_{\mathbf{b}_{1}, a} \circ \psi|_{\mathbf{a}_{1}^{-1}[\{a\}]} \circ \overline{\mathrm{pos}}_{\mathbf{a}_{1}, a}.
\]
To show that $[\varphi]_{\Phi^{\mathrm{fl}}} = [\psi]_{\Phi^{\mathrm{fl}}}$, we must prove that  $(\varphi, \psi) \in \equiv^{A}_{\scalebox{0.7}{$\left(\begin{smallmatrix}\mathbf{a}_{0}&\mathbf{b}_{0}\\\mathbf{a}_{1}& \mathbf{b}_{1}\end{smallmatrix}\right)$}}$, which requires:
\[
\psi \circ \sigma_{\mathbf{a}_{0}, \mathbf{a}_{1}} = \sigma_{\mathbf{b}_{0}, \mathbf{b}_{1}} \circ \varphi.
\]
Evaluating the $a$-fiber component of the left side:
\[
\psi|_{\mathbf{a}_{1}^{-1}[\{a\}]} \circ \sigma_{\mathbf{a}_{0}, \mathbf{a}_{1}, a} = \psi|_{\mathbf{a}_{1}^{-1}[\{a\}]} \circ \left( \overline{\mathrm{pos}}_{\mathbf{a}_{1}, a} \circ \overline{\mathrm{occ}}_{\mathbf{a}_{0}, a} \right).
\]
From our initial hypothesis $F(\varphi) = F(\psi)$, we can isolate $\psi|_{\mathbf{a}_{1}^{-1}[\{a\}]}$ by composing with bijections:
\[
\psi|_{\mathbf{a}_{1}^{-1}[\{a\}]} = \overline{\mathrm{pos}}_{\mathbf{b}_{1}, a} \circ \overline{\mathrm{occ}}_{\mathbf{b}_{0}, a} \circ \varphi|_{\mathbf{a}_{0}^{-1}[\{a\}]} \circ \overline{\mathrm{pos}}_{\mathbf{a}_{0}, a} \circ \overline{\mathrm{occ}}_{\mathbf{a}_{1}, a}.
\]
Note that $\overline{\mathrm{pos}}_{\mathbf{b}_{1}, a} \circ \overline{\mathrm{occ}}_{\mathbf{b}_{0}, a} = \sigma_{\mathbf{b}_{0}, \mathbf{b}_{1}, a}$. Substituting this back:
\[
\psi|_{\mathbf{a}_{1}^{-1}[\{a\}]} = \sigma_{\mathbf{b}_{0}, \mathbf{b}_{1}, a} \circ \varphi|_{\mathbf{a}_{0}^{-1}[\{a\}]} \circ \sigma_{\mathbf{a}_{1}, \mathbf{a}_{0}, a}.
\]
Since $\sigma_{\mathbf{a}_{1}, \mathbf{a}_{0}, a} = \sigma_{\mathbf{a}_{0}, \mathbf{a}_{1}, a}^{-1}$, we have:
\[
\psi|_{\mathbf{a}_{1}^{-1}[\{a\}]} \circ \sigma_{\mathbf{a}_{0}, \mathbf{a}_{1}, a} = \sigma_{\mathbf{b}_{0}, \mathbf{b}_{1}, a} \circ \varphi|_{\mathbf{a}_{0}^{-1}[\{a\}]}.
\]
This equality holds for every $a \in A$. Summing over the disjoint fibers, we get $\psi \circ \sigma_{\mathbf{a}_{0}, \mathbf{a}_{1}} = \sigma_{\mathbf{b}_{0}, \mathbf{b}_{1}} \circ \varphi$. Thus, $(\varphi, \psi) \in {\equiv}^{A}_{\scalebox{0.7}{$\left(\begin{smallmatrix}\mathbf{a}_{0}&\mathbf{b}_{0}\\ \mathbf{a}_{1}&\mathbf{b}_{1}\end{smallmatrix}\right)$}}$, hence $[\varphi]_{\Phi^{\mathrm{fl}}} = [\psi]_{\Phi^{\mathrm{fl}}}$ and, consequently,  $F^{\sharp}$ is faithful.

\textsf{Fullness.}
To prove the fullness of the induced functor $F^{\sharp} \colon \mathsf{C}(\mathbf{A}^{\star})/{\equiv^{A}} \mor\mathsf{Card}^{A}_{\mathrm{f}}$, we must show that for any morphism $g$ in the target category, there exists a corresponding morphism in the quotient category that $F^{\sharp}$ maps to $g$. Let $[\mathbf{a}]_{\Phi^{\mathrm{ob}}}$ and $[\mathbf{b}]_{\Phi^{\mathrm{ob}}}$ be objects in $\mathsf{C}(\mathbf{A}^{\star})/{\equiv^{A}}$, and let $g = (g_{a})_{a \in A} \colon F^{\sharp}([\mathbf{a}]_{\Phi^{\mathrm{ob}}}) \mor F^{\sharp}([\mathbf{b}]_{\Phi^{\mathrm{ob}}})$ be a morphism in $\mathsf{Card}^{A}_{\mathrm{f}}$. By definition of the objects in the target category, $g_{a}$ is a mapping:
\[
g_{a} \colon \bb{\mathbf{a}}_{a} \mor \bb{\mathbf{b}}_{a}.
\]
We define $\varphi\colon \mathbf{a}\mor \mathbf{b}$ as follows. For every index $i \in \bb{\mathbf{a}}$:
\[
\varphi(i) = \overline{\mathrm{pos}}_{\mathbf{b}, \mathbf{a}(i)} \left( g_{\mathbf{a}(i)} \left( \overline{\mathrm{occ}}_{\mathbf{a}, \mathbf{a}(i)}(i) \right) \right).
\]
Now we verify that $\varphi$ is a morphism in $\mathsf{C}(\mathbf{A}^{\star})$. 
For every $i \in \bb{\mathbf{a}}$, we evaluate $\mathbf{b}(\varphi(i))$ using the definition of $\varphi(i)$:
\[
\mathbf{b}(\varphi(i)) = \mathbf{b} \left( \overline{\mathrm{pos}}_{\mathbf{b}, \mathbf{a}(i)} \left( g_{\mathbf{a}(i)} \left( \overline{\mathrm{occ}}_{\mathbf{a}, \mathbf{a}(i)}(i) \right) \right) \right).
\]
By the definition of the position mapping $\overline{\mathrm{pos}}_{\mathbf{b}, a} \colon \bb{\mathbf{b}}_{a} \mor \mathbf{b}^{-1}[\{a\}]$, for any element $j \in \bb{\mathbf{b}}_{a}$, the index $k = \overline{\mathrm{pos}}_{\mathbf{b}, a}(j)$ is such that $\mathbf{b}(k) = a$. In our case, $a = \mathbf{a}(i)$ and $j = g_{\mathbf{a}(i)} \left( \overline{\mathrm{occ}}_{\mathbf{a}, \mathbf{a}(i)}(i) \right)$. Note that $j$ is indeed an element of $\bb{\mathbf{b}}_{\mathbf{a}(i)}$ because $g_{\mathbf{a}(i)} \colon \bb{\mathbf{a}}_{\mathbf{a}(i)} \to \bb{\mathbf{b}}_{\mathbf{a}(i)}$. Therefore:
\[
\mathbf{b} \left( \overline{\mathrm{pos}}_{\mathbf{b}, \mathbf{a}(i)} \left( g_{\mathbf{a}(i)} \left( \overline{\mathrm{occ}}_{\mathbf{a}, \mathbf{a}(i)}(i) \right) \right) \right) = \mathbf{a}(i).
\]
Thus, $\mathbf{b}(\varphi(i)) = \mathbf{a}(i)$ for all $i \in \bb{\mathbf{a}}$, proving that $\varphi$ is a morphism in $\mathsf{C}(\mathbf{A}^{\star})$.

We now check that, for every $a \in A$, the component $F(\varphi)_a$ equals $g_a$. Let $j \in \bb{\mathbf{a}}_a$. By the definition of the functor $F$ on arrows:
\[
F(\varphi)_a(j) = \overline{\mathrm{occ}}_{\mathbf{b}, a} \left( \varphi|_{\mathbf{a}^{-1}[\{a\}]} \left( \overline{\mathrm{pos}}_{\mathbf{a}, a}(j) \right) \right).
\]
Let $k = \overline{\mathrm{pos}}_{\mathbf{a}, a}(j)$. Since $k \in \mathbf{a}^{-1}[\{a\}]$, we have $\mathbf{a}(k) = a$. We substitute the definition of $\varphi(k)$:
\[
F(\varphi)_a(j) = \overline{\mathrm{occ}}_{\mathbf{b}, a} \left( \overline{\mathrm{pos}}_{\mathbf{b}, a} \left( g_{a} \left( \overline{\mathrm{occ}}_{\mathbf{a}, a} \left( \overline{\mathrm{pos}}_{\mathbf{a}, a}(j) \right) \right) \right) \right).
\]
Since $\overline{\mathrm{occ}}_{\mathbf{a}, a} \circ \overline{\mathrm{pos}}_{\mathbf{a}, a} = \mathrm{id}_{\bb{\mathbf{a}}_a}$, the inner term $\overline{\mathrm{occ}}_{\mathbf{a}, a} ( \overline{\mathrm{pos}}_{\mathbf{a}, a}(j) )$ reduces to $j$. The expression becomes $\overline{\mathrm{occ}}_{\mathbf{b}, a} \left( \overline{\mathrm{pos}}_{\mathbf{b}, a} \left( g_{a}(j) \right) \right)$. Since $g_a(j) \in \bb{\mathbf{b}}_a$ and $\overline{\mathrm{occ}}_{\mathbf{b}, a} \circ \overline{\mathrm{pos}}_{\mathbf{b}, a} = \mathrm{id}_{\bb{\mathbf{b}}_a}$, the expression reduces to $g_a(j)$. Hence, $F(\varphi)_a(j) = g_a(j)$ for all $j \in \bb{\mathbf{a}}_a$ and all $a \in A$. This implies $F(\varphi) = g$, and therefore $F^{\sharp}([\varphi]) = g$, proving that $F^{\sharp}$ is full.

Therefore $F^{\sharp} \colon \mathsf{C}(\mathbf{A}^{\star})/{\equiv^{A}} \mor \mathsf{Card}^{A}_{\mathrm{f}}$ is an isomorphism.

This, combined with the fact that $F = J\circ\vs^{A}{(\,\bigcdot\,)}$ is an equivalence, implies that the quotient functor $P_{\equiv^{A}} \colon \mathsf{C}(\mathbf{A}^{\star}) \mor \mathsf{C}(\mathbf{A}^{\star})/{\equiv^{A}}$ is an equivalence as well. Note that $\mathrm{Ob}(\mathsf{C}(\mathbf{A}^{\star})/{\equiv^{A}})$ is $A^{\star}/{\equiv^{A}}$, the underlying set of $\mathbf{A}^{\star}/{\equiv^{A}}$. Algebraically, however, $\mathbf{A}^{\star}$ and its homomorphic image $\mathbf{A}^{\star}/{\equiv^{A}}$ are not isomorphic (unless $\mathrm{card}(A)\leq 1$).

We next establish that $\mathsf{C}(\mathbf{A}^{\star})/{\equiv^{A}}$ is equivalent to $\mathsf{C}(\mathbf{A}^{\star})/{\equiv^{A\natural}}$, where $\equiv^{A\natural}$ is the strong generalized congruence on $\mathsf{C}(\mathbf{A}^{\star})$ canonically associated to $\equiv^{A}$. Indeed, consider the choice function $F$ for $(\mathrm{Hom}_{\mathsf{C}(\mathbf{A}^{\star})}(\mathbf{a},\mathbf{b}))_{(\mathbf{a},\mathbf{b}) \in \equiv^{A}}$ that sends $(\mathbf{a},\mathbf{b})$ in $\equiv^{A}$  to $\sigma_{\mathbf{a},\mathbf{b}}$, the canonical permutation for $(\mathbf{a},\mathbf{b})$, in $\mathrm{Hom}_{\mathsf{C}(\mathbf{A}^{\star})}(\mathbf{a},\mathbf{b})$. Note that $F$ is \emph{constructive}, as it is defined by an explicit rule. By definition, $F$ satisfies the conditions stated in Definition~\ref{ChfComp}, i.e., it is a $\equiv^{A}$-compatible choice function for $(\mathrm{Hom}_{\mathsf{C}(\mathbf{A}^{\star})}(\mathbf{a},\mathbf{b}))_{(\mathbf{a},\mathbf{b}) \in \equiv^{A}}$. Consequently, Corollary~\ref{EqvQRQN} applies, yielding that $\mathsf{C}(\mathbf{A}^{\star})/{\equiv^{A}}$ and $\mathsf{C}(\mathbf{A}^{\star})/{\equiv^{A\natural}}$ are equivalent.
\end{example}

\begin{remark}
The Riguet congruence $\equiv^{A}$ on $\mathsf{C}(\mathbf{A}^{\star})$ determines a double category, denoted by  
$\equiv^{A}\!\!\mathsf{C}(\mathbf{A}^{\star})$, whose components are as follows. The set of objects is $\mathrm{Ob}(\mathsf{C}(\mathbf{A}^{\star}))$. The set of horizontal morphisms between the previous objects is $\mathrm{Mor}(\mathsf{C}(\mathbf{A}^{\star}))$. These horizontal morphisms, together with the objects, form the category $\mathsf{H}_{0}(\equiv^{A}\!\!\mathsf{C}(\mathbf{A}^{\star}))$. The set of vertical morphisms between the same objects is $\mathrm{Mor}(\mathsf{C}(\mathbf{A}^{\star})_{\mathsf{cp}})$, where $\mathsf{C}(\mathbf{A}^{\star})_{\mathsf{cp}}$ is the wide subcategory of $\mathsf{C}(\mathbf{A}^{\star})$ whose morphisms are the canonical permutations associated with the elements of $\equiv^{A}$. This wide subcategory is isomorphic to the groupoid $\mathsf{C}(\equiv^{A})$ associated to the equivalence $\equiv^{A}$. These vertical morphisms, together with the objects, form the category $\mathsf{V}_{0}(\equiv^{A}\!\!\mathsf{C}(\mathbf{A}^{\star}))$. A double cell $(\varphi,\varphi')$ with boundary formed of the vertical morphisms $\sigma_{\mathbf{a},\mathbf{a}'}$ (left), $\sigma_{\mathbf{b},\mathbf{b}'}$ (right) and the horizontal morphisms $\varphi$ (top), $\varphi'$ (bottom), denoted by $(\sigma_{\mathbf{a},\mathbf{a}'},\sigma_{\mathbf{b},\mathbf{b}'},\varphi,\varphi')$, exists whenever $(\varphi,\varphi')\in \equiv^{A}_{\scalebox{0.7}{$\left(\begin{smallmatrix}\mathbf{a}&\mathbf{b}\\ \mathbf{a}'&\mathbf{b}'\end{smallmatrix}\right)$}}$. Double cells $(\varphi,\varphi')$, $(\psi,\psi')$, where the boundary of $(\psi,\psi')$ is $(\sigma_{\mathbf{b},\mathbf{b}'},\sigma_{\mathbf{c},\mathbf{c}'}, \psi,\psi')$, have a \emph{horizontal} composition $(\psi\circ \varphi,\psi'\circ \varphi')$, consistent with the horizontal composition of morphisms. This composition gives the category $\mathsf{H}_{1}(\equiv^{A}\!\!\mathsf{C}(\mathbf{A}^{\star}))$ of vertical morphisms and cells. Double cells $(\varphi,\varphi')$, $(\varphi',\varphi'')$, where the boundary of $(\varphi',\varphi'')$ is $(\sigma_{\mathbf{a}',\mathbf{a}''},\sigma_{\mathbf{b}',\mathbf{b}''}, \varphi', \varphi'')$, have a \emph{vertical} composition $(\varphi,\varphi'')$, consistent with the vertical composition of morphisms. This composition gives the category $\mathsf{V}_{1}(\equiv^{A}\!\!\mathsf{C}(\mathbf{A}^{\star}))$ of horizontal morphisms and cells. These two compositions satisfy the interchange laws. 

Moreover, since, for every injective mapping $f$ from $A$ to $B$, we have that  
\begin{enumerate}
\item for every pair of words $\mathbf{a}, \mathbf{a}' \in A^{\star}$ such that $\mathbf{a} \equiv^{A} \mathbf{a}'$, the canonical permutations $\sigma_{\mathbf{a}, \mathbf{a}'}$ and $\sigma_{f^{@}(\mathbf{a}), f^{@}(\mathbf{a}')}$ are equal, where $f^{@}$ is $(\eta_{B}\circ f)^{\sharp}$, the canonical extension of $\eta_{B}\circ f$, and 
\item for every words $\mathbf{a}, \mathbf{a}', \mathbf{b}, \mathbf{b}' \in A^{\star}$ such that $\mathbf{a} \equiv^{A} \mathbf{a}'$ and $\mathbf{b} \equiv^{A} \mathbf{b}'$, every mapping $\varphi$ from $\bb{\mathbf{a}}$ to $\bb{\mathbf{b}}$ such that $\mathbf{a}=\mathbf{b} \circ \varphi$ and every mapping $\varphi'$ from $\bb{\mathbf{a}'}$ to $\bb{\mathbf{b}'}$ such that $\mathbf{a}'=\mathbf{b}' \circ \varphi'$, if $\varphi' \circ \sigma_{\mathbf{a}, \mathbf{a}'} = \sigma_{\mathbf{b}, \mathbf{b}'} \circ \varphi$, then $\varphi' \circ \sigma_{f^{@}(\mathbf{a}), f^{@}(\mathbf{a}')} = \sigma_{f^{@}(\mathbf{b}), f^{@}(\mathbf{b}')} \circ \varphi$,
\end{enumerate}
we can assert that any injective mapping from $A$ to $B$ induces a double functor from $\equiv^{A}\!\!\mathsf{C}(\mathbf{A}^{\star})$ to $\equiv^{B}\!\!\mathsf{C}(\mathbf{B}^{\star})$. Therefore, there exists a functor from $\mathsf{Set}_{\mathrm{inj}}$, the category 
of sets and injective mappings, to the category $\mathsf{DblCat}$, of double categories and double functors.
\end{remark} 

\subsection{A Riguet Congruence on the Fundamental Groupoid of a Path-Connected Topological Space}\label{S:groupoid}

\begin{example}
Let $\mathbf{X}$ be a path-connected topological space, and let $\mathsf{\Pi}_{1}(\mathbf{X})$ denote its fundamental groupoid. The objects of $\mathsf{\Pi}_{1}(\mathbf{X})$ are the points of $X$, the underlying set of $\mathbf{X}$, and the morphisms are homotopy classes of paths. Since $\mathbf{X}$ is path-connected, $\mathrm{Hom}_{\mathsf{\Pi}_{1}(\mathbf{X})}(x, y) \neq \emptyset$ for all $x, y \in X$; thus, $\mathsf{\Pi}_{1}(\mathbf{X})$ is a strongly connected category. We define a Riguet congruence $\Phi$ on $\mathsf{\Pi}_{1}(\mathbf{X})$ as follows:
\begin{enumerate}
\item \textsf{On Objects:} We set $\Phi^{\mathrm{ob}} = \nabla_{X}$, so that all objects are identified. Fix a basepoint $x_{0} \in X$ and, for each $x \in X$, choose a homotopy class $\alpha_x \colon x_{0} \mor x$. For every pair $(x, y) \in X \times X$, we define the transport isomorphism as $\alpha_{x,y} = \alpha_{y} \circ \alpha_{x}^{-1} \colon x \mor y$.
    
\item \textsf{On Arrows:} For $f \colon x \mor y$ and $f' \colon x' \mor y'$, we define  $(f,f') \in \Phi^{\mathrm{fl}}_{\scalebox{0.7}{$\left(\begin{smallmatrix}x&y\\x'& y'\end{smallmatrix}\right)$}}$ if, and only if, $f' = \alpha_{y,y'} \circ f \circ \alpha_{x,x'}^{-1}$. Substituting our choice of paths, this is equivalent to $f' = (\alpha_{y'} \circ \alpha_{y}^{-1}) \circ f \circ (\alpha_{x} \circ \alpha_{x'}^{-1})$. 
\end{enumerate}
We restrict ourselves to verifying condition (f), leaving the verification of conditions (a) through (e) to the reader as they present no particular difficulty. The quadrangular completion condition $(\mathrm{f})$ is satisfied as follows. Let $f \colon x \mor y$ and let $(x,x'),(y,y') \in \Phi^{\mathrm{ob}}$. Define
\[
f' = \alpha_{y,y'} \circ f \circ \alpha_{x,x'}^{-1} \in \mathrm{Hom}_{\mathsf{\Pi}_{1}(\mathbf{X})}(x',y').
\]
Then, by definition, $(f,f') \in \Phi^{\mathrm{fl}}_{\scalebox{0.7}{$\left(\begin{smallmatrix}x&y\\x'& y'\end{smallmatrix}\right)$}}$, so condition $(\mathrm{f})$ holds.

We define a functor $F \colon \mathsf{\Pi}_{1}(\mathbf{X}) \mor \mathsf{B}\pi_1(\mathbf{X},x_{0})$, where 
$\mathsf{B}\pi_1(\mathbf{X},x_{0})$ is the one-object category associated with the group $\pi_{1}(\mathbf{X},x_{0})$, as follows: For an object $x$ of $\mathsf{\Pi}_{1}(\mathbf{X})$, $F(x)$ is the unique object $*$ of $\mathsf{B}\pi_1(\mathbf{X},x_{0})$. For a morphism $f \colon x \mor y$ of $ \mathsf{\Pi}_{1}(\mathbf{X})$, $F(f) = \alpha_y^{-1} \circ f \circ \alpha_x \in \pi_1(\mathbf{X},x_0)$. 
Note that $F$ realizes the classical equivalence of categories $\mathsf{\Pi}_1(\mathbf{X}) \simeq \mathsf{B}\pi_1(\mathbf{X}, x_0)$. We next prove that $F$ satisfies the hypotheses of Proposition~\ref{UPQ}:  
\begin{enumerate}
\item If $(x,x') \in \Phi^{\mathrm{ob}}$, then $F(x) = F(x') = *$.
\item If $(f,f') \in \Phi^{\mathrm{fl}}_{\scalebox{0.7}{$\left(\begin{smallmatrix}x&y\\x'& y'\end{smallmatrix}\right)$}}$, then, by definition of $\Phi$, we have that
$
f' = \alpha_{y,y'} \circ f \circ \alpha_{x,x'}^{-1}.
$
Applying $F$, we obtain
\[
\begin{aligned}
F(f') 
&= \alpha_{y'}^{-1} \circ f' \circ \alpha_{x'} \\
&= \alpha_{y'}^{-1} \circ (\alpha_{y,y'} \circ f \circ \alpha_{x,x'}^{-1}) \circ \alpha_{x'} \\
&= \alpha_{y'}^{-1} \circ (\alpha_{y'} \circ \alpha_y^{-1}) \circ f \circ (\alpha_x \circ \alpha_{x'}^{-1}) \circ \alpha_{x'} \\
&= \alpha_y^{-1} \circ f \circ \alpha_x \\
&= F(f).
\end{aligned}
\]
Hence $F(f) = F(f')$, as required.
\end{enumerate}
Thus, by Proposition~\ref{UPQ}, there exists a unique functor 
$F^{\sharp}$ from $\mathsf{\Pi}_{1}(\mathbf{X})/\Phi$ to $\mathsf{B}\pi_1(\mathbf{X},x_{0})$ 
such that $F^{\sharp} \circ P_\Phi = F$. 

To show that $F^{\sharp}$ is an isomorphism, we verify that it is bijective on objects, faithful, and full.

\textsf{Bijectivity on objects.}
The bijectivity on objects of $F^{\sharp}$ is trivial, since $\mathsf{\Pi}_1(\mathbf{X})/\Phi$ has a single object $[x]_{\Phi^{\mathrm{ob}}}$ and $\mathsf{B}\pi_1(\mathbf{X},x_{0})$ has a single object $*$.

\textsf{Faithfulness.}
Suppose that $[f]_{\Phi^{\mathrm{fl}}},[g]_{\Phi^{\mathrm{fl}}} \colon [x]_{\Phi^{\mathrm{ob}}} \mor  [y]_{\Phi^{\mathrm{ob}}}$ satisfy $F^\sharp([f]_{\Phi^{\mathrm{fl}}}) = F^\sharp([g]_{\Phi^{\mathrm{fl}}})$. Let $f_{0} \colon a_{0} \mor b_{0}$ and $f_{1} \colon a_{1} \mor b_{1}$ be representatives. By the definition of $F^\sharp$, we have
\[
\alpha_{b_{0}}^{-1} \circ f_{0} \circ \alpha_{a_{0}} = \alpha_{b_{1}}^{-1} \circ f_{1} \circ \alpha_{a_{1}} \in \pi_1(\mathbf{X},x_0),
\]
which implies $f_{1} = \alpha_{b_{0} ,b_{1} } \circ f_{0}  \circ \alpha_{a_{0} ,a_{1} }^{-1}$. 
But this is exactly the condition for $(f_{0},f_{1}) \in \Phi^{\mathrm{fl}}_{\scalebox{0.7}{$\left(\begin{smallmatrix}a_{0}&b_{0}\\a_{1}& b_{1}\end{smallmatrix}\right)$}}$, hence $[f_{0}]_{\Phi^{\mathrm{fl}}} = [f_{1}]_{\Phi^{\mathrm{fl}}}$, proving faithfulness.

\textsf{Fullness.}
Let $\gamma \in \pi_1(\mathbf{X},x_0)$ be any morphism in $\mathsf{B}\pi_1(\mathbf{X},x_0)$. 
We need to show that there exists a morphism 
$[f]_{\Phi^{\mathrm{fl}}} \colon [x]_{\Phi^{\mathrm{ob}}} \mor [y]_{\Phi^{\mathrm{ob}}}$ in 
$\mathsf{\Pi}_1(\mathbf{X})/\Phi$ such that 
$F^\sharp([f]_{\Phi^{\mathrm{fl}}}) = \gamma$. Consider the morphism $f \colon x \mor y$ in $\mathsf{\Pi}_1(\mathbf{X})$ defined by
\[
f = \alpha_y \circ \gamma \circ \alpha_x^{-1},
\]
where $\alpha_x \colon x_{0} \mor x$ and $\alpha_y \colon x_{0} \mor y$ are the chosen paths in $\mathbf{X}$ for the objects $x$ and $y$. 
By definition of $F^\sharp$, we have that
\[
F^\sharp([f]_{\Phi^{\mathrm{fl}}}) = \alpha_y^{-1} \circ f \circ \alpha_x 
= \alpha_y^{-1} \circ (\alpha_y \circ \gamma \circ \alpha_x^{-1}) \circ \alpha_x 
= \gamma.
\]
Hence every morphism in $\mathsf{B}\pi_1(\mathbf{X},x_0)$ has a preimage under $F^\sharp$, 
which proves that $F^\sharp$ is full.

Therefore, we conclude that the quotient category $\mathsf{\Pi}_{1}(\mathbf{X})/\Phi$ is canonically isomorphic to $\mathsf{B}\pi_1(\mathbf{X},x_{0})$, providing a concrete realization of the passage from the fundamental groupoid to the fundamental group via a Riguet congruence.

This approach refines the well-known fact that, for a path-connected topological space $\mathbf{X}$, the categories $\mathsf{\Pi}_{1}(\mathbf{X})$ and $\mathsf{B}\pi_{1}(\mathbf{X},x_{0})$ are equivalent, upgrading the classical equivalence to an actual isomorphism via an explicit Riguet congruence, specified at the level of both objects and morphisms.
\end{example}

\subsection{A Riguet Congruence for a Grothendieck Category with a Localized Subcategory}\label{S:grothendieck}

\begin{example}
Let $\mathsf{A}$ be a Grothendieck Category and $\mathsf{S}$ a localizing subcategory of $\mathsf{A}$. 
In this example we show that there exists a Riguet congruence $\Phi$ on $\mathsf{A}$ such that the categories $\mathsf{A}/\mathsf{S}$ and $\mathsf{A}/\Phi$ are equivalent. 
To establish this, we first recall the general construction of the quotient category for an arbitrary abelian category $\mathsf{A}$ and a Serre subcategory $\mathsf{S}$ of $\mathsf{A}$. We briefly outline the definition of $\mathsf{A}/\mathsf{S}$ and the universal property of the pair $(\mathsf{A}/\mathsf{S}, Q)$, where $Q \colon \mathsf{A} \mor \mathsf{A}/\mathsf{S}$ denotes the canonical localization functor. This background serves to highlight how the additional structure of Grothendieck categories---specifically the presence of a section functor---allows the Riguet congruence to internally realize this localization.  Next, by considering the standard skeleton $\mathsf{T}_{\cong}$ of $\mathsf{A}/\mathsf{S}$, we define a functor $F\colon \mathsf{A} \mor \mathsf{T}_{\cong}$ via the localization functor $Q$, and verify that $F$ satisfies conditions~(1) and~(2) of Proposition~\ref{UPQ}. Hence, Proposition~\ref{UPQ} yields a unique functor $F^{\sharp}\colon \mathsf{A}/\Phi \mor \mathsf{T}_{\cong}$ such that $F = F^{\sharp} \circ P_{\Phi}$. Moreover, $F^{\sharp}$ is an isomorphism, and hence $\mathsf{A}/\Phi$ and $\mathsf{A}/\mathsf{S}$ are equivalent.

The quotient category $\mathsf{A}/\mathsf{S}$ is constructed as the localization of $\mathsf{A}$ with respect to the set of $\mathsf{S}$-isomorphisms. A morphism $s \colon M \mor N$ in $\mathsf{A}$ is called an $\mathsf{S}$-isomorphism if $\mathrm{Ker}(s)$ and  $\mathrm{Coker}(s)\in \mathrm{Ob}(\mathsf{S})$. We denote the set of all such morphisms as $\Sigma_ {\mathsf{S}}$.
The quotient category $\mathsf{A}/\mathsf{S}$, also denoted by $\mathsf{A}[\Sigma^{-1}_{\mathsf{S}}]$, is defined by:
\begin{enumerate}
\item \textsf{Objects:} $\mathrm{Ob}(\mathsf{A}/\mathsf{S}) = \mathrm{Ob}(\mathsf{A})$.
\item \textsf{Morphisms:} A morphism $f \in \mathrm{Hom}_{\mathsf{A}/\mathsf{S}}(M, N)$ is an equivalence class of spans $(s,g)$ in $\mathsf{A}$ where $s\colon X\mor M \in \Sigma_{\mathsf{S}}$ and $g\colon X\mor N$ is any morphism in $\mathsf{A}$.
\end{enumerate}
Two spans $(s, g)$ and $(s', g')$, where $s'\colon X'\mor M \in \Sigma_{\mathsf{S}}$ and $g'\colon X'\mor N$ is any morphism in $\mathsf{A}$, from $M$ to $N$ represent the same morphism in $\mathsf{A}/\mathsf{S}$ if there exists an object $X'' \in \mathrm{Ob}(\mathsf{A})$ and morphisms $u\colon X'' \mor X$ and $v\colon X'' \mor X'$ such that $s \circ u = s' \circ v$, $g \circ u = g' \circ v$ and $s \circ u \in \Sigma_ {\mathsf{S}}$ (which implies $s' \circ v \in \Sigma_ {\mathsf{S}}$). We denote by $[s, g]$ the equivalence class of $(s,g)$.

In an abelian category $\mathsf{A}$, if $\mathsf{S}$ is a Serre subcategory, the set of $\mathsf{S}$-isomorphisms $\Sigma_{\mathsf{S}}$ automatically satisfies the Ore conditions. 
For every object $M \in \mathrm{Ob}(\mathsf{A}/\mathsf{S})$, the identity morphism at $M$ is $\mathrm{id}_{M} = [\mathrm{id}_{M},\mathrm{id}_{M}]$. Let $[s,g]$ be a morphism from $M$ to $N$ and $[t,h]$ a morphism from $N$ to $P$, where $t\colon Y\mor N \in \Sigma_{\mathsf{S}}$ and $h\colon Y\mor P$ is any morphism in $\mathsf{A}$. Their composition $[t,h] \circ [s,g]$ is $[s\circ s',h\circ g']$, where $(W, (s',g'))$ is the  pullback of $(g,t)$ in $\mathsf{A}$.

The localization functor $Q \colon \mathsf{A} \mor \mathsf{A}/\mathsf{S}$ is defined as follows. The object mapping of $Q$ is the identity at $\mathrm{Ob}(\mathsf{A})$ and the morphism mapping assigns to $f\colon M\mor N$ precisely $[\mathrm{id}_{M},f]$. The functor $Q$ is exact, sends every morphism $s \in \Sigma_ {\mathsf{S}}$ to an isomorphism and, for every object $M\in \mathrm{Ob}(\mathsf{A})$, $M\in \mathrm{Ob}(\mathsf{S})$ if, and only if, $Q(M)\cong 0$. Moreover, it satisfies the following universal property: for every abelian category $\mathsf{B}$ and exact functor $F \colon \mathsf{A} \mor \mathsf{B}$ such that $F(M)\cong 0$ for all $M \in \mathrm{Ob}(\mathsf{S})$, there exists a unique exact functor $F^{\natural}$ from $ \mathsf{A}/\mathsf{S}$ to $\mathsf{B}$ such that $F = F^{\natural}\circ Q$.

While the general construction of $\mathsf{A}/\mathsf{S}$ via fractions applies to any abelian category, the Riguet congruence $\Phi$ requires that every morphism in the quotient be representable by a morphism in $\mathsf{A}$. This condition is satisfied if, and only if, the localization functor $Q$ is full. In the specific context of Grothendieck categories, where $\mathsf{S}$ is a localizing subcategory, $Q$ admits a right adjoint $R \colon \mathsf{A}/\mathsf{S} \mor \mathsf{A}$ known as the section functor (see Gabriel~\cite{pG62} and Popescu~\cite{NP73}). The counit $\varepsilon \colon QR \mor \mathrm{Id}_{\mathsf{A}/\mathsf{S}}$ is a natural isomorphism, which implies that $R$ is fully faithful and, consequently, that $Q$ is full. This fullness is essential for verifying condition (f) in the definition of Riguet's congruence, as it ensures that the abstract fractions of the localization can be internally realized as morphisms in $\mathsf{A}$. Specifically, the unit of the adjunction $\eta \colon \mathrm{Id}_{\mathsf{A}} \mor RQ$ is such that for every $M \in \mathrm{Ob}(\mathsf{A})$, both $\mathrm{Ker}(\eta_M)$ and $\mathrm{Coker}(\eta_M)$ belong to $\mathsf{S}$. This structure allows $\Phi$ to provide a complete internal description of the quotient, identifying morphisms between $\mathsf{S}$-closed objects within the original category.

We define a Riguet congruence $\Phi$ on $\mathsf{A}$ as follows:
\begin{enumerate}
\item \textsf{On Objects:} $(M,N) \in \Phi^{\mathrm{ob}}$ if, and only if, $Q(M) \cong Q(N)$ in $\mathsf{A}/\mathsf{S}$. We then fix, for every pair $(M,N) \in \Phi^{\mathrm{ob}}$, an isomorphism $\alpha_{M,N} \colon Q(M) \mor Q(N)$ subject to the following compatibility conditions:
\begin{enumerate}
\item $\alpha_{M, M} = \mathrm{id}_{Q(M)}$ for every $M \in \mathrm{Ob}(\mathsf{A})$;
\item $\alpha_{N, M} = \alpha_{M,N}^{-1}$ for every $(M,N) \in \Phi^{\mathrm{ob}}$;
\item $\alpha_{M, P} = \alpha_{N, P} \circ \alpha_{M,N}$ for every $M,N,P \in \mathrm{Ob}(\mathsf{A})$ such that the pairs $(M,N)$, $(N,P)$ and $(M,P)$ are in $\Phi^{\mathrm{ob}}$.
\end{enumerate}
    
\item \textsf{On Arrows:} For every $f \in \mathrm{Hom}_{\mathsf{A}}(M,N)$ and $f' \in \mathrm{Hom}_{\mathsf{A}}(M', N')$ with $(M, M')$ and $(N, N') \in \Phi^{\mathrm{ob}}$, we declare $(f, f') \in \Phi^{\mathrm{fl}}_{\scalebox{0.7}{$\left(\begin{smallmatrix}M&N\\M'& N'\end{smallmatrix}\right)$}}$ if, and only if, in 
$\mathsf{A}/\mathsf{S}$, we have that $Q(f') = \alpha_{N, N'} \circ Q(f) \circ \alpha_{M, M'}^{-1}$.
\end{enumerate}
Let us note that $\Phi^{\mathrm{ob}}$ admits an equivalent characterization within the original category $\mathsf{A}$: $(M, N) \in \Phi^{\mathrm{ob}}$ if, and only if, $RQ(M) \cong RQ(N)$ in $\mathsf{A}$. To see that these definitions are indeed equivalent, suppose first that $Q(M) \cong Q(N)$ in $\mathsf{A}/\mathsf{S}$. By applying the section functor $R$, we immediately obtain $RQ(M) \cong RQ(N)$ in $\mathsf{A}$. Conversely, if $RQ(M) \cong RQ(N)$ in $\mathsf{A}$, applying the localization functor yields $QRQ(M) \cong QRQ(N)$. By virtue of the natural isomorphism $\varepsilon \colon QR \cong \mathrm{Id}_{\mathsf{A}/\mathsf{S}}$, it follows that $Q(M) \cong Q(N)$ in $\mathsf{A}/\mathsf{S}$

By definition $\Phi^{\mathrm{ob}}$ is an equivalence relation.

The conditions $(\mathrm{a})$ through $(\mathrm{e})$ are satisfied by the functoriality of $Q$ and the  properties imposed on the family $(\alpha_{M,N})_{(M.N)\in \Phi^{\mathrm{ob}}}$. 

(f) By the definition of $\Phi^{\mathrm{ob}}$, the pairs $(M, M')$ and $(N, N')$ correspond to fixed isomorphisms $\alpha_{M, M'} \colon Q(M) \mor Q(M')$ and $\alpha_{N, N'} \colon Q(N) \mor Q(N')$ in $\mathsf{A}/\mathsf{S}$. 

Consider the composed morphism in the quotient category:
\[ \gamma = \alpha_{N, N'} \circ Q(f) \circ \alpha_{M, M'}^{-1} \in \mathrm{Hom}_{\mathsf{A}/\mathsf{S}}(Q(M'), Q(N')). \]
Since $Q$ is a full functor (as established by the existence of the fully faithful right adjoint $R$ in a Grothendieck setting), the map:
\[ Q_{M', N'} \colon \mathrm{Hom}_{\mathsf{A}}(M', N') \mor \mathrm{Hom}_{\mathsf{A}/\mathsf{S}}(Q(M'), Q(N')) \]
is surjective. Therefore, there exists at least one $f' \in \mathrm{Hom}_{\mathsf{A}}(M', N')$ such that $Q(f') = \gamma$. By the definition of the congruence on arrows, this equality $Q(f') = \alpha_{N, N'} \circ Q(f) \circ \alpha_{M, M'}^{-1}$ is precisely the condition for $(f, f') \in \Phi^{\mathrm{fl}}$.

To prove that $\mathsf{A}/\Phi$ and $\mathsf{A}/\mathsf{S}$ are equivalent, we proceed as follows. 
If we let $M \cong N$ mean that there is an $\mathsf{A}/\mathsf{S}$-isomorphism from $M$ to $N$, then $\cong$ is an equivalence relation on $\mathrm{Ob}(\mathsf{A}/\mathsf{S})$. Hence, by the Axiom of Choice it has a transversal $T_{\cong}$; that is, for every equivalence class $[M]_{\cong}$, there exists a unique representative $\overline{M} \in T_{\cong}$ such that $[M]_{\cong} = [\overline{M}]_{\cong}$.  Let $\mathsf{T}_{\cong}$ be the full subcategory of $\mathsf{A}/\mathsf{S}$ that is generated by $T_{\cong}$. Since $\mathsf{T}_{\cong}$ is the standard skeleton of $\mathsf{A}/\mathsf{S}$, the categories $\mathsf{A}/\mathsf{S}$ and $\mathsf{T}_{\cong}$ are equivalent. 
We define a functor $F \colon \mathsf{A} \mor \mathsf{T}_{\cong}$ as follows: For $M \in \mathrm{Ob}(\mathsf{A})$, $F(M) = \overline{M}$, where $\overline{M} \in T_{\cong}$ is the unique representative such that $Q(M) \cong Q(\overline{M})$ in $\mathsf{A}/\mathsf{S}$.
For $f \in \mathrm{Hom}_{\mathsf{A}}(M, N)$, $F(f) \in \mathrm{Hom}_{\mathsf{T}_{\cong}}(\overline{M}, \overline{N})$ is defined as 
$
F(f) = \alpha_{N, \overline{N}} \circ Q(f) \circ \alpha_{M, \overline{M}}^{-1}.
$

We next verify for $F$ the conditions in Proposition~\ref{UPQ}. 
(1) If $(M, M') \in \Phi^{\mathrm{ob}}$, then $Q(M) \cong Q(M')$. Since they belong to the same isomorphism class, their skeletal representative is identical: $F(M) = \overline{M} = \overline{M'} = F(M')$. 
(2) If $(f, f') \in \Phi^{\mathrm{fl}}_{\scalebox{0.7}{$\left(\begin{smallmatrix}M&N\\M'& N'\end{smallmatrix}\right)$}}$, then $Q(f') = \alpha_{N, N'} \circ Q(f) \circ \alpha_{M, M'}^{-1}$. Thus:
    \begin{align*}
        F(f') &= \alpha_{N', \overline{N}} \circ Q(f') \circ \alpha_{M', \overline{M}}^{-1} \\
        &= \alpha_{N', \overline{N}} \circ \left(\alpha_{N, N'} \circ Q(f) \circ \alpha_{M, M'}^{-1}\right) \circ \alpha_{M', \overline{M}}^{-1} \\
        &= \left(\alpha_{N', \overline{N}} \circ \alpha_{N, N'}\right) \circ Q(f) \circ \left(\alpha_{M, M'}^{-1} \circ \alpha_{M', \overline{M}}^{-1}\right) \\
        &= \alpha_{N, \overline{N}} \circ Q(f) \circ \alpha_{M, \overline{M}}^{-1} 
        \\&= F(f).
    \end{align*}
Hence, by Proposition~\ref{UPQ}, there exists a unique functor $F^{\sharp} \colon \mathsf{A}/\Phi \mor \mathsf{T}_{\cong}$ such that $F^{\sharp} \circ P_{\Phi} = F$.

We next prove that $F^{\sharp}$ is an isomorphism of categories. The surjectivity of the object mapping of $F^{\sharp}$ holds by the definition of the transversal. For the injectivity  of the object mapping of $F^{\sharp}$, if $F^{\sharp}([M]_{\Phi^{\mathrm{ob}}}) = F^{\sharp}([N]_{\Phi^{\mathrm{ob}}})$ then $\overline{M} = \overline{N}$, thus  $Q(M) \cong Q(N)$, hence $(M,N) \in \Phi^{\mathrm{ob}}$, therefore $[M]_{\Phi^{\mathrm{ob}}} = [N]_{\Phi^{\mathrm{ob}}}$. 

To prove that $F^{\sharp} \colon \mathrm{Hom}_{\mathsf{A}/\Phi}([M]_{\Phi^{\mathrm{ob}}}, [N]_{\Phi^{\mathrm{ob}}}) \mor \mathrm{Hom}_{ \mathsf{T}_{\cong}}(\overline{M}, \overline{N})$ is a bijection, we show that it is both full and faithful. 

\textsf{Faithfulness.}
Suppose $F^{\sharp}([f]_{\Phi^{\mathrm{fl}}}) = F^{\sharp}([g]_{\Phi^{\mathrm{fl}}})$ for two morphisms $[f]_{\Phi^{\mathrm{fl}}}$, $[g]_{\Phi^{\mathrm{fl}}}$ from $[M]_{\Phi^{\mathrm{ob}}}$ to  $[N]_{\Phi^{\mathrm{ob}}}$. This implies that the representatives $f$ and $g$ in $\mathsf{A}$ may have different domains and codomains, say $f \colon M \mor N$ and $g \colon M' \mor N'$, provided that $(M, M') \in \Phi^{\mathrm{ob}}$ and $(N, N') \in \Phi^{\mathrm{ob}}$. By the definition of $\Phi^{\mathrm{ob}}$, we have that $Q(M) \cong Q(M')$ and $Q(N) \cong Q(N')$, which implies they share the same unique representatives in the transversal: $\overline{M} = \overline{M'}$ and $\overline{N} = \overline{N'}$. From the definition of $F$ on morphisms and the equality $F(f) = F(g)$, we have: \[
\alpha_{N, \overline{N}} \circ Q(f) \circ \alpha_{M, \overline{M}}^{-1} = \alpha_{N', \overline{N}} \circ Q(g) \circ \alpha_{M', \overline{M}}^{-1}.
\]
To prove that $[f]_{\Phi^{\mathrm{fl}}} = [g]_{\Phi^{\mathrm{fl}}}$, we must verify that $(f, g) \in \Phi^{\mathrm{fl}}_{\scalebox{0.7}{$\left(\begin{smallmatrix}M&N\\M'& N'\end{smallmatrix}\right)$}}$, which requires:
\[
Q(g) = \alpha_{N, N'} \circ Q(f) \circ \alpha_{M, M'}^{-1}.
\]
We isolate $Q(g)$ from the first equation by composing with the respective isomorphisms:
\[
Q(g) = \alpha_{N', \overline{N}}^{-1} \circ (\alpha_{N, \overline{N}} \circ Q(f) \circ \alpha_{M, \overline{M}}^{-1}) \circ \alpha_{M', \overline{M}}.
\]
By the compatibility conditions (b) and (c) imposed on the family $(\alpha_{M,N})_{(M,N) \in \Phi^{\mathrm{ob}}}$:
\allowdisplaybreaks
\begin{align*}
     \alpha_{N', \overline{N}}^{-1} \circ \alpha_{N, \overline{N}} &= \alpha_{\overline{N}, N'} \circ \alpha_{N, \overline{N}} = \alpha_{N, N'}. 
     \\\alpha_{M, \overline{M}}^{-1} \circ \alpha_{M', \overline{M}} &= \alpha_{\overline{M}, M} \circ \alpha_{M', \overline{M}} = \alpha_{M', M} = \alpha_{M, M'}^{-1}.
\end{align*}
Substituting these identities back into the expression for $Q(g)$, we obtain:
\[
Q(g) = \alpha_{N, N'} \circ Q(f) \circ \alpha_{M, M'}^{-1}.
\]
This is precisely the definition of the relation $(f, g) \in \Phi^{\mathrm{fl}}_{\scalebox{0.7}{$\left(\begin{smallmatrix}M&N\\M'& N'\end{smallmatrix}\right)$}}$, which means $[f]_{\Phi^{\mathrm{fl}}} = [g]_{\Phi^{\mathrm{fl}}}$. Thus, $F^{\sharp}$ is faithful.

\textsf{Fullness.}
Let $\gamma \in \mathrm{Hom}_{\mathsf{T}_{\cong}}(\overline{M}, \overline{N})$. Define $\psi = \alpha_{N, \overline{N}}^{-1} \circ \gamma \circ \alpha_{M, \overline{M}}$, which is a morphism from $Q(M)$  to $Q(N)$ in $\mathsf{A}/\mathsf{S}$. 
Since $\mathsf{A}$ is Grothendieck, the localization functor $Q$ is full. Thus, there exists an $f \in \mathrm{Hom}_{\mathsf{A}}(M, N)$ such that $Q(f) = \psi$. 
Then $F^{\sharp}([f]_{\Phi^{\mathrm{fl}}}) = \alpha_{N, \overline{N}} \circ (\alpha_{N, \overline{N}}^{-1} \circ \gamma \circ \alpha_{M, \overline{M}}) \circ \alpha_{M, \overline{M}}^{-1} = \gamma$.

Therefore, $F^{\sharp}$ is an isomorphism of categories. Since $\mathsf{T}_{\cong} \simeq \mathsf{A}/\mathsf{S}$, we have $\mathsf{A}/\Phi \simeq \mathsf{A}/\mathsf{S}$.

It is important to verify that the category $\mathsf{A}/\Phi$ is independent of the particular choice of the family of isomorphisms $\alpha = (\alpha_{M,N})_{(M,N) \in \Phi^{\mathrm{ob}}}$. Let $\alpha = (\alpha_{M,N})_{(M,N)\in \Phi^{\mathrm{ob}}}$ and $\beta = (\beta_{M,N})_{(M,N)\in \Phi^{\mathrm{ob}}}$ be two families of isomorphisms satisfying the compatibility conditions. These define two Riguet congruences, $\Phi_{\alpha}$ and $\Phi_{\beta}$, respectively. Recall that for every equivalence class $[M]_{\cong}$, there exists a unique representative $\overline{M} \in T_{\cong}$ such that $[M]_{\cong} = [\overline{M}]_{\cong}$. This unique correspondence allows us to define, for each object $M \in \mathrm{Ob}(\mathsf{A})$, a transition isomorphism $\theta_{M} \in \mathrm{Aut}_{\mathsf{A}/\mathsf{S}}(Q(M))$ given by: $\theta_{M} = \beta_{M, \overline{M}}^{-1} \circ \alpha_{M, \overline{M}}$. 
To prove that the categories are isomorphic, consider the assignment $F_{\alpha,\beta} \colon \mathsf{A}/\Phi_{\alpha} \mor \mathsf{A}/\Phi_{\beta}$ defined as: for every object $[M]_{\Phi^{\mathrm{ob}}_{\alpha}}$ of $\mathsf{A}/\Phi_{\alpha}$, $F_{\alpha,\beta}([M]_{\Phi^{\mathrm{ob}}_{\alpha}}) = [M]_{\Phi^{\mathrm{ob}}_{\beta}}$, and, for every morphism $[f]_{\Phi^{\mathrm{fl}}_{\alpha}}$ of $\mathsf{A}/\Phi_{\alpha}$, $F_{\alpha,\beta}([f]_{\Phi^{\mathrm{fl}}_{\alpha}}) = [f]_{\Phi^{\mathrm{fl}}_{\beta}}$. 
To show that $F_{\alpha,\beta}$ is well-defined, we must prove that for every $(M, M')$, $(N, N') \in \Phi^{\mathrm{ob}}_{\alpha}$ and every $f \in \mathrm{Hom}_{\mathsf{A}}(M, N)$ and $f' \in \mathrm{Hom}_{\mathsf{A}}(M', N')$, we have that $(f,f') \in \Phi^{\mathrm{fl}}_{\alpha\scalebox{0.7}{$\left(\begin{smallmatrix}M&N\\M'& N'\end{smallmatrix}\right)$}}$ if, and only if, $(f,f') \in \Phi^{\mathrm{fl}}_{\beta\scalebox{0.7}{$\left(\begin{smallmatrix}M&N\\M'& N'\end{smallmatrix}\right)$}}$. 
To incorporate the transition between families, we first justify the naturality $\theta_N Q(f) = Q(f) \theta_M$. We observe that for every $f \colon M \mor N$, both families $\alpha$ and $\beta$ must yield the same representative morphism in the skeleton $\mathsf{T}_{\equiv}$ because they target the same unique representatives $\overline{M}, \overline{N} \in T_{\cong}$:
\[
\alpha_{N, \overline{N}} Q(f) \alpha_{M, \overline{M}}^{-1} = \beta_{N, \overline{N}} Q(f) \beta_{M, \overline{M}}^{-1}. 
\]
Rearranging terms to isolate $Q(f)$ on the left-hand side, we have that:
\[ Q(f) = \left(\alpha_{N, \overline{N}}^{-1} \beta_{N, \overline{N}}\right) Q(f) \left(\beta_{M, \overline{M}}^{-1} \alpha_{M, \overline{M}}\right) = \theta_N^{-1} Q(f) \theta_M. \]
Hence $\theta_N Q(f) \theta_M^{-1} = Q(f)$. Now, by the compatibility conditions, the isomorphisms between any $M, M'$ and $N, N'$ can be decomposed as $\alpha_{M,M'} = \alpha_{M', \overline{M}}^{-1} \circ \alpha_{M, \overline{M}}$ and $\beta_{M,M'} = \beta_{M', \overline{M}}^{-1} \circ \beta_{M, \overline{M}}$. The condition $Q(f') = \alpha_{N, N'} \circ Q(f) \circ \alpha_{M, M'}^{-1}$ is equivalent to:
\[ 
Q(f') = \left(\alpha_{N', \overline{N}}^{-1} \alpha_{N, \overline{N}}\right) \circ Q(f) \circ \left(\alpha_{M', \overline{M}}^{-1} \alpha_{M, \overline{M}}\right)^{-1}.
\]
Substituting the identity $\alpha_{X, \overline{X}} = \beta_{X, \overline{X}} \theta_X$ into the arrow condition for $\Phi_\alpha$, we obtain:
\begin{align*}
Q(f') &= \left(\beta_{N', \overline{N}} \theta_{N'}\right)^{-1} \left(\beta_{N, \overline{N}} \theta_N\right) \circ Q(f) \circ \left( \left(\beta_{M', \overline{M}} \theta_{M'}\right)^{-1} \left(\beta_{M, \overline{M}} \theta_M\right) \right)^{-1} \\
&= \theta_{N'}^{-1} \left(\beta_{N', \overline{N}}^{-1} \beta_{N, \overline{N}}\right) \left(\theta_N Q(f) \theta_M^{-1}\right) \left(\beta_{M, \overline{M}}^{-1} \beta_{M', \overline{M}}\right) \theta_{M'} \\
&= \theta_{N'}^{-1} \beta_{N, N'} Q(f) \beta_{M, M'}^{-1} \theta_{M'},
\end{align*}
where we have applied the previously justified naturality $\theta_N Q(f) \theta_M^{-1} = Q(f)$. Furthermore, because $\theta_{N'}$ and $\theta_{M'}$ are fixed automorphisms determined by the unique representatives in the transversal, the relation remains invariant. This confirms that $(f,f') \in \Phi^{\mathrm{fl}}_{\alpha\scalebox{0.7}{$\left(\begin{smallmatrix}M&N\\M'& N'\end{smallmatrix}\right)$}}$ if, and only if, $(f,f') \in \Phi^{\mathrm{fl}}_{\beta\scalebox{0.7}{$\left(\begin{smallmatrix}M&N\\M'& N'\end{smallmatrix}\right)$}}$. Consequently, $F_{\alpha, \beta}$ is a well-defined functor and, being a bijection on both objects and morphisms, it constitutes an isomorphism of categories.
\end{example}

\subsection{A Riguet Congruence on the Category associated with a Deterministic Finite Automata}\label{S:automata}

\[
\]

\begin{example}
For a deterministic finite automaton (DFA) $\mathcal{A}$, after defining a Riguet congruence $\Phi$ on $\mathsf{C}(\mathcal{A})$, the category associated to $\mathcal{A}$, and the action category $\mathsf{Act}(\mathbf{T}(\mathcal{A}/{\sim}), Q/{\sim})$, we define a functor $F\colon\mathsf{C}(\mathcal{A})\mor \mathsf{Act}(\mathbf{T}(\mathcal{A}/{\sim}), Q/{\sim})$ and verify that $F$ satisfies conditions~(1) and~(2) of Proposition~\ref{UPQ}. Hence, Proposition~\ref{UPQ} yields a unique functor $F^{\sharp}\colon \mathsf{C}(\mathcal{A})/\Phi \mor \mathsf{Act}(\mathbf{T}(\mathcal{A}/{\sim}), Q/{\sim})$ such that $F = F^{\sharp} \circ P_{\Phi}$. Moreover, $F^{\sharp}$ is an isomorphism. 

We begin by introducing the notions and constructions required to achieve this goal. 
Let $\mathcal{A} = (Q, \Sigma, \delta, q_{0}, F)$ be a DFA and $\Sigma^{\star}$ the underlying set of $\boldsymbol{\Sigma}^{\star}$, the free monoid on $\Sigma$. For $\boldsymbol{\Sigma}^{\star}$, as previously, we denote by $\lambda$ the neutral element of $\boldsymbol{\Sigma}^{\star}$, by $\curlywedge$ the operation of concatenation of $\boldsymbol{\Sigma}^{\star}$ and by $\eta_{\Sigma}$, the canonical embedding of $\Sigma$ into $\Sigma^{\star}$ that sends $a\in \Sigma$ to $(a)$, usually identified with $a$ itself.  Before defining $\mathsf{C}(\mathcal{A})$ we recall that from the transition mapping $\delta$ from $Q\times \Sigma$ to $Q$ we get the extended transition mapping $\delta^{\star}$ from $Q\times \Sigma^{\star}$ to $Q$ defined by left recursion, for every $q\in Q$, every $a\in \Sigma$ and every $\mathbf{w}\in \Sigma^{\star}$ as: $\delta^{\star}(q,\lambda) = q$; and $\delta^{\star}(q,(a)\curlywedge \mathbf{w}) = \delta^{\star}(\delta(q,a),\mathbf{w})$. The mapping $\delta^{\star}$ has the following property: For every $q\in Q$ and every $\mathbf{u},\mathbf{v}\in\Sigma^{\star}$,
$\delta^{\star}(q,\mathbf{u}\curlywedge \mathbf{v})=\delta^{\star}(\delta^{\star}(q,\mathbf{u}),\mathbf{v})$. Its proof is by induction on the length of $\mathbf{u}$. The extended transition mapping $\delta^{\star}$ may equivalently be defined by right recursion, replacing the second clause above with $\delta^{\star}(q,\mathbf{w}\curlywedge (a)) = \delta(\delta^{\star}(q,\mathbf{w}),a)$.

We define $\mathsf{C}(\mathcal{A})$ as the category where $\mathrm{Ob}(\mathsf{C}(\mathcal{A})) = Q$, and, for $p$, $q\in Q$, the set $\mathrm{Hom}_{\mathsf{C}(\mathcal{A})}(p, q)$ of the morphisms from $p$ to $q$ is 
\[
\{p\}\times \{\mathbf{w} \in \Sigma^{\star}\mid \delta^{\star}(p, \mathbf{w}) = q\}\times \{q\}.
\] 
Thus, a morphism from $p$ to $q$ is an ordered triple $(p, \mathbf{w}, q)$ in which $\mathbf{w} \in \Sigma^{\star}$ such that $\delta^{\star}(p, \mathbf{w}) = q$. The composition of $(p, \mathbf{u}, q)$ and $(q, \mathbf{v}, r)$ is 
$(q, \mathbf{v}, r) \circ (p, \mathbf{u}, q) = (p, \mathbf{u} \curlywedge \mathbf{v}, r)$, and the identity at $q$ is $\mathrm{id}_{q} = (q, \lambda, q)$. Note that the category $\mathsf{C}(\mathcal{A})$ is the free category generated by $\mathsf{G}(\mathcal{A})$, where  $\mathsf{G}(\mathcal{A}) = (V, E, \mathrm{sc}, \mathrm{tg})$ is the quiver associate with $\mathcal{A}$ in which $V = Q$, $E = \{ (p, a, q) \in Q \times \Sigma \times Q \mid \delta(p, a) = q \}$, the source and target maps are $\mathrm{sc}(p, a, q) = p$ and $\mathrm{tg}(p, a, q) = q$. 

We denote by  $\mathcal{A}/{\sim} = (Q/{\sim}, \Sigma, \delta_{\sim}, [q_{0}]_{\sim}, F/{\sim})$ the minimal quotient automaton of  $\mathcal{A}$, where $\sim$ is the Myhill-Nerode congruence on $\mathcal{A}$,  defined, for every pair $(p, q) \in Q \times Q$, by
\[
(p, q) \in \sim \iff \forall \mathbf{w} \in \Sigma^{\star} (\delta^{\star}(p, \mathbf{w}) \in F \iff \delta^{\star}(q, \mathbf{w}) \in F).
\] 
Note that, for every $q\in Q$ and every $a\in \Sigma$, 
$\delta_{\sim}([q]_{\sim},a) = [\delta(q,a)]_{\sim}$. From the mapping $\delta_{\sim}$ we obtain the mapping $(\delta_{\sim}(\bigcdot,a))_{a\in \Sigma}$ from $\Sigma$ to $\mathrm{End}(Q/{\sim})$. Then, by the universal property of the free monoid on a set, there exists a unique homomorphism $(\delta_{\sim}(\bigcdot,a))_{a\in \Sigma}^{\sharp}$ from $\boldsymbol{\Sigma}^{\star}$ to $\mathbf{End}(Q/{\sim})$ such that 
$(\delta_{\sim}(\bigcdot,a))_{a\in \Sigma}^{\sharp}\circ \eta_{\Sigma} = (\delta_{\sim}(\bigcdot,a))_{a\in \Sigma}$. We denote by $\mathbf{T}(\mathcal{A}/{\sim})$ the image of $(\delta_{\sim}(\bigcdot,a))_{a\in \Sigma}^{\sharp}$ and we call it the transition monoid of $\mathcal{A}/{\sim}$. Thus $T(\mathcal{A}/{\sim})$, the underlying set of  $\mathbf{T}(\mathcal{A}/{\sim})$, consists of those mappings $f$ from $Q/{\sim}$ to $Q/{\sim}$ for which there exists a word $\mathbf{w} \in \Sigma^{\star}$ such that $f = f_{\mathbf{w}}$, where $f_{\mathbf{w}}$ is the mapping from $Q/{\sim}$ to $Q/{\sim}$ that sends $[q]_{\sim} \in Q/{\sim}$ to $f_{\mathbf{w}}([q]_{\sim}) = [\delta^{\star}(q, \mathbf{w})]_{\sim}$. Equivalently, $\mathbf{T}(\mathcal{A}/{\sim})$ is the submonoid of $\mathbf{End}(Q/{\sim})$ generated by $\{\delta_{\sim}(\bigcdot,a)\mid a\in \Sigma\}$. Moreover, we denote by $\mathbf{K}$ the monoid whose underlying set is $K=\{\mathbf w\in\Sigma^{\star}\mid f_{\mathbf w}=\mathrm{id}_{Q/\sim}\}$, i.e., the kernel of $(\delta_{\sim}(\bigcdot,a))_{a\in \Sigma}^{\sharp}$.

From now on, for every $p \in Q$ and every $(\mathbf{w},r)\in K\times [p]_{\sim}$, we assume that $\delta^{\star}(r,\mathbf w)\in [p]_{\sim}$ and that the mapping $\bigcdot$ from $K\times [p]_{\sim}$ to $[p]_{\sim}$ that sends $(\mathbf{w},r)\in K\times [p]_{\sim}$ to $\mathbf{w} \bigcdot r=\delta^{\star}(r,\mathbf w)$ acts transitively on $[p]_{\sim}$. Recall that the transitivity means that, for every pair of states $r,s\in [p]_{\sim}$, there exists a word $\mathbf w\in K$ such that $\delta^{\star}(r,\mathbf w)=s$. Note that the transitivity of the $\mathbf{K}$-action on the equivalence classes implies that, for every $p\in Q$, $[p]_{\sim} = \mathrm{Orb}_{\mathbf{K}}(p)$.

Under this assumption we have that, for every $p, q \in Q$, if $p \sim q$, then there exists a word $\mathbf{z} \in \Sigma^{\star}$ such that $\delta^{\star}(p, \mathbf{z}) = q$ and $\mathbf{z} \in K$, where the latter condition is, by definition of $K$, equivalent to $f_{\mathbf{z}} = \mathrm{id}_{Q/{\sim}}$.

Then we define the Riguet congruence $\Phi = (\Phi^{\mathrm{ob}}, \Phi^{\mathrm{fl}})$ on $\mathsf{C}(\mathcal{A})$ as follows: 
\begin{enumerate}
\item \textsf{On Objects:} $(p, q) \in \Phi^{\mathrm{ob}}$ if, and only if, $p$ and $q$ are strongly bisimilar, i.e., $\Phi^{\mathrm{ob}} = {\sim}$.
\item \textsf{On Arrows:} For every morphism $(p,\mathbf{w},q)$ in $\mathrm{Hom}_{\mathsf{C}(\mathcal{A})}(p, q)$ and 
every morphism $(p',\mathbf{w}',q')$ in $\mathrm{Hom}_{\mathsf{C}(\mathcal{A})}(p', q')$, with 
$(p, p') \in \Phi^{\mathrm{ob}}$ and $(q, q') \in \Phi^{\mathrm{ob}}$, we declare 
$((p,\mathbf{w},q), (p',\mathbf{w}',q'))\in \Phi^{\mathrm{fl}}_{\scalebox{0.7}{$\left(\begin{smallmatrix}p&q\\p'& q'\end{smallmatrix}\right)$}}$ if, and only if, for every $r \in Q$, we have that $[\delta^{\star}(r, \mathbf{w})]_\sim = [\delta^{\star}(r, \mathbf{w}')]_\sim$, or, what is equivalent, if, and only if,  $f_{\mathbf{w}} = f_{\mathbf{w}'}$ in $\mathbf{T}(\mathcal{A}/{\sim})$.
\end{enumerate}

By definition, $\Phi^{\mathrm{ob}} = {\sim}$ and $\sim$ is a congruence on $\mathcal{A}$, hence 
$\Phi^{\mathrm{ob}}$ is an equivalence relation on $\mathrm{Ob}(\mathsf{C}(\mathcal{A})) = Q$. The proof of conditions $(\mathrm{a})$ through $(\mathrm{e})$ is immediate from the definition of $\Phi^{\mathrm{fl}}$ and is omitted. 

(f) Let $(p, p') \in \sim$ and $(q, q') \in \sim$. Given $(p, \mathbf{w}, q) \in \mathrm{Hom}_{\mathsf{C}(\mathcal{A})}(p, q)$, we must find $(p', \mathbf{w}', q') \in \mathrm{Hom}_{\mathsf{C}(\mathcal{A})}(p', q')$ such that $f_{\mathbf{w}} = f_{\mathbf{w}'}$.
    
Since $p \sim p'$, we have that $f_{\mathbf{w}}([p]_{\sim}) = f_{\mathbf{w}}([p']_{\sim})$. Moreover, since $(p, \mathbf{w}, q)$ is a morphism from $p$ to $q$, we have that $\delta^{\star}(p, \mathbf{w}) = q$. Hence $[q]_{\sim} = f_{\mathbf{w}}([p]_{\sim})$. By the definition of $f_{\mathbf{w}}$, we have that $\delta^{\star}(p', \mathbf{w}) \in [q]_{\sim}$. Let $q'' = \delta^{\star}(p', \mathbf{w})$. Since $q'' \in [q]_{\sim}$ and $q \sim q'$, by transitivity $q'' \sim q'$. Then, by our previous results, for $q''$ and $q'$ we have that there exists a word $\mathbf{z} \in K$ such that $\delta^{\star}(q'', \mathbf{z}) = q'$ and $f_{\mathbf{z}} = \mathrm{id}_{Q/{\sim}}$. Define $\mathbf{w}' = \mathbf{w} \curlywedge \mathbf{z}$. Then:
 $\delta^{\star}(p', \mathbf{w}') = \delta^{\star}(\delta^{\star}(p', \mathbf{w}), \mathbf{z}) = \delta^{\star}(q'', \mathbf{z}) = q'$. Thus, $(p', \mathbf{w}', q') \in \mathrm{Hom}_{\mathsf{C}(\mathcal{A})}(p', q')$. In $\mathbf{T}(\mathcal{A}/{\sim})$ we have that $f_{\mathbf{w}'} = f_{\mathbf{z}} \circ f_{\mathbf{w}}$. Since $f_{\mathbf{z}} = \mathrm{id}_{Q/{\sim}}$, it follows  that $f_{\mathbf{w}'} = \mathrm{id}_{Q/{\sim}} \circ f_{\mathbf{w}} = f_{\mathbf{w}}$. This confirms that $((p, \mathbf{w}, q), (p', \mathbf{w}', q')) \in \Phi^{\mathrm{fl}}_{\scalebox{0.7}{$\left(\begin{smallmatrix}p&q\\p'& q'\end{smallmatrix}\right)$}}$, satisfying the quadrangular completion condition.

Since all conditions (a) through (f) are satisfied, $\Phi = (\sim, \Phi^{\mathrm{fl}})$ is a Riguet congruence on $\mathsf{C}(\mathcal{A})$.

We now establish that the quotient category $\mathsf{C}(\mathcal{A})/\Phi$ is isomorphic to the action category $\mathsf{Act}(\mathbf{T}(\mathcal{A}/{\sim}), Q/{\sim})$. The category $\mathsf{Act}(\mathbf{T}(\mathcal{A}/{\sim}), Q/{\sim})$ ha as objects the elements of $Q/{\sim}$, and as morphisms from $[s]_\sim$ to $[t]_\sim$ those ordered triples $([s]_\sim, f_{\mathbf{w}}, [t]_\sim)$ for which $f_{\mathbf{w}}\in T(\mathcal{A}/{\sim})$ and $f_{\mathbf{w}}([s]_\sim) = [t]_\sim$. We define a functor $F\colon \mathsf{C}(\mathcal{A})\mor \mathsf{Act}(\mathbf{T}(\mathcal{A}/{\sim}), Q/{\sim})$ as follows: For every state $p \in \mathrm{Ob}(\mathsf{C}(\mathcal{A}))$, we set:
$F(p) = [p]_{\sim}$. For every morphism $(p, \mathbf{w}, q) \in \mathrm{Hom}_{\mathsf{C}(\mathcal{A})}(p, q)$, we set: $F(p, \mathbf{w}, q) = ([p]_{\sim}, f_{\mathbf{w}}, [q]_{\sim})$. 
Note that $F(p, \mathbf{w}, q)$ is a well-defined morphism in the action category because $\delta^{\star}(p, \mathbf{w}) = q$ implies $f_{\mathbf{w}}([p]_{\sim}) = [\delta^{\star}(p, \mathbf{w})]_{\sim} = [q]_{\sim}$. 

We next verify for $F$ the conditions in Proposition~\ref{UPQ}. (1) If $(p, p') \in \Phi^{\mathrm{ob}}$, then $p \sim p'$, which implies $[p]_{\sim} = [p']_{\sim}$. Thus, $F(p) = F(p')$.
(2) Let $(p, \mathbf{w}, q) \in \mathrm{Hom}_{\mathsf{C}(\mathcal{A})}(p, q)$ and $(p', \mathbf{w}', q') \in \mathrm{Hom}_{\mathsf{C}(\mathcal{A})}(p', q')$. If $((p, \mathbf{w}, q), (p', \mathbf{w}', q')) \in \Phi^{\mathrm{fl}}_{\scalebox{0.7}{$\left(\begin{smallmatrix}p&q\\p'& q'\end{smallmatrix}\right)$}}$, then by definition, $[p]_{\sim} = [p']_{\sim}$ and $[q]_{\sim} = [q']_{\sim}$ (since $(p, p') \in \Phi^{\mathrm{ob}}$ and $(q, q') \in \Phi^{\mathrm{ob}}$). Hence $f_{\mathbf{w}} = f_{\mathbf{w}'}$ in $\mathbf{T}(\mathcal{A}/{\sim})$. Consequently, $F(p, \mathbf{w}, q) = ([p]_{\sim}, f_{\mathbf{w}}, [q]_{\sim}) = ([p']_{\sim}, f_{\mathbf{w}'}, [q']_{\sim}) = F(p', \mathbf{w}', q')$. 
Hence, by Proposition~\ref{UPQ}, there exists a unique functor 
$F^{\sharp} \colon \mathsf{C}(\mathcal{A})/\Phi \mor\mathsf{Act}(\mathbf{T}(\mathcal{A}/{\sim}), Q/{\sim})$
such that $F^{\sharp} \circ P_{\Phi} = F$. This functor is defined by $F^{\sharp}([p]_{\Phi^{\mathrm{ob}}}) = [p]_{\sim}$ and $F^{\sharp}([(p, \mathbf{w}, q)]_{\Phi^{\mathrm{fl}}}) = ([p]_{\sim}, f_{\mathbf{w}}, [q]_{\sim})$.

We next prove that $F^{\sharp}$ is an isomorphism of categories. The object mapping of $F^{\sharp}$ is bijective because it acts as the identity on the set of equivalence classes $Q/{\sim}$.

To prove that $F^{\sharp}\colon \mathrm{Hom}_{\mathsf{C}(\mathcal{A})/\Phi}([p]_{\sim},[q]_{\sim})\mor
\mathrm{Hom}_{\mathsf{Act}(\mathbf{T}(\mathcal{A}/{\sim}), Q/{\sim})}([p]_{\sim},[q]_{\sim})$  is a bijection, we show that it is both full and faithful. 

\textsf{Faithfulness.} Let $[p]_{\sim}$ and $[q]_{\sim}$ be two objects in $\mathsf{C}(\mathcal{A})/\Phi$ and  
$[(p_{0}, \mathbf{u}, q_{0})]_{\Phi^{\mathrm{fl}}}$, $[(p_{1}, \mathbf{v}, q_{1})]_{\Phi^{\mathrm{fl}}}$ two morphisms from $[p]_{\sim}$ to $[q]_{\sim}$. Then, by the definition of morphisms in the quotient category, we must have that $[p_{0}]_{\sim} = [p_{1}]_{\sim} = [p]_{\sim}$ and $[q_{0}]_{\sim} = [q_{1}]_{\sim} = [q]_{\sim}$, which implies that $(p_{0}, p_{1}) \in \Phi^{\mathrm{ob}}$ and $(q_{0}, q_{1}) \in \Phi^{\mathrm{ob}}$. 

Now, assume that $F^{\sharp}([(p_{0}, \mathbf{u}, q_{0})]_{\Phi^{\mathrm{fl}}}) = F^{\sharp}([(p_{1}, \mathbf{v}, q_{1})]_{\Phi^{\mathrm{fl}}})$. Then, by the definition of $F^{\sharp}$, we have that 
$
([p_{0}]_{\sim}, f_{\mathbf{u}}, [q_{0}]_{\sim}) = ([p_{1}]_{\sim}, f_{\mathbf{v}}, [q_{1}]_{\sim}). 
$
This equality in $\mathsf{Act}(\mathbf{T}(\mathcal{A}/{\sim}), Q/{\sim})$ implies, in particular, that $f_{\mathbf{u}} = f_{\mathbf{v}}$ in $\mathbf{T}(\mathcal{A}/{\sim})$.

Then, according to the definition of $\Phi^{\mathrm{fl}}$, $((p_{0}, \mathbf{u}, q_{0}), (p_{1}, \mathbf{v}, q_{1})) \in \Phi^{\mathrm{fl}}_{\scalebox{0.7}{$\left(\begin{smallmatrix}p_{1}&q_{1}\\p_{2}& q_{2}\end{smallmatrix}\right)$}}$ if, and only if, 
\begin{enumerate}
\item $(p_{0}, p_{1}) \in \Phi^{\mathrm{ob}}$ (which holds since $[p_{0}]_{\sim} = [p_{1}]_{\sim}$),
\item $(q_{0}, q_{1}) \in \Phi^{\mathrm{ob}}$ (which holds since $[q_{0}]_{\sim} = [q_{1}]_{\sim}$), and
\item $f_{\mathbf{u}} = f_{\mathbf{v}}$ (which holds since $F^{\sharp}([(p_{0}, \mathbf{u}, q_{0})]_{\Phi^{\mathrm{fl}}}) = F^{\sharp}([(p_{1}, \mathbf{v}, q_{1})]_{\Phi^{\mathrm{fl}}})$).
\end{enumerate}
Since all conditions are satisfied, $(p_{0}, \mathbf{u}, q_{0})$ and $(p_{1}, \mathbf{v}, q_{1})$ belong to the same equivalence class, i.e., $[(p_{0}, \mathbf{u}, q_{0})]_{\Phi^{\mathrm{fl}}} = [(p_{1}, \mathbf{v}, q_{1})]_{\Phi^{\mathrm{fl}}}$. Thus, $F^{\sharp}$ is faithful.

\textsf{Fullness.}  Let $([s]_{\sim}, f_{\mathbf{w}}, [t]_{\sim})$ be a morphism in the action category. Then we have $\delta^{\star}(s, \mathbf{w}) = t'$ for some $t' \sim t$. Hence, there exists a word $\mathbf{z}\in K$ such that $\delta^{\star}(t', \mathbf{z}) = t$ and $f_{\mathbf{z}} = \mathrm{id}_{Q/{\sim}}$. Taking $\mathbf{w}' = \mathbf{w} \curlywedge \mathbf{z}$, we have $\delta^{\star}(s, \mathbf{w}') = t$ and $f_{\mathbf{w}'} = f_{\mathbf{w}}$. Then $[(s, \mathbf{w}', t)]_{\Phi^{\mathrm{fl}}}$ is the preimage of $([s]_{\sim}, f_{\mathbf{w}}, [t]_{\sim})$.

Therefore, $F^{\sharp}$ is an isomorphism of categories:
\[ 
\mathsf{C}(\mathcal{A})/\Phi \cong \mathsf{Act}(\mathbf{T}(\mathcal{A}/{\sim}), Q/{\sim}). 
\]

\end{example}

\subsection{A Riguet Congruence on the Category of Sets and Relations}\label{S:rel}

\begin{example}
After defining a suitable Riguet congruence $\Phi$ on the category $\mathsf{Rel}$ of sets and relations, and the category $\mathsf{Card}_{\mathrm{orb}}$ of cardinal numbers $\kappa \in \mathbf{U}$ and orbits, we define a functor $F\colon \mathsf{Rel} \mor \mathsf{Card}_{\mathrm{orb}}$ and verify that $F$ satisfies conditions~(1) and~(2) of Proposition~\ref{UPQ}. Hence, Proposition~\ref{UPQ} yields a unique functor $F^{\sharp}\colon \mathsf{Rel}/\Phi \mor \mathsf{Card}_{\mathrm{orb}}$ such that $F = F^{\sharp} \circ P_{\Phi}$. Moreover, $F^{\sharp}$ is an isomorphism.

We define a Riguet congruence $\Phi = (\Phi^{\mathrm{ob}}, \Phi^{\mathrm{fl}})$ on $\mathsf{Rel}$ as follows:
\begin{enumerate}
\item \textsf{On Objects:} $(X,Y) \in \Phi^{\mathrm{ob}}$ if, and only if, $X \cong Y$ in $\mathsf{Rel}$.
\item \textsf{On Arrows:} For every $R \in \mathrm{Hom}_{\mathsf{Rel}}(X, Y)$ and $R' \in \mathrm{Hom}_{\mathsf{Rel}}(X', Y')$ with $(X, X') \in \Phi^{\mathrm{ob}}$ and $(Y, Y') \in \Phi^{\mathrm{ob}}$ we have that $(R, R') \in \Phi^{\mathrm{fl}}_{\scalebox{0.7}{$\left(\begin{smallmatrix}X&Y\\X'& Y'\end{smallmatrix}\right)$}}$ if, and only if, there exists isomorphisms $\sigma \colon X \mor X'$ and $\tau \colon Y\mor Y'$ such that
$R' = \tau \circ R \circ \sigma^{-1}$.
\end{enumerate}

By definition, $\Phi^{\mathrm{ob}}$ is an equivalence relation on $\mathrm{Ob}(\mathsf{Rel})$. The proof of conditions $(\mathrm{a})$ through $(\mathrm{f})$ is immediate from the definition of $\Phi^{\mathrm{fl}}$ and is omitted. 

We next prove that the quotient category $\mathsf{Rel}/\Phi$ is isomorphic to the category $\mathsf{Card}_{\mathrm{orb}}$, defined as follows. The objects of $\mathsf{Card}_{\mathrm{orb}}$ are the cardinal numbers $\kappa \in \mathbf{U}$. For cardinals $\kappa, \lambda \in \mathbf{U}$, the morphisms from $\kappa$ to $\lambda$ are the orbits of the set $\mathrm{Sub}(\kappa \times \lambda)$ under the action of the product of symmetric groups $\boldsymbol{\mathfrak{S}}_{\kappa} \times \boldsymbol{\mathfrak{S}}_{\lambda}$. Specifically, for every $(\sigma, \tau) \in \mathfrak{S}_\kappa \times \mathfrak{S}_\lambda$ and $R \in \mathrm{Sub}(\kappa \times \lambda)$, the action is defined by:
$(\sigma, \tau) \bigcdot R = \tau \circ R \circ \sigma^{-1}$, and $\mathrm{Orb}(R) = \{ \tau \circ R \circ \sigma^{-1} \mid (\sigma, \tau) \in \mathfrak{S}_{\kappa} \times \mathfrak{S}_{\lambda} \}$, the orbit of $R$ under this group action, is a morphism from $\kappa$ to $\lambda$. The identity at $\kappa$ is $\mathrm{id}_{\kappa} = \mathrm{Orb}(\Delta_{\kappa})$. For every $\mu \in \mathbf{U}$ and $S \in \mathrm{Sub}(\lambda \times \mu)$, the composition is given by $\mathrm{Orb}(S) \circ \mathrm{Orb}(R) = \mathrm{Orb}(S \circ R)$. 
The composition in $\mathsf{Card}_{\mathrm{orb}}$ is well-defined because the relational composition $S \circ R$ is equivariant with respect to the action of the symmetric groups. To verify this, let $R, R' \in \mathrm{Sub}(\kappa \times \lambda)$ and $S, S' \in \mathrm{Sub}(\lambda \times \mu)$ be such that $\mathrm{Orb}(R) = \mathrm{Orb}(R')$ and $\mathrm{Orb}(S) = \mathrm{Orb}(S')$. By the definition of the action, there exist permutations $(\sigma, \tau) \in \mathfrak{S}_{\kappa} \times \mathfrak{S}_{\lambda}$ and $(\rho, \theta) \in \mathfrak{S}_{\lambda} \times \mathfrak{S}_{\mu}$ such that $R' = \tau \circ R \circ \sigma^{-1}$ and $S' = \theta \circ S \circ \rho^{-1}$. The composition of $R'$ and $S'$ yields: 
\[
S' \circ R' = (\theta \circ S \circ \rho^{-1}) \circ (\tau \circ R \circ \sigma^{-1}) = \theta \circ (S \circ \eta \circ R) \circ \sigma^{-1},
\]
where $\eta = \rho^{-1} \circ \tau \in \mathfrak{S}_{\lambda}$ is an internal permutation. Specifically, any such internal permutation satisfies:
\[
\mathrm{Orb}(S \circ \eta \circ R) = \mathrm{Orb}(S \circ R).
\]
This ensures that $S' \circ R'$, being a transform of $S \circ \eta \circ R$ by the external permutations $\theta$ and $\sigma^{-1}$, will always lie in $\mathrm{Orb}(S \circ R)$, regardless of the specific choices of $\sigma, \tau, \rho$, and $\theta$. Thus, the composition is well-defined.

We define a functor $F \colon \mathsf{Rel} \mor \mathsf{Card}_{\mathrm{orb}}$ as follows. For every object $X$ of $\mathsf{Rel}$, we set: $F(X) = \mathrm{card}(X)$. To define $F$ on morphisms, we fix for each set $X$ an isomorphism $i_{X} \colon \mathrm{card}(X) \mor X$. For every morphism $R \colon X \mor Y$ of $\mathsf{Rel}$, we first define its lifting to the cardinal level $\overline{R} \in \mathrm{Sub}(\mathrm{card}(X) \times \mathrm{card}(Y))$ as:
\[
\overline{R} = (i_{X} \times i_{Y})^{-1}[R].
\]
We then set: $F(R) = \mathrm{Orb}(\overline{R})$. 

The functor $F$ is such that, for every $X$ and $Y$, if $X\cong Y$, then $F(X) = F(Y)$. Moreover, for every $R\colon X\mor Y$ and every $R'\colon X'\mor Y'$, if $X\cong X'$, $Y\cong Y'$ and $(R, R') \in \Phi^{\mathrm{fl}}_{\scalebox{0.7}{$\left(\begin{smallmatrix}X&Y\\X'& Y'\end{smallmatrix}\right)$}}$, then $\mathrm{Orb}(\overline{R}) = \mathrm{Orb}(\overline{R'})$. Fix bijections $j_{X'} \colon \mathrm{card}(X') \mor X'$ and $j_{Y'} \colon \mathrm{card}(Y') \mor Y'$. The liftings are:
\[
\overline{R} = (i_{X} \times i_{Y})^{-1}[R] \quad \text{and} \quad \overline{R'} = (j_{X'} \times j_{Y'})^{-1}[R'].
\] 
By hypothesis, there exist bijections $\sigma \colon X \mor X'$ and $\tau \colon Y \mor Y'$ such that $R' = \tau \circ R \circ \sigma^{-1}$. Substituting this into the expression for $\overline{R'}$ we obtain: 
\[
\overline{R'} = (j_{X'} \times j_{Y'})^{-1}[\tau \circ R \circ \sigma^{-1}].
\] 
Since for any relation $S \subseteq V \times W$ and bijections $f$ from $\mathrm{card}(V)$ to $V$ and $g$ from  $\mathrm{card}(W)$ to $W$, we have that $(f \times g)^{-1}[S] = g^{-1} \circ S \circ f$, we can rewrite $\overline{R'}$ as a composition of relations. By inserting the identities $\mathrm{id}_Y = i_Y \circ i_Y^{-1}$ and $\mathrm{id}_X = i_X \circ i_X^{-1}$, we obtain:
\[\overline{R'} = (j_{Y'}^{-1} \circ \tau \circ i_{Y}) \circ \underbrace{(i_{Y}^{-1} \circ R \circ i_{X})}_{\overline{R}} \circ (i_{X}^{-1} \circ \sigma^{-1} \circ j_{X'}),\]
\[\overline{R'} = (j_{Y'}^{-1} \circ \tau \circ i_{Y}) \circ \overline{R} \circ (j_{X'}^{-1} \circ \sigma \circ i_{X})^{-1}.
\]
Let $\hat{\tau} = j_{Y'}^{-1} \circ \tau \circ i_{Y} \in \mathfrak{S}_{\mathrm{card}(Y)}$ and $\hat{\sigma} = j_{X'}^{-1} \circ \sigma \circ i_{X} \in \mathfrak{S}_{\mathrm{card}(X)}$. Since $\overline{R'} = \hat{\tau} \circ \overline{R} \circ \hat{\sigma}^{-1}$, it follows that $\overline{R'}$ and $\overline{R}$ lie in the same orbit under the action of $\mathfrak{S}_{\mathrm{card}(X)} \times \mathfrak{S}_{\mathrm{card}(Y)}$. Thus, $\mathrm{Orb}(\overline{R}) = \mathrm{Orb}(\overline{R'})$.  Therefore, by Proposition~\ref{UPQ}, the universal property of the quotient, there exists a unique functor $F^{\sharp}$ from $\mathsf{Rel}/\Phi$ to $\mathsf{Card}_{\mathrm{orb}}$ such that $F = F^{\sharp}\circ P_{\Phi}$. 

Next we prove that $F^{\sharp}$ is an isomorphism. The objects of $\mathsf{Rel}/\Phi$ are equivalence classes $[X]_{\Phi^{\mathrm{ob}}}$, which correspond to the set of all sets in $\mathbf{U}$ with the same cardinality. Since the objects of $\mathsf{Card}_{\mathrm{orb}}$ are the cardinal numbers in $\mathbf{U}$, the mapping $F^{\sharp}([X]) = \mathrm{card}(X)$ is a bijection on objects.

To prove that $F^{\sharp} \colon \mathrm{Hom}_{\mathsf{Rel}/\Phi}([X]_{\Phi^{\mathrm{ob}}}, [Y]_{\Phi^{\mathrm{ob}}}) \mor \mathrm{Hom}_{\mathsf{Card}_{\mathrm{orb}}}(\kappa, \lambda)$ is a bijection, we show that it is both full and faithful. 

\textsf{Fullness.} 
Let $\mathrm{Orb}(R) \in \mathrm{Hom}_{\mathsf{Card}_{\mathrm{orb}}}(\kappa, \lambda)$ be an arbitrary morphism, where $R \in \mathrm{Sub}(\kappa \times \lambda)$. We define the relation $S \in \mathrm{Sub}(X \times Y)$ as the transport $S = (i_{X} \times i_{Y})[R]$. Its lifting to the cardinal level is $\overline{S} = (i_{X} \times i_{Y})^{-1}[S] = R$. Thus, $F(S) = \mathrm{Orb}(\overline{S}) = \mathrm{Orb}(R)$. Since $F = F^{\sharp} \circ P_{\Phi}$, it follows that $F^{\sharp}([S]_{\Phi^{\mathrm{fl}}}) = \mathrm{Orb}(R)$.
 
\textsf{Faithfulness.}
Suppose $F^{\sharp}([R]_{\Phi^{\mathrm{fl}}}) = F^{\sharp}([S]_{\Phi^{\mathrm{fl}}})$ for two relations $R, S \subseteq X \times Y$. This implies $\mathrm{Orb}(\overline{R}) = \mathrm{Orb}(\overline{S})$. Then there exist permutations $\hat{\sigma} \in \mathfrak{S}_{\kappa}$ and $\hat{\tau} \in \mathfrak{S}_{\lambda}$ such that $\overline{S} = \hat{\tau} \circ \bar{R} \circ \hat{\sigma}^{-1}$. Applying the transport $(i_{X} \times i_{Y})$ to both sides, we obtain $S = (i_{Y} \circ \hat{\tau} \circ i_{Y}^{-1}) \circ R \circ (i_{X} \circ \hat{\sigma} \circ i_{X}^{-1})^{-1}$. Setting $\sigma = i_{X} \circ \hat{\sigma} \circ i_{X}^{-1}$ and $\tau = i_{Y} \circ \hat{\tau} \circ i_{Y}^{-1}$, we have $S = \tau \circ R \circ \sigma^{-1}$, hence $[R]_{\Phi^{\mathrm{fl}}} = [S]_{\Phi^{\mathrm{fl}}}$.

Since $F^{\sharp}$ is a bijection on objects and a bijection on every hom-set, it is an isomorphism of categories $\mathsf{Rel}/\Phi \cong \mathsf{Card}_{\mathrm{orb}}$.

The construction of the congruence $\Phi$ reflects the principle of \textit{difunctionality}---a concept introduced by Riguet \cite{JR48} for binary relations via the inclusion $R \circ R^{-1} \circ R \subseteq R$, and, to the best of our knowledge, independently by Mac Lane \cite{sM61} in the study of additive relations via the identity $f \circ f^{\sharp} \circ f = f$ (note that for binary relations, the inclusion $R\subseteq R \circ R^{-1} \circ R$ always holds). In the present categorical setting, this property is manifested within each component 
$\Phi^{\mathrm{fl}}_{\scalebox{0.7}{$\left(\begin{smallmatrix}X&Y\\X'& Y'\end{smallmatrix}\right)$}}$.

While the well-definedness of the quotient $\mathsf{Rel}/\Phi$ is established by the equivariance of relational composition under the action of the symmetric groups, the difunctionality of the component relation
$\Phi^{\mathrm{fl}}_{\scalebox{0.7}{$\left(\begin{smallmatrix}X&Y\\X'& Y'\end{smallmatrix}\right)$}} \circ 
\left(\Phi^{\mathrm{fl}}_{\scalebox{0.7}{$\left(\begin{smallmatrix}X&Y\\X'& Y'\end{smallmatrix}\right)$}}\right)^{-1} \circ 
\Phi^{\mathrm{fl}}_{\scalebox{0.7}{$\left(\begin{smallmatrix}X&Y\\X'& Y'\end{smallmatrix}\right)$}} \subseteq 
\Phi^{\mathrm{fl}}_{\scalebox{0.7}{$\left(\begin{smallmatrix}X&Y\\X'& Y'\end{smallmatrix}\right)$}}$
ensures that the transport of relations is independent of the choice of isomorphism witnesses $(\sigma, \tau)$. This structural property allows the Riguet congruence to perform a \textit{skeletization of morphisms}, collapsing the hom-sets of $\mathsf{Rel}$ into the orbits of $\mathsf{Card}_{\mathrm{orb}}$. This captures the pure algebraic essence of relations.

Let us note that the formalism developed for the category $\mathsf{Rel}$ and its quotient $\mathsf{Card}_{\mathrm{orb}}$ extends naturally to the category $\mathsf{Mod}(\mathbf{R})_{\mathrm{addrel}}$ of (left) $\mathbf{R}$-modules and additive relations. In this context, every additive relation $X$ is intrinsically difunctional ($X \circ X^{-1} \circ X = X$) due to its submodule structure. This ensures that the Riguet congruence $\Phi$ identifies additive relations that share the same homological essence---domain of definition, kernel, image, and indeterminacy---as characterized by Mac Lane~\cite{sM61,sM63}.

\end{example}

\section*{Acknowledgements}
The second and third authors were supported by grant PID2024-159495NB-I00 from the Ministerio de Ciencia e Innovaci\'{o}n, Spain, and the CIAICO/2023/007 grant  from the Conselleria d'Educaci\'{o}, Universitats i Ocupaci\'{o}, Generalitat Valenciana. The second author has held a Specially Appointed Professor position at Nantong University during the completion of this work. The third author was supported by the  CIACIF/2022/489 grant from the Conselleria d'Educaci\'{o}, Universitats i Ocupaci\'{o}, Generalitat Valenciana co-funded by the European Social Fund.

\end{document}